\documentclass[11pt]{amsart}
\usepackage{blindtext}
\usepackage{url}
\usepackage{tabularx}
\usepackage[utf8]{inputenc}
\usepackage[margin=2.5cm]{geometry}
\usepackage[utf8]{inputenc}
\usepackage{amsmath}
\usepackage{amsfonts}
\usepackage{color}
\usepackage{slashbox}
\usepackage{amssymb}
\usepackage{amsthm}
\usepackage{array}
\usepackage{rotating}
\usepackage{physics}
\usepackage{tikz}
\usepackage{braids}
\usepackage{tikz-cd}
\tikzset{%
	symbol/.style={%
		draw=none,
		every to/.append style={%
			edge node={node [sloped, allow upside down, auto=false]{$#1$}}}
	}
}
\usepackage[all]{xy}
\usepackage{cite}
\usepackage{mathdots}
\usepackage{pstricks}
\usepackage{mathrsfs}
\usepackage{setspace}
\usepackage{hyperref}
\theoremstyle{theorem}
\newtheorem{theorem}{Theorem}

\numberwithin{theorem}{subsection}
\newtheorem{lemma}[theorem]{Lemma}
\newtheorem{proposition}[theorem]{Proposition}
\newtheorem{corollary}[theorem]{Corollary}
\theoremstyle{definition}
\newtheorem{definition}[theorem]{Definition}
\newtheorem{example}[theorem]{Example}
\theoremstyle{remark}
\newtheorem{remark}[theorem]{Remark}

\theoremstyle{notation}
\newtheorem{notation}[theorem]{Notation}

\newcommand{\Gr}{\mathbf{Gray}}
\newcommand{\A}{\mathcal{A}}

\newcommand{\Cat}{\mathbf{Cat}}

\title{Tricategorical universal properties via enriched homotopy theory}
\author{Adrian Miranda}
\thanks{The material is based on research conducted while supported by MQRES PhD Scholarship 20192497, and written while supported by EPSRC under grant EP/V002325/2. I thank Steve Lack for his guidance while conducting this research and Nicola Gambino for helpful discussions while I was preparing this paper. I also thank Bojana Femic for useful discussions which led to the inclusion of Remark \ref{Remark lax double functors}.}
\address{Department of Mathematics, University of Manchester, United Kingdom}
\email{adrian.miranda@manchester.ac.uk}
\begin{document}

\maketitle

\begin{abstract}
	We develop the theory of tricategorical limits and colimits, and show that they can be modelled up to biequivalence via certain homotopically well-behaved limits and colimits enriched over the monoidal model category $\mathbf{Gray}$ of $2$-categories and $2$-functors. This categorifies the relationship that bicategorical limits and colimits have with the so called `flexible' enriched limits in $2$-category theory. As examples, we establish the tricategorical universal properties of Kleisli constructions for pseudomonads, Eilenberg-Moore and Kleisli constructions for (op)monoidal pseudomonads, centre constructions for $\mathbf{Gray}$-monoids, and strictifications of bicategories and pseudo-double categories. 
\end{abstract}

\tableofcontents

\section{Introduction}

\subsection{Context and motivation}

 Many interesting constructions on categories enjoy universal properties involving not just functors but also natural transformations, making them two-dimensional in nature. Various constructions on two-dimensional categorical structures have also been considered in the literature. These are motivated by applications to algebraic geometry \cite{Campbell PhD} and quantum field theory \cite{BDSV Extended 3D bordism theory of modular objects, BDSV Modular categories as representations of the 3-dimensional bordism 2-category}, and an analysis of their universal properties is desirable. While enriched universal properties in the sense of \cite{Kelly Basic Concepts of Enriched Category Theory} successfully capture some such phenomena \cite{Coherent Approach to Pseudomonads}, it has been observed that they do not capture other constructions of interest such as Kleisli constructions for pseudomonads \cite{Cheng Hyland Power Pseudodistributive Laws, Formal Theory of Pseudomonads, Miranda Enriched Kleisli objects for pseudomonads}. The inadequacy of strict notions is unsurprising, since when we work with low dimensional higher categories, we are used to familiar equations between morphisms of codimension greater than one often only holding up to coherent higher dimensional data. This is just as true for the commutativity conditions involved in the universal properties of limits and colimits as it is for conditions such as the category axioms, functoriality or naturality.
 \\
 \\
 \noindent The theory of bicategorical limits and colimits treats this phenomenon in the two-dimensional setting \cite{Bird Kelly Power Street, Power Coherence for categories with finite bilimits}. Here the weights on a bicategory $\mathcal{A}$ comprise the $2$-category $\mathbf{Bicat}(\mathcal{A}, \Cat)$ of pseudofunctors into $\Cat$, pseudonatural transformations between them, and modifications between them. These are notions enriched over the monoidal bicategory $\Cat$ \cite{Garner Shulman Enriched Categories as a Free Cocompletion}, but they are also closely related to structures enriched over the one-dimensional structure of categories and functors, with the higher dimensional structure encoded via the canonical monoidal model structure on $\Cat$ whose weak equivalences are the equivalences of categories and whose fibrations are those functors which lift isomorphisms; the \emph{isofibrations} \cite{Gambino Homotopy Limits for 2-categories}. Specifically, bicategorical limits weighted by a pseudofunctor $W: \A \to \Cat$ correspond to `up to equivalence' versions of limits weighted by some $W': \A' \to \Cat$, where there is a biequivalence $E: \A' \sim \A$ and a pseudonatural equivalence $W' \simeq WE$ with $\A'$ a $2$-category and $W'$ a $2$-functor that satisfies the following equivalent conditions.
 
 \begin{enumerate}
 	\item $W'$ is cofibrant in the projective model structure in the enriched functor category $[\A', \Cat]$ \cite{Lack Homotopy Theoretic Aspects of 2-monads}.
 	\item $W'$ is a flexible algebra \cite{Bird Kelly Power Street, Two Dimensional Monad Theory} for the $2$-monad whose forgetful $2$-functor is given by restriction along the inclusion $\mathbf{Ob}(\A') \to \A'$, and whose free $2$-functor is given by left Kan extension.
 \end{enumerate}

\noindent These $\Cat$-enriched weights are also precisely the flexible ones in the sense of two-dimensional monad theory, and moreover any bicategory with finite bilimits is biequivalent to a $2$-category with finite flexible limits. Computing projective cofibrantly weighted, i.e. flexible, limits (resp. colimits) only requires imposing equations (resp. relations) between morphisms in categories, rather than between objects. In this way, bicategorical limits and colimits are captured by up-to-equivalence versions of these simpler, enriched notions. For the limit case these are generated by products, inserters, equifiers, and pseudo splittings of idempotents. Other prominent examples include descent objects of truncated cosimplicial objects, Eilenberg-Moore objects of monads, comma objects, inverters, and pseudo, lax and oplax limits. Non-examples include equalisers, pullbacks and identifiers.
 \\
 \\
 The definition of tricategorical limits and colimits first appeared in the literature as Definition 7.3 of \cite{Three dimensional monad theory}, and their theory is developed further in Section 3.3 of \cite{Campbell PhD}. We continue this research program, characterising the enriched notions which correspond to tricategorical universal constructions. In particular, we work over the category $\mathbf{Gray}$ of $2$-categories and $2$-functors equipped with the $\mathbf{Gray}$-tensor product, for which the corresponding closed structure consists of $2$-functors, pseudonatural transformations and modifications, and the Lack model structure, for which the weak equivalences are the biequivalences and the fibrations are those $2$-functors which both lift adjoint equivalences and isomorphisms in hom-categories; the \emph{equiv-fibrations}. The goal of this paper is to tricategorical universal properties in these terms, via homotopically well-behaved enriched limits and colimits for which the universal property holds up to weak equivalence, i.e. biequivalence, rather than isomorphism.
 \\
 \\
 \noindent Compared to the two-dimensional setting \cite{Gambino Homotopy Limits for 2-categories, Power Coherence for categories with finite bilimits}, there are extra subtleties with the relationship between the following aspects of weights 
 
 \begin{itemize}
 	\item their flexibility when viewing them as algebras for the appropriate enriched monad,
 	\item their good homotopical behaviour, and
 	\item their ability to model tricategorical universal properties.
 \end{itemize}

\noindent Part of this is due to the fact that in dimension three the two kinds of weights are received by different bases $\mathbf{Bicat}$ and $\mathbf{Gray}$, while in the two -dimensional setting they are both received by $\mathbf{Cat}$. The semi-strictification triadjunction of \cite{Campbell Strictification} is used to address this issue. A more substantial cause of subtlety is that not every $2$-category is cofibrant in the Lack model structure of \cite{Quillen 2-cat}. For this reason, our results are are not the most straightforward categorifications of their two-dimensional counterparts. Rather extra clauses to do with pointwise cofibrancy of weights or hom-wise cofibrancy of $\mathbf{Gray}$-categories are required. Such clauses are exactly what is needed to ensure that equations are only imposed on cells of the highest dimension available, since non-cofibrant $2$-categories appearing in the homs of diagram $\Gr$-categories or in the image of weights would introduce such relations. The appropriate homotopical properties are listed below. The second and third of these are new to the tricategorical setting since unlike categories, not all $2$-categories are cofibrant. Note though that $(1)$ and $(2)$ together imply $(3)$.

\begin{enumerate}
	\item projective cofibrancy of the weight, just as in the bicategorical context,
	\item hom-wise cofibrancy of the shape $\mathfrak{A}$ on which the diagram $F: \mathfrak{A} \rightarrow \mathfrak{B}$ and weight $W: \mathfrak{A} \to \mathbf{Gray}$ are indexed,
	\item pointwise cofibrancy of the weight.
\end{enumerate} 

\noindent Theorem \ref{tricolimit reduction to strict 3 k transfors via LW}, Theorem \ref{trinatural transformation classifier cofibrancy} and Corollary \ref{Reduction to Gray natural biequivalence} give more detail on the precise sense in which tricategorical notions are modelled by homotopical $\mathbf{Gray}$-enriched notions. This information is synthesised in Theorem \ref{coherence for tricategories with trilimits}, which says that any tricategory with finite trilimits $\mathfrak{B}$ is triequivalent to a $\mathbf{Gray}$-category $\overline{\mathfrak{B}}$ admitting finite, flexibly weighted $\mathbf{Gray}$-enriched limits on hom-wise cofibrant shapes. This categorifies Power's coherence theorem for bicategories with finite bilimits \cite{Power Coherence for categories with finite bilimits}. We use semi-strictification results on $\left(3, k\right)$-transfors from \cite{Miranda strictifying operational coherences}, as well as the $\mathbf{Gray}$-monad perspective developed in \cite{Buhne PhD} and \cite{Buhne Gray Homo}. Our proofs are similar to their two-dimensional analogues except for the extra need for hom-wise cofibrancy. We also exhibit tricategorical completeness and cocompleteness of examples such as $\mathbf{Bicat}$, and describe concrete examples of tricategorical limits and colimits that are of interest in the theory of pseudomonads, two-dimensional monoidal structures, and strictification.

\subsection{Main results}

\begin{enumerate}
	\item Section \ref{Section reduction to strict pointwise cofibrant weights on cofibrant Gray categories} reduces the tricategories involved in the shape of a tricategorical (co)limit to cofibrant $\mathbf{Gray}$-categories, and the trihomomorphisms involved in their weights and diagrams to $\mathbf{Gray}$-functors. Further assumptions on the weight, such as pointwise cofibrancy, are also shown to not lose generality. These results are collected in Proposition \ref{special weights and diagrams for tricolimits} and summarised in Remark \ref{summarising reduction to pointwise cofibrant weights on cofibrant Gray categories}. The work in Section \ref{Section reduction to strict pointwise cofibrant weights on cofibrant Gray categories} still leaves a universal property in terms of trinatural transformations and trimodifications; data that are weaker than in the $\mathbf{Gray}$-enriched setting.
	\item Section \ref{Weak 3 k transfors as higher cells of pseudoalgebras} simplifies tricategorical universal properties to ones expressible in terms of enriched notions over $\mathbf{Gray}$ as a monoidal model category. Such universal properties involve cells in enriched functor categories $[\mathfrak{A}^\text{op}, \mathbf{Gray}]$ rather than trinatural transformations or trimodifications. We leverage the theory of monads in the $2$-category $\mathbf{Gray}$-$\mathbf{CAT}_{2}:= \mathcal{V}\text{-}\Cat$ with $\mathcal{V} = \Gr$, and the extension of this theory to pseudoalgebras and weak higher cells developed in part 3 of \cite{Gurski Coherence in Three Dimensional Category Theory}, Chapter 1 of \cite{Buhne PhD} and \cite{Buhne Gray Homo}, contributing to this theory by establishing a triadjunction between the base of a $\Gr$-monad and its $\Gr$-category of pseudoalgebras in Proposition \ref{pseudoalgebra triadjunction}. By these means we achieve our goal of capturing tricategorical limits and colimits entirely in terms of enrichment over the monoidal model category $\mathbf{Gray}$. The first main result is Theorem \ref{tricolimit reduction to strict 3 k transfors via LW}, which reduces the universal property of $\mathfrak{B}\left(W\odot F, Y\right)$ to one in terms of a biequivalence in the enriched functor $\mathbf{Gray}$-category $[\mathfrak{A}^\text{op}, \mathbf{Gray}]$ rather than in the tricategory $\mathbf{TRICAT}\left(\mathfrak{A}^\text{op}, \mathbf{Bicat}\right)$. However, this biequivalence still varies only trinaturally in $Y \in \mathfrak{B}$. This last aspect of weakness outside the scope of the monoidal model category $\mathbf{Gray}$ will be removed in Corollary \ref{Reduction to Gray natural biequivalence} once projective cofibrant weights are related to flexible enriched weights on hom-wise cofibrant $\mathbf{Gray}$-categories in Theorem \ref{trinatural transformation classifier cofibrancy}.
	\item In Section \ref{Examples trilimits and tricolimits} we give examples of tricategorical limits and colimits, and show in Proposition \ref{main examples have all trilimits and tricolimits} that several tricategories of interest have tricategorical limits and colimits. We give an explicit formula for trilimits in $\mathbf{Bicat}$ in Proposition \ref{trilimits in Hom}, and use this to analyse the interaction of representables and the Yoneda embedding with trilimits, culminating in Theorem \ref{coherence for tricategories with trilimits} which says that tricategories with finite trilimits are triequivalent to $\mathbf{Gray}$-categories with finite flexible $\mathbf{Gray}$-enriched limits.
	
	\item Our leading examples of tricategorical limits include the following.

	\begin{itemize}
		\item We view a variant of the centre construction for $\mathbf{Gray}$-monoids considered in \cite{Baez Neuchl Braided Monoidal 2-categories} and \cite{Crans Generalised Centers} as a tricategorical descent object in Example \ref{center as descent object}.
		\item We capture the lifted monoidal structures of Eilenberg-Moore pseudoalgebras studied in \cite{Miranda Opmonoidal Pseudomonads} as certain tricategorical Eilenberg-Moore objects.
	\end{itemize}
	
	\item Of tricolimits, our leading examples include the following.
	
	\begin{itemize}
		\item Example \ref{Example strictification of pseudo-double categories} exhibits strictification of pseudo-double categories as a codescent object.
		\item Strictification of bicategories as a tricopower of $\mathbf{1}$, which is considered in Example \ref{strictification as a tricolimit}.
		\item Trikleisli pseudoadjunctions for pseudomonads, which are defined in Definition \ref{trikleisli object definition} and characterised in the tricategories $\mathbf{Bicat}$ and $\mathbf{Gray}$ in Theorem \ref{biessential surjectivity on objects characterises trikleisli pseudoadjunctions in Gray or Hom} as those pseudoadjunctions which have biessentially surjective on objects left pseudoadjoints. This strengthens previous results in \cite{Miranda Enriched Kleisli objects for pseudomonads}, resolving the problems posed in the closing remarks of \cite{Formal Theory of Pseudomonads} and in Section 4 of \cite{Cheng Hyland Power Pseudodistributive Laws} regarding a full tricategorical treatment of the universal property of Kleisli bicategories. 
	\end{itemize}
\end{enumerate} 

\begin{remark}
	By \emph{tricategory}, we will always mean one that is algebraic in the sense of \cite{Gurski PhD}. Such tricategories are more structured that those originally defined in \cite{GPS tricategory}. In particular, they are equipped with specified pseudonatural adjoint equivalences to each of the constraints $\alpha_{h, g, f}: (h.g).f \to h.(g.f)$, $\lambda_{g}: g \to 1_{Y}.g$ and $\rho_{g}: g.1_{X} \to g$. By a \emph{triadjunction}, we will mean a pair of trihomomorphisms $R: \mathfrak{A} \to \mathfrak{B}$, $L: \mathfrak{B} \to \mathfrak{A}$, and a trinatural family of biequivalences $\mathfrak{A}(LX, Y) \sim \mathfrak{B}(X, RY)$ indexed by $(X, Y) \in \mathfrak{A}\times \mathfrak{B}$ (cf. page 35 of \cite{Campbell PhD}).
\end{remark}

\begin{remark}
	Recall \cite{Gurski Coherence in Three Dimensional Category Theory} that if $\mathfrak{A}$ is a tricategory and $\mathfrak{B}$ is either a $\Gr$-category or $\mathbf{Bicat}$ then there is a $\Gr$-category (resp. tricategory) $\mathbf{Tricat}(\mathfrak{A}, \mathfrak{B})$ of trihomomorphisms, trinatural transformations, trimodifications and perturbations. Following \cite{crans tensor of gray categories}, we may refer to the $k$-cells in such a tricategory as \emph{$(3, k)$-transfors}. On the other hand, data appearing in enriched functor $\Gr$-categories $[\mathfrak{A}, \mathfrak{B}]$ when $\mathfrak{A}$ and $\mathfrak{B}$ are $\Gr$-categories will be called \emph{$\Gr$-functors}, \emph{$\Gr$-natural transformations}, and \emph{$\Gr$-modifications}. In both settings the $k$-cells are the perturbations. 
\end{remark}

\section{Background}

\noindent To achieve these results in this paper we must overcome the extra conceptual difficulties stemming from non-cofibrant $2$-categories and the extra computational difficulties associated with weak three-dimensional category theory. In this Section we review some background results which help us to overcome these difficulties, and in particular allow us to avoid many cumbersome calculations in three-dimensional categorical structures.

\subsection{Yoneda for tricategories}

\begin{lemma}\label{Yoneda lemma for tricategories}
	(Theorem 2.12 of \cite{Buhne PhD}) Let $\mathfrak{A}$ be a tricategory, $W: \mathfrak{A}^\text{op} \rightarrow \mathbf{Bicat}$ a trihomomorphism and $X \in \mathfrak{A}$ an object.
	
	\begin{enumerate}
		\item There is a strict homomorphism of bicategories $\varepsilon_{X, V}$ as depicted below, which is moreover a biequivalence.
		
		$$\begin{tikzcd}[font=\fontsize{9}{6}]
			\mathbf{Tricat}\left(\mathfrak{A}^\text{op}{,} \mathbf{Bicat}\right)\left(\mathfrak{A}\left(-{,} X\right){,} V\right) 
			\arrow[rrr, "\varepsilon_{X{,}V}"]
			&&& V\left(X\right)
			\\
			p
			\arrow[dd, rightarrow, bend right, "\sigma"'name=A]
			\arrow[dd, rightarrow, bend left, "\tau"name=B]
			\arrow[from = A,to=B, Rightarrow, "\Phi"]
			&&&
			p_{X}\left(1_{X}\right)
			\arrow[dd, rightarrow, bend right = 60, "{\left(\sigma_{X}\right)}_{1_{X}}"'name=C]
			\arrow[dd, rightarrow, bend left = 60, "{\left(\tau_{X}\right)}_{1_{X}}"name=D]
			\arrow[from = C,to=D, Rightarrow, "{\left(\Phi_{X}\right)}_{1_{X}}"]
			\\
			&\mapsto
			\\
			q &&&q\left(1_{X}\right)
		\end{tikzcd}$$
		\item There is a tri-fully faithful trihomomorphism $\mathfrak{Y}: \mathfrak{A} \rightarrow \mathbf{Tricat}\left(\mathfrak{A}^{\text{op}}, \mathbf{Bicat}\right)$ given by sending $X$ to the representable trihomomorphism $\mathfrak{A}\left(-, X\right): \mathfrak{A}^{\text{op}} \rightarrow \mathbf{Bicat}$.
		\item The biequivalence $\varepsilon_{X, V}$ is the component at $?=X$ of a trinatural biequivalence from the trihomomorphism depicted below, to $V$.
		
		\begin{tikzcd}
			\mathfrak{A}^\text{op}
			\arrow[rr, "\mathfrak{Y}^\text{op}"]
			&& {\mathbf{Tricat}\left(\mathfrak{A}^{\text{op}}{,}\mathbf{Bicat}\right)}^{\text{op}}
			\arrow[rrrr, "\mathbf{Tricat}\left(\mathfrak{A}^{\text{op}}{,}\mathbf{Bicat}\right)\left(-{,}W\right)"]
			&&&& \mathbf{Bicat}
		\end{tikzcd}
		
		\item The assignment $X \mapsto \varepsilon_{X, V}$ defines a strict trinatural transformation.
	\end{enumerate}
\end{lemma}

\begin{proof}
	Parts (1) and (3) are Theorem 2.12 of \cite{Buhne PhD}. Part (2) follows from part (1) by taking $W:= \mathfrak{A}\left(-, X\right)$. Part (4) is easy to see by inspecting the definition of the tricategory $\mathbf{Tricat}\left(\mathfrak{A}^\text{op}, \mathbf{Bicat}\right)$. Let $q: V \rightarrow W$ be another trinatural transformation. Then the component pseudofunctor ${\left(qp\right)}_{X}: \mathfrak{A}\left(X, X\right) \rightarrow W\left(X\right)$ is by definition the same as the composite pseudofunctor $q_{X}p_{X}$, so both of these will have the same output on the object $1_{X} \in \mathfrak{A}\left(X, X\right)$. Similarly if $\phi: q \rightarrow q'$ is a trimodification then the component pseudonatural transformation ${\left(\sigma.p\right)}_{X}$ is by definition equal to the whiskering of $\sigma_{X}$ with $p_{X}$, and hence both will have the same component at $1_{X}$. The argument for perturbations is similar. Finally, the strictness aspect of part (4) follows from the fact that $\mathbf{Bicat}$ has $2$-functorial restrictions.
\end{proof}

\subsection{Strictification and cofibrancy}

\begin{proposition}\label{Hom triequivalent to Gray c}
	Let $\mathbf{Bicat}$ be the tricategory of bicategories, pseudofunctors, pseudonatural transformations and modifications, and let $\Gr$ be its locally-full sub-$\Gr$-category determined by $2$-categories and $2$-functors. Consider the inclusion $I_{2}: \Gr \to \mathbf{Bicat}$.
	
	\begin{enumerate}
		\item $I_{2}$ is triessentially surjective on objects.
		\item $I_{2}$ has a tri-fully faithful left triadjoint $\mathbf{st}_{2}$.
		\item The triessential image $\Gr_\text{c}$ of $\mathbf{st}_{2}$ consists of cofibrant $2$-categories, or those $2$-categories whose underlying categories are free on graphs.
	\end{enumerate} 
\end{proposition}

\begin{proof}
	See \cite{Bicat Not Triequivalent To Gray} for part (1), Section 3.1 of \cite{Campbell PhD} for part (2), and Theorem 4.8 of \cite{Quillen 2-cat} for part (3).
\end{proof}

\begin{proposition}\label{Proposition 3D strictification tetraadjunction}
	Let $\mathfrak{A}$ be a tricategory and $\mathfrak{B}$ be a $\Gr$-category.
	
	\begin{enumerate}
		\item There is a cofibrant $\Gr$-category $\mathbf{st}_{3}(\mathfrak{A})$ and a triequivalence $E_{\mathfrak{A}}: \mathfrak{A} \rightarrow \mathbf{st}_{3}(\mathfrak{A})$, which is moreover a biequivalence in the tricategory of small tricategories defined in \cite{Low dimensional structures formed by tricategories}, and a unit for the triadjunction $\mathbf{st}_{3} \dashv I_{3}$ of Proposition 4.3.4 part (4) in \cite{Miranda strictifying operational coherences}.
		\item If $\mathfrak{A}$ is a $\Gr$-category whose underlying sesquicategory is free on a $2$-computad then \begin{enumerate}
			\item  $E_\mathfrak{A}$ has a section given by a $\Gr$-functor, which is also a triequivalence.
			\item Any trihomomorphism $\mathfrak{A} \to \mathfrak{B}$ is trinaturally equivalent to a $\Gr$-functor.
			\item Any trinatural transformation $p: F \Rightarrow G$ between $\Gr$-functors from $\mathfrak{A}$ to $\mathfrak{B}$ admits an equivalence trimodification $p \cong p'$ such that $p'$ is semi-strict in the sense of Definition 3.2.3 of \cite{Miranda strictifying operational coherences}.
		\end{enumerate}
		\item There is a $\Gr$-category $\mathbf{Tricat}(\mathfrak{A}, \mathfrak{B})$ of trihomomorphisms, trinatural transformations, trimodifications and perturbations.
		\item If $\mathfrak{A}$ is a $\Gr$-category then there is a sub-$\Gr$-category $I_{\mathfrak{A}, \mathfrak{B}}: \mathbf{Gray}$-$\mathbf{Cat}(\mathfrak{A}, \mathfrak{B}) \rightarrow \mathbf{Tricat}(\mathfrak{A}, \mathfrak{B})$, consisting of $\Gr$-functors and semi-strictly decomposable trinatural transformations (Definition 5.3.1 of \cite{Miranda strictifying operational coherences}).
		\item The $\Gr$-functor displayed below is a triequivalence.
	\end{enumerate}

$$\begin{tikzcd}
	\mathbf{\Gr}\text{-}\Cat(\mathbf{st}_{3}(\mathfrak{A}){,}\mathfrak{B}) \arrow[rrr, "	\mathbf{\Gr}\text{-}\Cat(E_{\mathfrak{A}}{,}\mathfrak{B})"] &&& 	\mathbf{\Gr}\text{-}\Cat(\mathfrak{A}{,}\mathfrak{B}) \arrow[r, "I_{\mathfrak{A}{,}\mathfrak{B}}"] &\mathbf{Tricat}(\mathfrak{A}{,}\mathfrak{B})
\end{tikzcd}$$
\end{proposition}

\begin{proof}
	See Theorem 10.11 of \cite{Gurski Coherence in Three Dimensional Category Theory} and \cite{Gray categories model algebraic tricategories} for part (1), Section 9 of \cite{Quillen Gray Cat} for part (2a), Corollary 4.1 of \cite{Gray categories model algebraic tricategories} for part (2b), and Proposition 4.3.2 of \cite{Miranda strictifying operational coherences} for part (2c). Part (3) is Theorem 9.4 of \cite{Gurski Coherence in Three Dimensional Category Theory} and part (4) is an easy consequence. Finally, part (5) is Theorem 5.3.4 of \cite{Miranda strictifying operational coherences}.
\end{proof}

\noindent Although explicit descriptions of $\mathbf{st}_{2}$ and $\mathbf{st}_{3}$ can be inferred from their universal properties as weak functor classifiers, as in \cite{Miranda strictifying operational coherences}, such descriptions will not be needed in this paper. We will however give an explicit description of the generalisation of $\mathbf{st}_{2}$ to the pseudo-double category setting in Example \ref{Example strictification of pseudo-double categories}.

\begin{proposition}\label{Proposition cofibrancy localises}
	Let $\mathfrak{A}$ be a cofibrant $\mathbf{Gray}$-category and consider a pair of objects $X, Y \in \mathfrak{A}$. Then the hom-$2$-category $\mathfrak{A}(X, Y)$ is cofibrant. 
\end{proposition}

\begin{proof}
	Suppose that $\mathfrak{A}$ is a $\mathbf{Gray}$-category whose underlying sesquicategory is free on a $2$-computad $\mathbb{A}$. We need to show that hom-$2$-categories of $\mathfrak{A}$ also have underlying categories that are free on graphs. Indeed, they admit presentations in terms of $2$-computads $\overline{\mathbb{A}}\left(X, Y\right)$, with
	\begin{itemize}
		\item vertices given by morphisms in $\mathfrak{A}$, hence paths in  the free sesquicategory on $\mathbb{A}$,
		\item edges given by generating $2$-cells in $\mathbb{A}$ that may have been whiskered on either side,
		\item Two different kinds of $2$-cells. \begin{itemize}
			\item One kind given by $3$-cells in $\mathbb{A}$ that may have been whiskered on either side by $1$-cells in $\mathfrak{A}$.
			\item The other kind given by $2$-cells $\left(y, \psi, w, \phi, x\right)$ and ${\left(y, \psi, w, \phi, x\right)}^{-1}$, as depicted below right, given data in $\mathfrak{A}\left(X, Y\right)\left(yhwfx, ykwgx\right)$ as depicted below left. Here $\phi$ and $\psi$ are generating $2$-cells in $\mathbb{A}$, but all $1$-cells may be paths.
		\end{itemize} 
	\end{itemize} 
	$$\begin{tikzcd}[column sep = 18, font=\fontsize{9}{6}]
		X \arrow[r, "x"]
		& X' \arrow[rr, bend right, "g"' name = g]\arrow[rr, bend left, "f" name = f]
		&&W \arrow[r, "w"]
		& W'\arrow[rr, bend left, "h" name = h]\arrow[rr, bend right, "k"' name = k]
		&& Y'\arrow[r, "y"]
		&Y&{}
		\arrow[from = f, to = g, Rightarrow, shorten = 5, "\phi"]
		\arrow[from = h, to = k, Rightarrow, shorten = 5, "\psi"]
	\end{tikzcd}\begin{tikzcd}[font=\fontsize{9}{6}]
		yhwfx \arrow[rr, "y\psi wfx"] \arrow[dd, "yhw\phi x"']
		&{}\arrow[dd, Rightarrow,shorten = 8, shift right = 5]
		\arrow[dd, Leftarrow, shift left=5, shorten = 8]
		& ykwfx
		\arrow[dd, "ykw\phi x"]
		\\
		\\
		yhwgx\arrow[rr, "y\psi wgx"'] &{}
		& ykwgx
	\end{tikzcd}$$
	\begin{itemize}
		\item For every relation on $3$-cells between objects $X$ and $Y$ in the presentation of $\mathfrak{A}$, and pair of morphisms $W \rightarrow X$ and $Y \rightarrow Z$ in $\mathfrak{A}$, there is a `whiskered relation' on corresponding $2$-cells in $\mathfrak{A}\left(X, Y\right)$.
		\item There are extra relations giving the $\mathbf{Gray}$-category axioms that are not also axioms of their hom-$2$-categories. These correspond to relations $\left([l],\Phi, [r]\right)\sim \left([l], \Psi, [r]\right)$ whenever $\Phi = \Psi$ in the free $\mathbf{Gray}$-category. 
	\end{itemize}
	
	\noindent One observes that the only relations above are on $3$-cells. The correctness of this description may be observed by appealing to a combinatorial description of free $\Gr$-categories on $\Gr$-computads, such as what is given in Section 2 of \cite{Morehouse 2-categories from a Gray Perspective} or Remark 3.2.15 of \cite{Miranda PhD}. Alternatively, one can explicitly construct the $\Gr$-category structure on this $2$-$\Cat$ enriched graph and then check the required universal property as is done for a related structure in Section 3.2 of \cite{Miranda PhD}. We leave these details to the interested reader.
\end{proof}

\subsection{Other key tools}

\noindent We recall that trinatural biequivalences can be detected pointwise.

\begin{lemma}\label{pointwise biequivalence is biequivalence}
	Let $\mathfrak{B}$ be either $\mathbf{Bicat}$ or a cocomplete $\mathbf{Gray}$-category and let $p: F \rightarrow G: \mathfrak{A} \rightarrow \mathfrak{B}$ be a trinatural transformation between trihomomorphisms. If all components of $p$ are part of biequivalences internal to $\mathfrak{B}$, then $p$ is itself part of a biequivalence internal to $\mathbf{TRICAT}\big(\mathfrak{A}, \mathfrak{B}\big)$.
\end{lemma}

\begin{proof}
	See \cite{Buhne PhD} Proposition A.16.
\end{proof}

\noindent The following says that left triadjoints preserve tricolimits. It appears as Lemma 3.3.2 of \cite{Campbell PhD}.

\begin{lemma}\label{left triadjoint preserves tricolimits}
	Let $L: \mathfrak{B} \rightsquigarrow \mathfrak{C}$ be a trihomomorphism between tricategories and suppose $R: \mathfrak{C} \rightsquigarrow \mathfrak{B}$ is right triadjoint to $L$. If $W \odot F$ is a tricolimit of $F:\mathfrak{A}\rightsquigarrow \mathfrak{B}$ weighted by $W: \mathfrak{A}^\text{op} \rightsquigarrow \mathbf{Bicat}$ then $L\left(W\odot F\right)$ is a tricolimit of $LF$ weighted by $W$.
\end{lemma}

\begin{proof}
	Observe that there is the following chain of biequivalences, each of which is trinatural in $X$. The first is the biequivalence given as part of the triadjunction $L \dashv R$, the second is the definition of $W \odot F$, and the third is post-composition along a certain biequivalence. Specifically, the restriction along $F$ of the trinatural biequivalence $\mathfrak{B}\left(-. RX\right) \sim \mathfrak{C}\left(L-, X\right)$ in $\mathbf{TRICAT}\big(\mathfrak{B}^\text{op}, \mathbf{Bicat}\big)$ given by the hom-wise triadjunction $L \dashv R$.
	
	\begin{align*}
		\mathfrak{C}\left(L\left(W\odot F\right), X\right) &\sim \mathfrak{B}\left(W \odot F, RX\right) \\
		&\sim \mathbf{TRICAT}\big(\mathfrak{A}^\text{op}, \mathbf{Bicat}\big)\left(W, \mathfrak{B}\left(F-, RX\right)\right) \\
		&\sim  \mathbf{TRICAT}\big(\mathfrak{A}^\text{op}, \mathbf{Bicat}\big)\left(W, \mathfrak{C}\left(LF-, X\right)\right)
	\end{align*}
	
\end{proof}

\section{Definitions and simplifications of the data}\label{Section reduction to strict pointwise cofibrant weights on cofibrant Gray categories}

\noindent This Section will make heavy use of the Yoneda Lemma for tricategories (Lemma \ref{Yoneda lemma for tricategories}), the strictification triadjunction recalled in Proposition \ref{Hom triequivalent to Gray c}, and the pointwise nature of biequivalences in functor tricategories recalled in Lemma \ref{pointwise biequivalence is biequivalence}).
\\
\\
 \noindent In Subsection \ref{Subsection definition of trilimit and immediate corollaries} we restrict our attention to corollaries of the definition which follow without using sophisticated coherence results. These results are summarised below. They are all easy three-dimensional versions of basic results in category theory.

\begin{itemize}
	\item Corollary \ref{trilimits stable under biequivalence} is stability under biequivalence.
	\item Corollary \ref{tri(co)limit weighted by representable is evaluation} shows that tricategorical limits and colimits weighted by representables can be given by evaluating the diagram at the representing object.
	\item Corollary \ref{tri(co)limit weighted by representable is evaluation} shows that any weight can be seen as a self-weighted tricolimit of a diagram of representables.
	\item Proposition \ref{taking tricolimits is cocontinuous in the weight} records that taking tricolimits is tricocontinuous in the weight, and Corollary \ref{triessential surjectivity on objects aspect of universal property of tricocompletion} uses this to establish the triessential surjectivity on objects fragment of the expected universal property of free tricocompletions.
	\item Proposition \ref{tricolimits preserves by I: Gray --> Hom} gives sufficient conditions on the weight for a tricolimit for it to be preserved by $I_{2}: \mathbf{Gray} \rightarrow \mathbf{Bicat}$.
\end{itemize}

\noindent Subsection \ref{subsection reduction to simpler shapes targets weights and diagrams} will use the results of \cite{Miranda strictifying operational coherences} to show that certain simplifying assumptions on weights and diagrams involved in tricategorical limits and colimits can be made without loss of generality. It will also make use of Proposition \ref{Hom triequivalent to Gray c} and Proposition \ref{Proposition 3D strictification tetraadjunction}.

\begin{itemize}
	\item Corollary \ref{reduction to cofibrant targets} will mean that the existence of tricategorical limits and colimits can be detected in triequivalent $\mathbf{Gray}$-categories.
	\item Lemma \ref{varying the weight and the diagram} allows weights and diagrams to vary along trinatural biequivalences. The dual statement for tricolimits follows similarly, using $F^\text{op}$ in place of $F$.
	\item Corollary \ref{reduction to pointwise cofibrant weights} allows one to consider weights which factor through the trihomomorphism \begin{tikzcd}
		\mathbf{Bicat} \arrow[r, "\mathbf{st}_{2}"] &\mathbf{Gray} \arrow[r, "I"] &\mathbf{Bicat}
	\end{tikzcd}. Among other things, such weights are pointwise cofibrant.
	\item Corollary \ref{reduction to cofibrant shapes} allows us to restrict attention to shapes which are cofibrant $\mathbf{Gray}$-categories.
	\item Finally, Proposition \ref{special weights and diagrams for tricolimits} combines these results, and also simplifies the trihomomorphisms involved in the weight and the diagram to $\mathbf{Gray}$-functors.
\end{itemize}

\subsection{The definition, and immediate corollaries}\label{Subsection definition of trilimit and immediate corollaries}

\noindent We recall the definition of tricategorical limits and colimits in Definition \ref{trilimits and tricolimits definition}, before which we fix notation for $\mathbf{Gray}$-enriched limits and colimits.

\begin{notation}\label{Notation Gray enriched limits and colimits}
	When $W: \mathfrak{A} \rightarrow \mathbf{Gray}$, $W': \mathfrak{A}^\text{op} \rightarrow \mathbf{Gray}$ and $F: \mathfrak{A} \rightarrow \mathfrak{B}$ are $\mathbf{Gray}$-functors, the $\mathbf{Gray}$-enriched weighted limits and colimits are denoted $\{W, F\}$ and $W'\cdot F$.
\end{notation}

\begin{definition}\label{trilimits and tricolimits definition}
	Let $\mathfrak{A}$, $\mathfrak{B}$ be tricategories and let $W: \mathfrak{A} \rightarrow \mathbf{Bicat}$, $W': \mathfrak{A}^\text{op} \rightarrow \mathbf{Bicat}$, and $F: \mathfrak{A} \rightarrow \mathfrak{B}$ be trihomomorphisms.
	
	\begin{enumerate}
		\item A \emph{trilimit} of $F$ weighted by $W$, if it exists, will be a representing object for the trihomomorphism depicted below. This representing object will be denoted $<W, F>$.
		
		$$\begin{tikzcd}
			\mathfrak{B} 
			\arrow[r, "\mathfrak{Y}"]
			& \mathbf{TRICAT}\left({\mathfrak{B}}^{\text{op}}{,}\mathbf{Bicat}\right)
			\arrow[r, "-\circ F^\text{op}"]
			&  \mathbf{TRICAT}\left({\mathfrak{A}}^{\text{op}}{,}\mathbf{Bicat}\right)
			\arrow[rrrr, "\mathbf{Tricat}\left(\mathfrak{A}^\text{op}{,}\mathbf{Bicat}\right)\left(W{,}-\right)"] 
			&&&& \mathbf{Bicat}
		\end{tikzcd}$$
		
		\item A \emph{tricolimit} of $F$ weighted by $W'$ will be, if it exists, a trilimit $<W', F^\text{op}>$ of $F^\text{op}$ weighted by $W'$. This will be denoted $W' \odot F$.
	\end{enumerate}
	In both cases $\mathfrak{A}$ will be called the \emph{shape}, $\mathfrak{B}$ will be called the \emph{target}, $F$ will be called the \emph{diagram} and $W$ (resp. $W'$) will be called the \emph{weight}.
\end{definition}

\begin{corollary}\label{trilimits stable under biequivalence}
	\hspace{1mm}
	\begin{enumerate}
		\item The property of being a trilimit is stable under extension along maps $w: Y \rightarrow Z$ such that $\mathfrak{B}\left(X, w\right): \mathfrak{B}\left(X, Y\right) \rightarrow \mathfrak{B}\left(X, Z\right)$ is a biequivalence for all $X$.
		\item The property of being a tricolimit is stable under restriction along maps $w: X \rightarrow Y$ such that $\mathfrak{B}\left(w, Z\right): \mathfrak{B}\left(Y, Z\right) \rightarrow \mathfrak{B}\left(X, Z\right)$ is a biequivalence for all $Z$.
		\item The property of being a tricategorical limit of colimit is stable under internal biequivalence.
	\end{enumerate} 
\end{corollary}

\begin{proof}
	For part (1), we need only show trinaturality. This follows from Lemma \ref{Yoneda lemma for tricategories}. Part (2) is dual, while part (3) is a special case since representables preserve internal biequivalences.
\end{proof}

\noindent We now show that when the weight is of the form $\mathfrak{A}\left(-, X\right): \mathfrak{A}^\text{op} \rightarrow \mathbf{Bicat}$, trilimits and tricolimits exist and can be given by evaluation of the diagram on $X$. The generalisation to representable weights, that is ones which are merely biequivalent to such in $\mathbf{Tricat}\left(\mathfrak{A}^\text{op}, \mathbf{Bicat}\right)$, will follow from Lemma \ref{varying the weight and the diagram} part (1).

\begin{corollary}\label{tri(co)limit weighted by representable is evaluation}
	Let $F: \mathfrak{A} \rightarrow \mathfrak{B}$ be a trihomomorphism and let $X$ be an object of $\mathfrak{A}$. \begin{enumerate}
		\item The object $FX$ provides a tricolimit of $F$ weighted by $\mathfrak{A}\left(-, X\right): \mathfrak{A}^\text{op} \rightarrow \mathbf{Bicat}$,
		\item The object $FX$ provides a trilimit of $F$ weighted by $\mathfrak{A}\left(X, -\right): \mathfrak{A} \rightarrow \mathbf{Bicat}$.
	\end{enumerate}
\end{corollary}

\begin{proof}
	For part (1), the defining biequivalence at $Y \in \mathfrak{B}$ is $\varepsilon_{X, V}$ of Lemma \ref{Yoneda lemma for tricategories} part (1), with $V := \mathfrak{B}\left(F-, Y\right)$, and trinaturality in $Y$ is part of the Yoneda embedding for $\mathfrak{B}$. Part (2) is part (1) applied to $F^\text{op}$.
\end{proof}

\noindent Next, we show that every weight $V$ is itself a tricolimit, namely of the Yoneda embedding $\mathfrak{Y}$ weighted by $V$ itself.

\begin{corollary}\label{every presheaf is a tricolimit of yoneda}
	Let $V: \mathfrak{A}^\text{op} \rightarrow \mathbf{Bicat}$ be a trihomomorphism. Then $V$ provides a tricolimit of the diagram $\mathfrak{Y}: \mathfrak{A} \rightarrow \mathbf{TRICAT}\left(\mathfrak{A}^\text{op},\mathbf{Bicat}\right)$ weighted by $V$, so that the defining biequivalences of the universal property  vary strictly naturally in $W \in \mathbf{TRICAT}\left(\mathfrak{A}^\text{op}, \mathbf{Bicat}\right)$.
\end{corollary}

\begin{proof}
	The required biequivalence is displayed below. This is given by extension along $\varepsilon_{X, W}$ of Lemma \ref{Yoneda lemma for tricategories} part (1). Lemma \ref{Yoneda lemma for tricategories} part (3) is the assertion of strict naturality in $W$.
	
	$$\mathbf{TRICAT}\left(\mathfrak{A}^\text{op}, \mathbf{Bicat}\right)\left(V, W\right) \sim \mathbf{TRICAT}\left(\mathfrak{A}^\text{op}, \mathbf{Bicat}\right)\left(V, \mathbf{TRICAT}\left(\mathfrak{A}^\text{op}, \mathbf{Bicat}\right)\left(\mathfrak{Y}-, W\right)\right)$$
\end{proof}

\begin{proposition}\label{taking tricolimits is cocontinuous in the weight}
	Suppose that $\mathfrak{B}$ has small tricolimits $W \odot F$ for some trihomomorphism $F: \mathfrak{A}\rightsquigarrow\mathfrak{B}$. Consider the trihomomorphism depicted below.

	$$\begin{tikzcd}
		\mathfrak{B} \arrow[r, "\mathfrak{Y}"] & \mathbf{TRICAT}\left(\mathfrak{B}^\text{op}{,}\mathbf{Bicat}\right) \arrow[rrr, "\mathbf{TRICAT}\left({F}^\text{op}{,}\mathbf{Bicat}\right)"] &&& \mathbf{TRICAT}\left(\mathfrak{A}^\text{op}{,}\mathbf{Bicat}\right)
	\end{tikzcd}$$
	
	\noindent This trihomomorphism has a left triadjoint given 
	\begin{itemize}
		\item on objects by $W \mapsto W \odot F$,
		\item between hom-bicategories as depicted below. Here $\sim$ denotes the biequivalence defining the respective tricolimits, and $\eta_{W'}$ corresponds to $1_{W'\odot F}$ under $\sim$.

	\end{itemize}
	$$\begin{tikzcd}[font=\fontsize{9}{6}]
		\mathbf{TRICAT}\left(\mathfrak{A}^\text{op}{,}\mathbf{Bicat}\right)\left(W{,}W'\right)
		\arrow[rr, "\eta_{W'}\circ -"] 
		&& \mathbf{TRICAT}\left(\mathfrak{A}^\text{op}{,}\mathbf{Bicat}\right)\left(W{,}\mathfrak{B}\left(F-{,}W'\odot F\right)\right) 
		\arrow[rr, "\sim"]
		&& \mathfrak{B}\left(W\odot F{,}W'\odot F\right)
	\end{tikzcd}$$
\end{proposition}

\begin{proof}
	The remaining data of the trihomomorphism can be constructed from the data of the biequivalence defining the tricolimit. The universal property is immediate from the definition of tricolimit.
\end{proof}

\begin{corollary}\label{triessential surjectivity on objects aspect of universal property of tricocompletion}
	Let $\mathfrak{B}$ be a tricategory with small tricolimits, and let $\mathfrak{A}$ be a small tricategory. For any trihomomorphism $F: \mathfrak{A} \rightsquigarrow \mathfrak{B}$, the left triadjoint $-\odot F: \mathbf{TRICAT}\left(\mathfrak{A}^\text{op}, \mathfrak{B}\right)$ of Proposition \ref{taking tricolimits is cocontinuous in the weight} restricts along $\mathfrak{Y}_{\mathfrak{A}}: \mathfrak{A} \rightsquigarrow \mathbf{TRICAT}\left(\mathfrak{A}^\text{op}, \mathbf{Bicat}\right)$ to a trihomomorphism that is trinaturally biequivalent to $F$.
\end{corollary}

\begin{proof}
	The component on $X \in \mathfrak{A}$ of the trinatural biequivalence is given by the universal property of $\mathfrak{A}\left(-, X\right) \odot F$, as per Corollary \ref{tri(co)limit weighted by representable is evaluation} part (1).
\end{proof}

\noindent Proposition \ref{tricolimits preserves by I: Gray --> Hom}, to follow, will later be strengthened once we show that $\mathbf{Gray}$ and $\mathbf{Bicat}$ both have all tricolimits, and give conditions on $W$ which suffice for $W \odot F$ to also be a cofibrant $2$-category. 

\begin{proposition}\label{tricolimits preserves by I: Gray --> Hom}
	Let $F: \mathfrak{A} \rightsquigarrow \mathbf{Gray}$ be a trihomomorphism such that $FX$ is a cofibrant $2$-category for every $X \in \mathfrak{A}$, and let $W: \mathfrak{A}^\text{op} \rightsquigarrow \mathbf{Bicat}$ be such that the tricolimit $W \odot F$ exists and is also a cofibrant. Then $I$ preserves the tricolimit $W \odot F$.
\end{proposition}

\begin{proof}
	Let $\mathcal{B}$ be a bicategory. The claim follows from the following chain of biequivalences, which are trinatural in $\mathcal{B}$. Note that the first biequivalence uses the assumption that $W \odot F$ is cofibrant while the last biequivalence uses the assumption that $F$ is pointwise cofibrant.
	
	{\small
		\begin{align*}
			\mathbf{Bicat}\left(I\left(W\odot F\right){,}\mathcal{B}\right)
			&\sim \mathbf{Gray}_\text{c}\left(W\odot F{,}\mathbf{st}_{2}\mathcal{B}\right) 
			&\text{Proposition }\ref{Hom triequivalent to Gray c}\text{ part (2)}
			\\
			&= \mathbf{Gray}\left(W\odot F{,}\mathbf{st}_{2}\mathcal{B}\right) 
			&\text{full sub-$\mathbf{Gray}$-category}
			\\
			&\sim \mathbf{TRICAT}\left(\mathfrak{A}^\text{op}{,}\mathbf{Bicat}\right)\left(W{,}\mathbf{Gray}\left(F-{,}\mathbf{st}_{2}\mathcal{B}\right)\right)
			&\text{Definition of} W\odot F
			\\
			&\sim \mathbf{TRICAT}\left(\mathfrak{A}^\text{op}{,}\mathbf{Bicat}\right)\left(W{,}\mathbf{Bicat}\left(IF-{,}\mathcal{B}\right)\right)
			&\text{Proposition } \ref{Hom triequivalent to Gray c}\text{ part (2)}
		\end{align*}
	}
\end{proof}

\subsection{Reductions to simpler shapes, weights and diagrams}\label{subsection reduction to simpler shapes targets weights and diagrams}

\begin{corollary}\label{reduction to cofibrant targets}
	Let $F: \mathfrak{A} \rightarrow \mathfrak{B}$ and $W : \mathfrak{A} \rightarrow \mathbf{Bicat}$ be trihomomorphisms between tricategories. Let $E_{\mathfrak{B}}: \mathfrak{B} \rightsquigarrow \mathbf{st}_{3}\left(\mathfrak{B}\right)$ be the unit from Proposition \ref{Proposition 3D strictification tetraadjunction} part (5). Then the trilimit $<W, F>$ exists if and only if the trilimit $<W, E_\mathfrak{B}F>$ exists, in which case $<W, E_\mathfrak{B}F>\sim E_\mathfrak{B}\left(<W, F>\right)$
\end{corollary}

\begin{proof}
	Observe that $E_\mathfrak{B}$ and its triequivalence inverse are in particular right triadjoints. The result follows from Lemma \ref{left triadjoint preserves tricolimits}.
\end{proof}

\begin{lemma}\label{varying the weight and the diagram}
	Let $W: \mathfrak{A} \rightarrow \mathbf{Bicat}$ and $F: \mathfrak{A} \rightarrow \mathfrak{B}$ be trihomomorphisms. \begin{enumerate}
		\item Let $w: W \rightarrow W'$ be part of a biequivalence in $\mathbf{TRICAT}\big(\mathfrak{A}, \mathbf{Bicat}\big)$. Then $<W, F>$ exists if and only if $<W', F>$ exists, in which case they are biequivalent.
		\item Let $f: F \rightarrow F'$ be part of a biequivalence in $\mathbf{TRICAT}\big(\mathfrak{A}, \mathfrak{B}\big)$. Then $<W, F>$ exists if and only if $<W, F'>$ exists, in which case they are biequivalent.
	\end{enumerate}
\end{lemma}

\begin{proof}
	For part (1), restriction along the biequivalence $w: W \rightarrow W'$ in $\mathbf{Tricat}(\mathfrak{A}, \mathbf{Bicat})$ defines a biequivalence $\mathbf{TRICAT}[\mathfrak{A}, \mathbf{Bicat}]\left(W', \mathfrak{B}\left(F-, X\right)\right) \rightarrow \mathbf{TRICAT}[\mathfrak{A}, \mathbf{Bicat}]\left(W, \mathfrak{B}\left(F-, X\right)\right) $ which is trinatural in $X$. For part (2) we similarly extend along the biequivalence $\mathfrak{B}\left(f, X\right): \mathfrak{B}\left(F'-, X\right) \rightarrow \mathfrak{B}\left(F-, X\right)$.
\end{proof}

\noindent Just as Corollary \ref{reduction to cofibrant targets} allowed us to restrict attention to targets which are cofibrant $\mathbf{Gray}$-categories, Corollaries \ref{reduction to pointwise cofibrant weights} and \ref{reduction to cofibrant shapes} to follow will together allow us to further restrict our attention to pointwise cofibrant weights on cofibrant shapes.

\begin{corollary}\label{reduction to pointwise cofibrant weights}
	\hspace{1mm}
	\begin{enumerate}
		\item A tricategory has the trilimit $<W, F>$ if and only if it has the trilimit $<W', F>$ where $W':= \left(I\mathbf{st}_{2}\right)W$.
		\item A tricategory has all trilimits if it has all trilimits weighted by $W: \mathfrak{A} \rightarrow \mathbf{Bicat}$ such that \begin{itemize}
			\item $WA$ is a cofibrant $2$-category for all $A \in \mathfrak{A}$,
			\item $Wa: WA \rightarrow WA'$ is a $2$-functor whose underlying functor is free on a morphism of graphs, for all $a \in \mathfrak{A}\left(A, A'\right)$, and
			\item $W \phi: W a \Rightarrow Wa'$ is a pseudonatural transformation whose $1$-cell components are generators.
		\end{itemize}
	\end{enumerate}
\end{corollary}

\begin{proof}
	For part (1), apply Lemma \ref{varying the weight and the diagram} part (1) to the biequivalence depicted below, where $s$ is the unit of the semi-strictification triadjunction of Proposition 4.3.4 part (4) of \cite{Miranda Enriched Kleisli objects for pseudomonads} and $\lambda$ is the left unitor in the tricategory $\mathbf{TRICAT}_{3}$ of large tricategories, which differs from \cite{Low dimensional structures formed by tricategories} only in size. Part (2) follows by observing that data in the image of $\mathbf{st}_{2}$ have the properties mentioned.
	
	$$\begin{tikzcd}
		W
		\arrow[r, "\lambda_{W}"]
		& 1_{\mathbf{Bicat}}.W
		\arrow[r, "s.1_{W}"]
		&
		\left(I\mathbf{st}_{2}\right).W
	\end{tikzcd}$$
\end{proof}

\begin{corollary}\label{reduction to cofibrant shapes}
	\hspace{1mm}
	\begin{enumerate}
		\item A tricategory has the trilimit $<W, F>$ if and only if it has the trilimit $<WE_\mathfrak{A}^{*}, FE_\mathfrak{A}^{*}>$, where $E_{\mathfrak{A}}^{*}$ is biequivalence inverse in $\mathbf{Tricat}_{3}$ to the unit of $\mathbf{st}_{3} \dashv I_{3}$.
		\item A tricategory has all trilimits if and only if it has all trilimits weighted by $W: \mathfrak{A} \rightsquigarrow \mathbf{Bicat}$ such that $\mathfrak{A}$ is a cofibrant $\mathbf{Gray}$-category and $W$ factors through $I\mathbf{st}_{2}: \mathbf{Bicat} \rightsquigarrow \mathbf{Bicat}$.
	\end{enumerate}
\end{corollary}

\begin{proof}
	Part (1) follows from two-out-of-three, after observing that restriction along $\mathfrak{E}_{*}: \mathbf{st}_{3}\left(\mathfrak{A}\right)\rightarrow \mathfrak{A}$ defines a trinatural biequivalence depicted below. Part (2) follows by first using Corollary \ref{reduction to pointwise cofibrant weights} and then observing that $\mathbf{st}_{3}\left(\mathfrak{A}\right)$ are cofibrant $\mathbf{Gray}$-categories.
	
	$$\mathbf{TRICAT}[\mathfrak{A}, \mathbf{Bicat}]\left(W, \mathfrak{B}\left(F-, X\right)\right) \rightarrow \mathbf{TRICAT}\left(\mathbf{st}_{3}\left(\mathfrak{A}\right), \mathbf{Bicat}\right)\left(WE_\mathfrak{A}, \mathfrak{B}\left(FE_\mathfrak{A}-, X\right)\right)$$ 
\end{proof}

\noindent By Proposition 7.4 of \cite{Three dimensional monad theory}, a $\mathbf{Gray}$-category has all trilimits if and only if it has trilimits whose diagram and weight are both $\mathbf{Gray}$-functors. This reduction uses two steps. First, it composes a general weight $W: \mathfrak{A} \rightsquigarrow \mathbf{Bicat}$ with the left triadjoint $\mathbf{st}_{2}: \mathbf{Bicat} \rightarrow \mathbf{Gray}$. Second, it uses the triequivalence $ \mathbf{st}_{3}\left(\mathfrak{A}^\text{op}\right)\rightsquigarrow \mathfrak{A}^\text{op}$, and the reflection of weak transfors between tricategories into $\mathbf{Gray}$-enriched analogues. We will also use this reflection in Section \ref{Weak 3 k transfors as higher cells of pseudoalgebras}, but can already strengthen Proposition 7.4 of \cite{Three dimensional monad theory} using Proposition 4.1.1 of \cite{Miranda strictifying operational coherences}.

\begin{proposition}\label{special weights and diagrams for tricolimits}
	\hspace{1mm}
	\begin{enumerate}
		\item A tricategory $\mathfrak{B}$ has the trilimit of a general diagram $F: \mathfrak{A} \rightarrow \mathfrak{B}$ and weight $W: \mathfrak{A} \rightarrow \mathbf{Bicat}$ if and only if for any triequivalence $G: \mathfrak{B} \rightarrow \mathfrak{C}$ with $\mathfrak{C}$ a $\mathbf{Gray}$-category, $\mathfrak{C}$ has the trilimit of the $\mathbf{Gray}$ functor $F'$ equivalent to the composite trihomomorphism depicted below via Proposition 4.1.1 of \cite{Miranda strictifying operational coherences}.
		
		$$\begin{tikzcd}
			\mathbf{st}_{3}\left(\mathfrak{A}\right) \arrow[r, "E_\mathfrak{A}"] & \mathfrak{A} \arrow[r, "F"] & \mathfrak{B} \arrow[r, "G"] & \mathfrak{C}
		\end{tikzcd}$$
		
		\noindent weighted by the $\mathbf{Gray}$-functor $IW'$, where $W'$ is the $\mathbf{Gray}$-functor equivalent to the composite trihomomorphism depicted below via Proposition 4.1.1 of \cite{Miranda strictifying operational coherences}. 
		
		$$\begin{tikzcd}
			\mathbf{st}_{3}\left(\mathfrak{A}\right) \arrow[r, "E_\mathfrak{A}"] & \mathfrak{A} \arrow[r, "W"] & \mathbf{Bicat} \arrow[r, "\mathbf{st}_{2}"] & \mathbf{Gray}
		\end{tikzcd}$$
		\item A $\mathbf{Gray}$-category $\mathfrak{B}$ has all trilimits if it has all trilimits of the following form.
		
		\begin{itemize}
			\item The shape $\mathfrak{A}$ is a cofibrant $\mathbf{Gray}$-category.
			\item The diagram $F: \mathfrak{A} \rightarrow \mathfrak{B}$ is a $\mathbf{Gray}$-functor.
			\item The weight $W: \mathfrak{A} \rightarrow \mathbf{Bicat}$ factors through $I\mathbf{st}_{2}: \mathbf{Bicat} \rightarrow \mathbf{Bicat}$.
		\end{itemize}
	\end{enumerate}
\end{proposition}

\begin{proof}
	We first note that since $\mathfrak{C}$ and $\mathbf{Gray}$ are $\mathbf{Gray}$-categories, there is no ambiguity in the trihomomorphisms in part (1). The proof of part (1) is as follows.
	
	\begin{align*}
		&\mathfrak{B} \text{     has     } <W{,} F> 
		\\
		\iff &\mathfrak{C} \text{     has     } <W{,} GF> &\text{Proposition }\ref{reduction to cofibrant targets} 
		\\
		\iff & \mathfrak{C} \text{     has     } <WE_{*}{,} \left(GF\right)E_{*}> &\text{Proposition }\ref{reduction to cofibrant shapes}
		\\
		\iff &  \mathfrak{C} \text{     has     } <\left(I\mathbf{st}_{2}\right)\left(WE_{*}\right){,} \left(GF\right)E_{*}> & \text{Proposition }\ref{reduction to pointwise cofibrant weights}
		\\
		\iff &\mathfrak{C} \text{     has     } <IW'{,}F'> & \text{Lemma }\ref{varying the weight and the diagram}
	\end{align*}
	
	\noindent Part (2) follows by applying part (1) to each pair consisting of a weight and diagram, and noting that $IW'$ and $F'$ are of the required form.
\end{proof}

\begin{remark}\label{summarising reduction to pointwise cofibrant weights on cofibrant Gray categories}
	At this stage we have reduced the study of tricategorical limits and colimits to special cases in which the shape and target are cofibrant $\mathbf{Gray}$-categories, and the trihomomorphisms involved are $\mathbf{Gray}$-functors. The weight may also be assumed to satisfy the following conditions. The first of these seems to be of most importance.
	\begin{itemize}
		\item For all $X \in \mathfrak{A}$, $WX$ is a cofibrant $2$-category, 
		\item For all $f \in \mathfrak{A}\left(X, Y\right)$, $Wf: WX \rightarrow WY$ is a $2$-functor between cofibrant $2$-categories, whose functor between underlying categories is also free.
		\item For all $2$-cells $\phi: f \Rightarrow g$ in $\mathfrak{A}$, $W\phi: Wf \Rightarrow Wg$ is a pseudonatural transformation whose $1$-cell components are all generators in $WY$.
	\end{itemize}  
	
	\noindent It remains to simplify the universal property to be one in terms of data in $\mathbf{Gray}$-enriched functor categories, rather than weaker data such as trinatural transformations or trimodifications. Section \ref{Weak 3 k transfors as higher cells of pseudoalgebras} uses three-dimensional monad theory to bridge this gap.
\end{remark}

\section{Simplifying the universal property}\label{Weak 3 k transfors as higher cells of pseudoalgebras}

\noindent Subsection \ref{Subsection Gray monads and their algebras} reviews the abstract notions of $\mathbf{Gray}$ monads and their $\mathbf{Gray}$-categories of pseudoalgebras. We extend existing theory in Proposition \ref{pseudoalgebra triadjunction} by exhibiting a triadjunction between the base $\mathbf{Gray}$-category and the $\mathbf{Gray}$-category of pseudoalgebras. Subsection \ref{subsection universal property with respect to Gray enriched 3 k transfors} will specialise to $\mathbf{Gray}$-monads whose categories of strict algebras are enriched functor categories. As shown in \cite{Buhne Gray Homo}, in this context the $\mathbf{Gray}$-category of pseudoalgebras is given by locally strict trihomomorphisms and general trinatural transformations or trimodifications between them. The $\mathbf{Gray}$-adjunction between strict algebras and pseudoalgebras will then be used in Theorem \ref{tricolimit reduction to strict 3 k transfors via LW} to characterise homs out of tricolimits as being trinaturally biequivalent to a hom in an enriched functor $\mathbf{Gray}$-category. Finally, Subsection \ref{explicit Gray monad T and weak morphism classifier L} will characterise the $\mathbf{Gray}$-enriched weights needed to express tricategorical notions in terms of algebraic and homotopical properties. This will be used to further simplify the universal property, resulting in a statement entirely expressible in the language of enrichment over $\mathbf{Gray}$ as a monoidal model category. In particular, we will not need to mention weaker data than those in enriched internal homs $[\mathfrak{A}, \mathfrak{B}]$. On the other hand, biequivalences, which are the weak-equivalences in $\mathbf{Gray}$, will still be used.

\subsection{$\mathbf{Gray}$-monads and their pseudoalgebras}\label{Subsection Gray monads and their algebras}

The aim of this Subsection is to describe the triadjunction between the pseudoalgebras of a $\Gr$-monad $(\mathfrak{B}, T, \eta, \mu)$ and the base $\Gr$-category $\mathfrak{B}$. As a consequence of this triadjunction, it follows that the left triadjoint preserves tricategorical colimits.

\begin{definition}\label{Definition Gray monads and pseudoalgebras}
	\hspace{1mm}
	\begin{itemize}
		\item 
		$\mathbf{Gray}$-$\mathbf{CAT}_{2}$ will denote the $2$-category of large $\mathbf{Gray}$-categories, $\mathbf{Gray}$-functors and $\mathbf{Gray}$-natural transformations. $\mathbf{Gray}\text{-}\mathbf{Cat}_{2}$ will denote the full sub-$2$-category on small $\mathbf{Gray}$-categories.
		\item Adjunctions in the $2$-categories $\mathbf{Gray}\text{-}\mathbf{CAT}_{2}$, or $\mathbf{Gray}\text{-}\mathbf{Cat}_{2}$ are called \emph{$\mathbf{Gray}$-adjunctions}.
		\item Monads in the $2$-categories $\mathbf{Gray}\text{-}\mathbf{CAT}_{2}$, or $\mathbf{Gray}\text{-}\mathbf{Cat}_{2}$ are called \emph{$\mathbf{Gray}$-monads} and will typically have their structure denoted as $\left(\mathfrak{B}, T, \eta, \mu\right)$.
	\end{itemize}
\end{definition}

\begin{definition}\label{triadjunction definition}
	A \emph{coherent triadjunction} consists of
	
	\begin{itemize}
		\item Two $\Gr$-categories $\mathfrak{A}$, and $\mathfrak{B}$.
		\item Two $\mathbf{Gray}$-functors $F: \mathfrak{A} \rightarrow \mathfrak{B}$ and $U: \mathfrak{B} \rightarrow \mathfrak{A}$,
		\item Two trinatural transformations $\eta: 1_\mathfrak{A} \rightarrow UF$ and $\varepsilon: FU \rightarrow 1_\mathfrak{B}$ which are semi-strict in the sense of Definition 3.2.3 of \cite{Miranda strictifying operational coherences},
		\item Two adjoint equivalence trimodifications $\rho: 1_{U} \Rightarrow U\varepsilon.\eta U$, $\lambda: \varepsilon F.F\eta \Rightarrow 1_{F}$,
		\item Two invertible perturbations $l$ and $r$. The perturbation $l$ has source depicted as below left and target $1_\varepsilon$. Meanwhile, $r$ has target depicted below right and source $1_{\eta}$.
	\end{itemize}
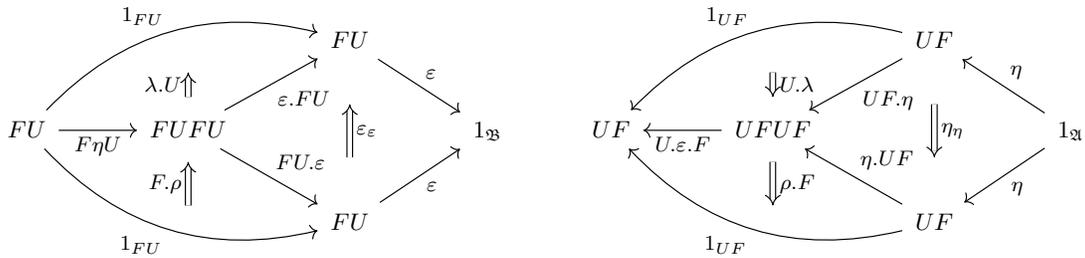
	$$\begin{tikzcd}[column sep = 15, font=\fontsize{9}{6}]
		&&{}\arrow[d, Leftarrow, shorten = 5, "\lambda.U"']
		&& FU \arrow[rrd,"\varepsilon"] \arrow[dd,Leftarrow,shorten=6mm,"\varepsilon_{\varepsilon}"]
		\\
		FU \arrow[rr,"F{\eta}U"']
		\arrow[rrrru, bend left, "1_{FU}"]
		\arrow[rrrrd, bend right, "1_{FU}"']
		&& FUFU
		\arrow[d, shorten = 5, Leftarrow, "F.\rho"']
		\arrow[rru,"\varepsilon.{FU}"'] \arrow[rrd,"FU.\varepsilon"]
		&&&& 1_\mathfrak{B} &{}
		\\
		&&{}
		&& FU \arrow[rru,"\varepsilon"']
	\end{tikzcd}	\begin{tikzcd}[column sep = 15, font=\fontsize{9}{6}]
		&&{}\arrow[d, Rightarrow, shorten = 5, "U.\lambda"]
		&& UF \arrow[rrd,leftarrow,"\eta"] \arrow[dd,Rightarrow,shorten=6mm,"\eta_{\eta}"]
		\\
		UF \arrow[rr,leftarrow, "U.\varepsilon.F"']
		\arrow[rrrru,leftarrow, bend left, "1_{UF}"]
		\arrow[rrrrd,leftarrow, bend right, "1_{UF}"']
		&& UFUF
		\arrow[d, shorten = 5, Rightarrow, "\rho.F"]
		\arrow[rru,leftarrow,"UF.\eta"'] \arrow[rrd,leftarrow,"\eta.UF"]
		&&&& 1_\mathfrak{A}
		\\
		&&{}
		&& UF \arrow[rru,leftarrow,"\eta"']
	\end{tikzcd}$$

	\noindent Subject to the requirement that the following two pastings of perturbations in the $2$-categories $\mathbf{Gray}$-$\mathbf{Cat}(\mathfrak{B}, \mathfrak{A})(F, F)$ and $\mathbf{Gray}$-$\mathbf{Cat}(\mathfrak{A}, \mathfrak{B})(U, U)$ are identities on $\lambda$ and $\rho$ respectively.

	$$\begin{tikzcd}[column sep=9, row sep = 7, font=\fontsize{9}{6}]
		&&&&
		{\varepsilon}F.F\eta
		\arrow[lllldddddddddddddddddddd,bend right=45, "1"']
		\arrow[dddd, "1.F\rho_F.1"]
		\arrow[rrrrdddddddddddddddddddd,bend left=45, "1"]
		\\
		\\
		\\
		\\
		&&&&
		{\varepsilon}F.FU{\varepsilon}F.F{\eta}UF.F{\eta}
		\arrow[ldddd, "1.1.F\eta_\eta"']
		\arrow[rdddd, "\varepsilon_{\varepsilon_{F}}.1.1"]
		\\
		\\
		\\
		\\
		{}\arrow[rrr, shorten = 15, shift right = 30, Rightarrow, "1.Fr"]
		&&&
		{\varepsilon}F.FU{\varepsilon}F.FUF{\eta}.F{\eta}
		\arrow[llldddddddddddd, "1.FU\lambda.1"]
		\arrow[rdddd, "\varepsilon_{\varepsilon_{F}}.1.1"']
		\arrow[rr, Rightarrow, shorten = 18, "{\left(\varepsilon_{\varepsilon_{F}}\right)}_{\left(F\eta_\eta\right)}"]
		&&
		{\varepsilon}F.{\varepsilon}FUF.F{\eta}UF.F{\eta}
		\arrow[ldddd, "1.1.F\eta_\eta"]
		\arrow[rrrdddddddddddd, "1.\lambda_{UF}.1"']
		\arrow[rrr, shorten = 15, shift right = 30, Rightarrow, "lF.1"]
		&&&{}
		\\
		\\
		\\
		\\
		&&&&
		{\varepsilon}F.{\varepsilon}FUF.FUF{\eta}.F{\eta}
		\arrow[dddd, "1.\varepsilon_{F\eta}.1"]
		\\
		\\
		\\
		\\
		&&{}\arrow[rr, Rightarrow, shorten = 40, "\varepsilon_\lambda.1"]
		&&
		{\varepsilon}F.F{\eta}.{\varepsilon}F.F{\eta}
		\arrow[rr, Rightarrow, shorten = 40, "1.\lambda_\eta"]
		\arrow[lllldddd, "\lambda.1.1"]
		\arrow[rrrrdddd, "1.1.\lambda"']
		&&{}
		\\
		\\
		\\
		\\
		{\varepsilon}F.F\eta
		\arrow[rrrrdddddd,bend right=15, "\lambda"']
		\arrow[rrrrrrrr,Rightarrow, shorten = 150, shift right = 8, "\lambda_\lambda"]
		&&&&&&&&
		{\varepsilon}F.F\eta
		\arrow[lllldddddd,bend left=15, "\lambda"]
		\\
		\\
		\\
		\\
		\\
		\\
		&&&&
		1_F
	\end{tikzcd}$$
	
	$$\begin{tikzcd}[column sep=9, row sep = 7, font=\fontsize{9}{6}]
		&&&&
		U\varepsilon.\eta U
		\arrow[lllldddddddddddddddddddd,leftarrow,bend right=45, "1"']
		\arrow[dddd, leftarrow,"1.U\lambda_{U}.1"]
		\arrow[rrrrdddddddddddddddddddd,bend left=45,leftarrow, "1"]
		\\
		\\
		\\
		\\
		&&&&
		U{\varepsilon}.U{\varepsilon}FU.UF{\eta}U.{\eta}U
		\arrow[ldddd, leftarrow,"U\varepsilon_\varepsilon.1.1"']
		\arrow[rdddd,leftarrow, "1.1.\eta_{\eta_{U}}"]
		\\
		\\
		\\
		\\
		{}\arrow[rrr, shorten = 15, shift right = 30, Leftarrow, "Ul.1"]
		&&&
		U{\varepsilon}.UFU{\varepsilon}.UF{\eta}U.{\eta}U
		\arrow[llldddddddddddd, leftarrow, "1.UF\rho.1"]
		\arrow[rdddd,leftarrow, "1.1.\eta_{\eta_U}"']
		\arrow[rr, Leftarrow, shorten = 18, "{\left(U\varepsilon_{\varepsilon}\right)}_{\left(\eta_{\eta_U}\right)}"]
		&&
		U{\varepsilon}.U{\varepsilon}FU.{\eta}UFU.{\eta}U
		\arrow[ldddd, leftarrow,"U\varepsilon_\varepsilon.1.1"]
		\arrow[rrrdddddddddddd, leftarrow, "1.\rho_{FU}.1"']
		\arrow[rrr, shorten = 15, shift right = 30, Leftarrow, "1.rU"]
		&&&{}
		\\
		\\
		\\
		\\
		&&&&
		U\varepsilon.UFU\varepsilon.\eta UFU.\eta U
		\arrow[dddd,leftarrow, "1.\eta_{U\varepsilon}.1"]
		\\
		\\
		\\
		\\
		&&{}\arrow[rr, Leftarrow, shorten = 40, "1.\eta_\rho"]
		&&
		U\varepsilon . \eta U.U\varepsilon .\eta U
		\arrow[rr, Leftarrow, shorten = 40, "\rho_\varepsilon .1"]
		\arrow[lllldddd, leftarrow,"1.1.\rho"]
		\arrow[rrrrdddd, leftarrow,"\rho.1.1"']
		&&{}
		\\
		\\
		\\
		\\
		U\varepsilon. \eta U
		\arrow[rrrrdddddd,bend right=15,leftarrow, "\rho"']
		\arrow[rrrrrrrr,Leftarrow, shorten = 150, shift right = 8, "\rho_\rho"]
		&&&&&&&&
		U\varepsilon.\eta U
		\arrow[lllldddddd,leftarrow,bend left=15, "\rho"]
		\\
		\\
		\\
		\\
		\\
		\\
		&&&&
		1_U
	\end{tikzcd}$$
\end{definition}

\begin{remark}\label{Remark coherent vs hom-wise triadjunctions}
	We explain some of the data in the pasting diagrams defining a coherent triadjunction, with reference to fragments of the semi-strictly generated closed structure of \cite{Miranda semi-strictly generated closed structure on Gray-Cat}.
	\begin{itemize}
		\item The invertible perturbations $\eta_{\rho}$ and $\varepsilon_{\lambda}$ are given as described in Proposition 3.2.6 of \cite{Miranda semi-strictly generated closed structure on Gray-Cat}.
		\item The invertible perturbations $\rho_\varepsilon$ and $\lambda_{\eta}$ are given as described in Proposition 3.3.1 of \cite{Miranda semi-strictly generated closed structure on Gray-Cat}.
		\item The trimodifications $\eta_\eta$, $\varepsilon_\varepsilon$, $\eta_{U\varepsilon}$ and $\varepsilon_{F\eta}$ are the left adjoint equivalences given as described in Proposition 3.2.3 of \cite{Miranda semi-strictly generated closed structure on Gray-Cat}.
		\item The invertible perturbations $\rho_\rho$, 	 $\lambda_\lambda$,	  ${\left(U\varepsilon_\varepsilon\right)}_{\left(\eta_{\eta_U}\right)}$ 
		and	   ${\left(F\eta_\eta\right)}_{\left(\varepsilon_{\varepsilon_{F}}\right)}$ are the interchangers in hom-$\mathbf{Gray}$-categories described in Remark 8.2.1 part (4e) of \cite{Miranda semi-strictly generated closed structure on Gray-Cat}.
	\end{itemize} 

\noindent By Lemma 3.2.1 of \cite{Campbell PhD}, a coherent triadjunction, or even just the equivalences $\lambda$ and $\rho$, suffice to construct a hom-wise triadjunction. 
\end{remark} 

\noindent For space reasons, we will not recall the definitions of pseudoalgebras and their higher dimensional cells, but assume familiarity with Definitions 13.4, 13.5, 13.6 and 13.7 of \cite{Gurski Coherence in Three Dimensional Category Theory}, as well as with the $\mathbf{Gray}$-category structure described in Section 13.3 of \cite{Gurski Coherence in Three Dimensional Category Theory}. Note that they introduce these concepts in the lax setting, whereas for our purposes data is pseudo. The difference is described explicitly in Definitions 13.8, 13.9 and 13.10 of \cite{Gurski Coherence in Three Dimensional Category Theory}.

\begin{notation}\label{notation for pseudoalgebras of Gray monads}
	\hspace{1mm}
	\begin{enumerate}
		\item The monadic adjunction associated to $\left(\mathfrak{B}, T, \eta, \mu\right)$ in $\mathbf{Gray}\text{-}\mathbf{Cat}_{2}$ will be denoted $F^{T} \dashv U^{T}$.
		\item The $\mathbf{Gray}$-category of pseudoalgebras, pseudomorphisms, pseudotransformations, and $3$-cells is described in detail in Sections 13.2 and 13.3 of \cite{Gurski Coherence in Three Dimensional Category Theory}, and will be denoted as $\mathbf{PsAlg}\left(\mathfrak{B}, T\right)$.
		\item Components of $k$-cells in $\mathbf{PsAlg}\left(\mathfrak{B}, T\right)$ that are denoted with subscripts will appear as lowest dimensional data in $\left(m+k\right)$-dimensional components of $\left(3, m\right)$-transfors featuring in the triadjunction of Proposition \ref{pseudoalgebra triadjunction}. Other $k$-cells in $\mathbf{PsAlg}\left(\mathfrak{B}, T\right)$ are denoted with superscripts. In more detail, our notation differs from that in Section 13.2 of \cite{Gurski Coherence in Three Dimensional Category Theory} in the following ways. 
		
		\begin{itemize}
			\item Pseudoalgebras will be denoted $\left(X, \varepsilon_{X}, \mathfrak{m}, u_{X}, \pi^{X}, \Lambda_{X}, \rho^{X}\right)$, with the only deviation from the notation of Definition 13.4 of \cite{Gurski Coherence in Three Dimensional Category Theory} being in the components $\varepsilon_{X}: TX \rightarrow X$, $u_{X}: \varepsilon_{X}.\eta \Rightarrow 1_{X}$ and capitalised $\Lambda$.
			\item When a strict algebra is being considered as a pseudoalgebra, its identity components will be omitted.
			\item Pseudomorphisms will be denoted as depicted below. This deviates from Definition 13.5 of \cite{Gurski Coherence in Three Dimensional Category Theory} as their $F$ is our $\varepsilon_{f}: f.\varepsilon_{X} \Rightarrow \varepsilon_{Y}.Tf$, their $\mathfrak{h}$ is our $u_{f}$ and their $\mathfrak{m}$ is our $\mathfrak{f}$.
			
			$$\left(f, \varepsilon_{f}, u_{f}, \mathfrak{f}\right): \left(X, \varepsilon_{X}, \mathfrak{m}^{X}, u_{X}, \pi^X, \Lambda_{X}, \rho\right) \rightarrow \left(Y, \varepsilon_{Y}, \mathfrak{m}^{Y}, u_{Y}, \pi^Y, \lambda_{Y}, \rho\right)$$ 
			\item Pseudotransformations $\left(f, \varepsilon_{f}, u_{f}, \mathfrak{f}\right) \Rightarrow \left(g, \varepsilon_{g}, u_{g}, \mathfrak{g}\right)$, defined in Definition 13.6 of \cite{Gurski Coherence in Three Dimensional Category Theory}, will be denoted as $\left(\phi, \varepsilon_{\phi}\right)$. Here $\phi: f \Rightarrow g$ is the $2$-cell component and $\varepsilon_{\phi}$ is the $3$-cell component, as depicted below.
			
			$$\begin{tikzcd}[column sep = 18, row sep = 15, font=\fontsize{9}{6}]
				TX
				\arrow[rr,bend left = 45, "Tf"]
				\arrow[dd,"\varepsilon_{X}"']
				&
				{}
				\arrow[d,Rightarrow, shorten = 2,"\varepsilon_{f}"]
				&
				TY
				\arrow[dd,"\varepsilon_Y"]
				\\
				&{}
				\\
				X
				\arrow[rr,bend left = 45,"f" {name = C}]
				\arrow[rr, bend right = 45, "g"' {name = D}]
				&
				{}
				&
				Y
				\arrow[from =C, to =D, Rightarrow, shorten = 10, "\phi"]
			\end{tikzcd}\begin{tikzcd}
				\\
				\cong^{\varepsilon_{\phi}}
			\end{tikzcd}	\begin{tikzcd}[column sep = 18, row sep = 15, font=\fontsize{9}{6}]
				TX
				\arrow[rr,bend left = 45, "Tf" {name = A}]
				\arrow[rr, bend right = 45, "Tg"' {name = B}]
				\arrow[dd,"\varepsilon_{X}"']
				&
				&
				TY
				\arrow[dd,"\varepsilon_{Y}"]
				\\
				&{}
				\arrow[d,Rightarrow, shorten = 5,"\varepsilon_{g}"]
				\\
				X\
				\arrow[rr, bend right = 45, "g"']
				&
				{}
				&
				Y &{}
				\arrow[from =A, to =B, Rightarrow, shorten = 10, shift right = 1, "T\phi"]
			\end{tikzcd} $$		
			
			\item $3$-cells of pseudoalgebras will be denoted just as in Definition 13.7 of \cite{Gurski Coherence in Three Dimensional Category Theory}.
		\end{itemize}
		\item The inclusion of strict algebras into pseudoalgebras will be denoted $i: \mathfrak{B}^{T} \rightarrow \mathbf{PsAlg}\left(\mathfrak{B}{,}T\right)$.
		\item If it exists, the left $\mathbf{Gray}$-adjoint to $i$ will be denoted $L$.
		\item $F:=$ \begin{tikzcd}
			\mathfrak{B} \arrow[r, "F^{T}"] & \mathfrak{B}^{T} \arrow[r, "i"] & \mathbf{PsAlg}\left(\mathfrak{B}{,}T\right)
		\end{tikzcd}
	\end{enumerate}
\end{notation}

\begin{proposition}\label{pseudoalgebra triadjunction}
	Let $\left(\mathfrak{B}, T, \eta, \mu\right)$ be a $\mathbf{Gray}$-monad.
	\begin{enumerate}
		\item There is a $\mathbf{Gray}$-functor $U: \mathbf{PsAlg}\left(\mathfrak{B}{,}T\right) \rightarrow \mathfrak{B}$ which sends a $k$-cell to its underlying $k$-cell in $\mathfrak{B}$.
		\item There is a semi-strict trinatural transformation $\varepsilon: FU \rightarrow 1_{\mathbf{PsAlg}\left(\mathfrak{B}{,}T\right)}$ defined in the following way.
		\begin{itemize}
			\item The component on a pseudoalgebra $\left(X, \varepsilon_{X}, u_{X}, \mathfrak{m}^{X}, \Lambda_{X}, \rho^{X}, \pi^{X}\right)$ is given by the pseudomorphism depicted below. 
			
			$$\left(\varepsilon_{X}, \mathfrak{m}^{X}, \rho^{X}, \pi^{X}\right): \left(TX, \mu_{X}\right) \rightarrow \left(X, \varepsilon_{X}, u_{X}, \mathfrak{m}^{X}, \Lambda_X, \rho^X, \pi^X\right)$$
			\item The component on a pseudomorphism $\left(f, \varepsilon_{f}, u_{f}, \mathfrak{f}\right)$ is the pseudotransformation $\left(\varepsilon_{f}, \mathfrak{f}\right)$.
			\item The component on a pseudotransformation $\left(\phi, \varepsilon_{\phi}\right)$ is $\varepsilon_{\phi}$.
		\end{itemize}
		\item There is an adjoint equivalence of trinatural transformations whose left adjoint is given by the trimodification $u: 1_{U} \Rightarrow U\varepsilon.\eta_{U}$ defined by mapping $k$-cells to their component denoted with $u$ as per Notation \ref{notation for pseudoalgebras of Gray monads}.
		\item The trinatural transformation $\varepsilon$ whiskers with $F$ and $U$ as $U\varepsilon_{F} = \mu$.
		\item There is an equation of semi-strict trinatural transformations ${\varepsilon}_{F}.F\eta = 1_{{F}}$.
		\item The restriction of $u$ along ${F}$ is an identity.
		\item The assignment which maps a pseudoalgebra $\left(X, \varepsilon_{X}, u_{X}, \pi^{X}, \Lambda_{X}, \rho^{X}\right)$ to $\Lambda_{X}$ defines an invertible perturbation from the trimodification depicted below to the identity trimodification on $\varepsilon$.

		$$\begin{tikzcd}[font=\fontsize{9}{6}, row sep = 15, column sep = 15]
			&&{}\arrow[dd, Rightarrow, shorten = 5, "Fu"]
			&&FU
			\arrow[rrdd, "\varepsilon"]
			\arrow[dddd, Rightarrow, shorten = 20, "\varepsilon_{\varepsilon}"]
			\\
			\\
			FU
			\arrow[rr, "F\eta U"]
			\arrow[rrrruu, bend left, "1_{FU}"]
			\arrow[rrrrdd, bend right, "1_{FU}"']
			&& FUFU
			\arrow[rrdd, "\varepsilon_{FU}"]
			\arrow[rruu, "FU\varepsilon"'] &&&&1_{\mathbf{PsAlg}\left(\mathfrak{B}{,}T\right)}
			\\
			&&=
			\\
			&&&&FU
			\arrow[rruu, "\varepsilon"']
		\end{tikzcd}$$
		
		\item These data form a triadjunction as described in Definition \ref{triadjunction definition}. 
		\item $F: \mathfrak{B} \to \mathbf{PsAlg}(\mathfrak{B}, T)$ preserves tricategorical colimits.
	\end{enumerate} 
\end{proposition}

\begin{proof}
	Part (1) is clear from the description of the $\mathbf{Gray}$-category structure on $\mathbf{PsAlg}\left(\mathfrak{B}, T\right)$ given in Section 13.3 of \cite{Gurski Coherence in Three Dimensional Category Theory}. The rest of the proofs are mostly an exercise in unwinding definitions and matching axioms to properties which need to be proved. Strictness and the axioms for the $\mathbf{Gray}$-monad $\left(\mathfrak{B}, T, \eta, \mu\right)$ make the proofs substantially easier than they would otherwise be. We focus on the properties of the data in $\mathbf{PsAlg}\left(\mathfrak{B}, T\right)$ which are needed in the proofs.
	\\
	\\
	\noindent For part (2), we first check that the components of $\varepsilon$ are well-defined as data in $\mathbf{PsAlg}\left(\mathfrak{B}, T\right)$. For $1$-cell components on objects, axioms (1), (2) and (3) for pseudomorphisms use axioms (1), (2) and (3) for pseudoalgebras. For $2$-cell components on morphisms, axioms (1) and (2) for pseudotransformations follow from axioms (3) and (1) for pseudomorphisms respectively. Finally, the axiom for $3$-cells uses axiom (2) for pseudotransformations. Naturality in $3$-cells for $\varepsilon$ uses the axiom for $3$-cells, while semi-strictness of $\varepsilon$ follows from the $\mathbf{Gray}$-category structure.
	\\
	\\
	\noindent For part (3), the local modification axiom for $u$ follows from axiom (1) for the pseudotransformation on which it needs to be checked. The unit and composition axioms are the immediate from the $\mathbf{Gray}$-category structure of $\mathbf{PsAlg}\left(\mathfrak{B}, T\right)$. Parts (4), (5) and (6) all use the structure of free algebras $\left(TX, \mu_{X}\right)$. Part (4) uses the fact that the $1$-cell components of free algebras are given by the multiplication $\mu$, while parts (5) and (6) use the unit laws for $\left(\mathfrak{B}, T, \eta, \mu\right)$. For part (7) the perturbation axiom follows from axiom (2) for pseudomorphisms.
	\\
	\\
	\noindent For the first axiom for triadjunctions in Example \ref{triadjunction definition}, it is easy to see that every perturbation is in fact the identity, using strictness of the $\mathbf{Gray}$-monad and of free algebras. In particular, the trimodification corresponding to $\lambda$ in Example \ref{triadjunction definition} is the identity by part (5), $F\eta_\eta$ is the identity since $\eta$ is $\mathbf{Gray}$-natural, the perturbation corresponding to $r$ is the identity by part (6), and the perturbation corresponding to $l_F$ is the identity by part (7) and the fact that the component $\Lambda_{X}$ is the identity for free algebras. The second axiom for triadjunctions follows from axiom (4) for pseudoalgebras. Finally, part (9) follows from part (8) by Lemma \ref{left triadjoint preserves tricolimits}.
\end{proof}

\subsection{Universal property up to a trinatural biequivalence with an enriched functor category}\label{subsection universal property with respect to Gray enriched 3 k transfors}

\noindent Lemma \ref{Tricat reflection into GrayHom} recalls some important facts about the $\mathbf{Gray}$-monads for enriched functor $\mathbf{Gray}$-categories $[\mathfrak{A}^\text{op}, \mathbf{Gray}]$, and their associated $\mathbf{Gray}$-categories of pseudoalgebras. Combining the universal property of Lemma \ref{Tricat reflection into GrayHom} part (2) with Lemma \ref{strictifying components of a trinatural transformation into Hom}., we will be able to prove Theorem \ref{tricolimit reduction to strict 3 k transfors via LW}. This will further reduce tricategorical limits and colimits to the point of being able to recognise homs out of tricolimits (resp. into trilimits) as being trinaturally biequivalent to enriched functor categories.

\begin{definition}
	A trihomomorphism between $\Gr$-categories is called \emph{locally strict} if it is given by $2$-functors between hom-$2$-categories.
\end{definition}

\begin{lemma}\label{Tricat reflection into GrayHom}
	Let $\mathfrak{A}$ be a $\mathbf{Gray}$-category, and let $O:\mathfrak{A}_{0}\rightarrow \mathfrak{A}^\text{op}$ denote the $\mathbf{Gray}$-functor which includes the set of objects.
	
	\begin{enumerate}
		\item The functor $[O, \mathbf{Gray}]: [\mathfrak{A}^\text{op}, \mathbf{Gray}] \rightarrow [\mathfrak{A}_{0}, \mathbf{Gray}]$ has both adjoints, and is monadic in the $2$-category $\mathbf{GRAY}$-$\mathbf{CAT}_{2}$ of enriched categories, enriched functors and enriched natural transformations.
		\item The full sub-$\mathbf{Gray}$-category of $\mathbf{TRICAT}\left(\mathfrak{A}^\text{op}, \mathbf{Gray}\right)$ on locally strict trihomomorphisms, is the $\mathbf{Gray}$-category of pseudoalgebras $\mathbf{PsAlg}_{T}$ of Proposition 13.15 of \cite{Gurski Coherence in Three Dimensional Category Theory}.
		\item The inclusion $i: [\mathfrak{A}^\text{op}, \mathbf{Gray}] \rightarrow \mathbf{TRICAT}_\text{ls}[\mathfrak{A}^\text{op}, \mathbf{Gray}]$ has a left adjoint $L$ in $\mathbf{GRAY}$-$\mathbf{CAT}_{2}$.
		\item Let $\left(T, \eta, \mu\right)$ denote the monad induced on $[\mathfrak{A}_{0}, \mathbf{Gray}]$. Then $L$ is given on objects by mapping $W \in \mathbf{TRICAT}_\text{ls}\left(\mathfrak{A}^\text{op}, \mathbf{Gray}\right)$ to the codescent object as in Definition 11.7 of \cite{Gurski Coherence in Three Dimensional Category Theory}, of the $\mathbf{Gray}$-functor $\underline{W}$ described on objects below, and in full in Proposition 14.2 of \cite{Gurski Coherence in Three Dimensional Category Theory}. 
		
		$$\begin{tikzcd}
			\Delta_\text{ps}^{G}
			\arrow[rr, "\underline{W}"]
			&& {[}\mathfrak{A}^\text{op}{,} \mathbf{Gray}{]}
			\\
			{[}n{]}
			\arrow[rr, mapsto]
			&& \left({T}^{n}W{,}\mu_{T^{n-1}W}\right)
		\end{tikzcd}$$ 
		\item The unit of $L \dashv i$ is a biequivalence internal to $\mathbf{TRICAT}_\text{ls}\left(\mathfrak{A}^\text{op}, \mathbf{Gray}\right)$.
	\end{enumerate}
\end{lemma}

\begin{proof}
	Part (1) is Theorem 1.4 of \cite{Buhne PhD}, in the case $\mathcal{L} = \mathcal{V} = \mathbf{Gray}$. Part (2) is Theorem 1.6 of \cite{Buhne PhD}. Parts (3), (5) follow from Theorem 14.10 and Corollary 15.14 of \cite{Gurski Coherence in Three Dimensional Category Theory}, noting that the base $[\mathfrak{A}_{0}, \mathbf{Gray}]$ is $\mathbf{Gray}$-cocomplete and $T$ is $\mathbf{Gray}$-cocontinuous. The specific case of this $\mathbf{Gray}$-monad $T$ also appears in Theorem 8 of \cite{Buhne Gray Homo}. For part (4), the explicit description of $L$ is given at the start of Section 14.1 of \cite{Gurski Coherence in Three Dimensional Category Theory}.
\end{proof}

\noindent An explicit description of the monad $\left(T, \eta, \mu\right)$ can be derived from the general enriched category theory setting, while an explicit description of $L$ can be extracted from Chapter 14 of \cite{Gurski Coherence in Three Dimensional Category Theory}. These are discussed in Subsection \ref{explicit Gray monad T and weak morphism classifier L} as they are needed to understand the triessential image of $L$, which is tri-fully faithful by Lemma \ref{Tricat reflection into GrayHom} part (5). We first use Lemma \ref{Tricat reflection into GrayHom} part (3) to establish that homs out of tricolimits (resp. into trilimits) are trinaturally biequivalent to data involving $\mathbf{Gray}$-enriched transfors. Later we will improve this to a $\mathbf{Gray}$-natural biequivalence. On top of Lemma \ref{Tricat reflection into GrayHom} part (3), Theorem \ref{tricolimit reduction to strict 3 k transfors via LW} will also require Lemma \ref{strictifying components of a trinatural transformation into Hom}, which allows certain trinatural transformations into $\mathbf{Bicat}$ to be `component-wise strictified'.

\begin{lemma}\label{strictifying components of a trinatural transformation into Hom}
	Let $\mathfrak{A}$ be a $\mathbf{Gray}$-category, and let $W, G: \mathfrak{A} \rightarrow \mathbf{Bicat}$ be strict homomorphisms of tricategories. Suppose that for every $f: X\rightarrow Y \in \mathfrak{A}$, $Wf: WX \rightarrow WY$ is a $2$-functor between cofibrant $2$-categories and $Gf: GX \rightarrow GY$ is a $2$-functor. Consider the $2$-fully faithful map $S$ depicted below. Here $\mathbf{TRICAT}_\text{x}\left(\mathfrak{A}, \mathbf{Bicat}\right)\left(W, G\right)$ for $\text{x} \in \{\text{strict}, \text{pseudo}\}$ have the same objects, $2$-cells and $3$-cells, but differ on morphisms. Specifically, $S$ includes those trinatural transformations whose $1$-cell components are $2$-functors into the larger tricategory whose morphisms are arbitrary trinatural transformations $p: W \rightarrow G$, and hence may have $1$-cell components given by pseudofunctors.
	
	$$S: \mathbf{TRICAT}_\text{strict}\left(\mathfrak{A}, \mathbf{Bicat}\right)\left(W, G\right) \hookrightarrow \mathbf{TRICAT}_\text{pseudo}\left(\mathfrak{A}, \mathbf{Bicat}\right)\left(W, G\right)$$
	
	\begin{enumerate}
		\item $S$ is a biequivalence.
		\item $S$ is strictly natural as $G$ varies along any $g: G \rightarrow G'$ where $G'X$ are also $2$-categories and $g_{X}$ are also $2$-functors.
	\end{enumerate}
\end{lemma}

\begin{proof}
	For part (1) it suffices to check biessential surjectivity on objects, so let $p: W \rightarrow G$ be an arbitrary trinatural transformation. We need to show that $p$ is equivalent to one whose $1$-cell components are $2$-functors, rather than merely pseudofunctors. By Remark 2.3.4 of \cite{Miranda strictifying operational coherences}, since $WX$ is a cofibrant $2$-category and $GX$ is a $2$-category, there is an invertible icon $e_{X}: p_{X} \rightarrow \overline{p}_{X}$ for each $X \in \mathfrak{A}$, with $\overline{p}_{X}$ being a $2$-functor. The assignment $X \mapsto \overline{p}_{X}$ extends to a trinatural transformation by assigning $f: X \rightarrow Y$ to the pseudonatural transformation depicted below. 
	
	$$\begin{tikzcd}
		\overline{p}_{Y}.Wf \arrow[r, "e_{X}^{-1}.1"] & p_{Y}.Wf \arrow[r, "p_{f}"] & Gf.p_{X} \arrow[r, "1.e_{X}"] & Gf.\overline{p}_{X}
	\end{tikzcd}$$ 
	
	\noindent The free and operational $3$-cell coherences of $p$ are whiskered by $e_{X}$ in the evident way to produce the analogous components for $\overline{p}$, possibly also cancelling $e_{X}$ with its inverse. The trinaturality axioms follow from those for $p$. The assignment $X \mapsto e_{X}$ extends to a strict trimodification, again by cancelling $e_{X}$ with its inverse. The trimodification axioms are then easy to check. Part (2) is also an easy inspection.
\end{proof}

\noindent In Theorem \ref{tricolimit reduction to strict 3 k transfors via LW} to follow, the simplifying assumptions on $F$ and $\mathfrak{A}$ do not come at the cost of generality, since by the results of Subsection \ref{subsection reduction to simpler shapes targets weights and diagrams} more general data can always be replaced with data of this form. Indeed, a more elaborate statement without any simplifying assumptions is also possible. On the other hand, we allow $W$ to be more general than needed by Proposition \ref{special weights and diagrams for tricolimits} since it does not make the statement any more complicated. Theorem \ref{tricolimit reduction to strict 3 k transfors via LW} part (1) will be improved upon by Corollary \ref{Reduction to Gray natural biequivalence}, which will characterise tricategorical limits and colimits entirely in the language of $\mathbf{Gray}$ as a monoidal model category.

\begin{theorem}\label{tricolimit reduction to strict 3 k transfors via LW}
	Let $\mathfrak{A}$ be a cofibrant $\mathbf{Gray}$-category, $W: \mathfrak{A}^\text{op} \rightsquigarrow \mathbf{Gray}_{c}$ be a locally strict trihomomorphism, let \begin{tikzcd}
		\mathbf{Gray}_{c} \arrow[r, "J"] & \mathbf{Gray} \arrow[r, "I"] & \mathbf{Bicat}
	\end{tikzcd} denote the inclusions, and let $L \dashv i$ be the $\mathbf{Gray}$-adjunction of Lemma \ref{Tricat reflection into GrayHom} part (3). Let $F: \mathfrak{A} \rightarrow \mathfrak{B}$ be a $\mathbf{Gray}$-functor.
	\begin{enumerate}
		\item If $\mathfrak{B}$ is hom-wise cofibrant then an object $IW\odot F \in \mathfrak{B}$ is a tricolimit if and only if there is a biequivalence of the form $\mathfrak{B}\left(IW\odot F, ? \right) \rightarrow [\mathfrak{A}^\text{op}, \mathbf{Gray}]\left(LW, \mathfrak{B}\left(F-, ?\right)\right)$ in the $\mathbf{Gray}$-category $\mathbf{TRICAT}_{\text{s}}\left(\mathfrak{B}, \mathbf{Gray}\right)$.
		\item Suppose $\mathfrak{B}$ is a hom-wise cofibrant $\mathbf{Gray}$-category that has $\mathbf{Gray}$-enriched limits (resp. colimits) for all weights $Li\left(W'\right)$ where $W'$ is pointwise cofibrant and the shape $\underline{\mathfrak{A}}$ is cofibrant. Then $\mathfrak{B}$ has all trilimits  (resp. all tricolimits) with defining biequivalence in $\mathbf{TRICAT}_{s}\big(\mathfrak{B}, \mathbf{Gray}\big)$ (resp. $\mathbf{TRICAT}_{s}\big(\mathfrak{B}^\text{op}, \mathbf{Gray}\big)$).
		\item Parts (1) and (2) remain true without the assumption that $\mathfrak{B}$ is hom-wise cofibrant, except that in part (1) the maps $\mathfrak{B}\left(IW\odot F, Y \right) \rightarrow [\mathfrak{A}^\text{op}, \mathbf{Gray}]\left(LW, \mathfrak{B}\left(F-, Y\right)\right)$ may be pseudofunctors, so that the biequivalence is now in $\mathbf{TRICAT}_{s}\big(\mathfrak{B}, \mathbf{Bicat}\big)$.
	\end{enumerate}
\end{theorem}

\begin{proof}
	Part (1) follows from the following chain of biequivalences, which are all trinatural in $Y$ either by definition, by Yoneda, or by triadjointness.
	{\small
		\begin{align*}
			\mathfrak{B}\left(IJW\odot F{,}Y\right) &\sim\mathbf{TRICAT}\left(\mathfrak{A}^\text{op}{,}\mathbf{Bicat}\right)\left(IJW{,}I\mathfrak{B}\left(F-{,}Y\right)\right) &\text{Definition }\ref{trilimits and tricolimits definition} 
			\\
			& \sim \mathbf{TRICAT}\left(\mathfrak{A}^\text{op}{,}\mathbf{Gray}_\text{c}\right)\left(W{,}\mathfrak{B}\left(F-{,}Y\right)\right) &\text{Proposition }\ref{Hom triequivalent to Gray c} \text{ part (2)}
			\\
			& = \mathbf{TRICAT}_\text{ls}\left(\mathfrak{A}^\text{op}{,}\mathbf{Gray}\right)\left(JW{,}i\mathfrak{B}\left(F-{,}Y\right)\right) &J\text{ fully faithful,}
			\\
			&&W\text{ locally strict}
			\\
			&&\mathfrak{B} \text{ hom-wise cofibrant}
			\\
			&\cong {[}\mathfrak{A}^\text{op}{,}\mathbf{Gray}{]}\left(L\left(JW\right){,}\mathfrak{B}\left(F-{,}Y\right)\right) & \text{Lemma }\ref{Tricat reflection into GrayHom}\text{ part (3)}
		\end{align*}
	}
	\noindent Part (2) follows from part (1) and Proposition \ref{special weights and diagrams for tricolimits}. Given a general diagram $F : \mathfrak{A} \rightsquigarrow \mathfrak{B}$ and weight $W: \mathfrak{A}^\text{op} \rightsquigarrow \mathbf{Bicat}$, we first use Proposition \ref{special weights and diagrams for tricolimits} to replace these data with $\mathbf{Gray}$-functors $W': {{\mathfrak{A}}}^\text{op} \rightarrow \mathbf{Gray}$ and $F': {\mathfrak{A}} \rightarrow \mathfrak{B}$ where $\underline{\mathfrak{A}}$ is a cofibrant $\mathbf{Gray}$-category and $W'$ is pointwise cofibrant. Then the $\mathbf{Gray}$-enriched colimit $Li\left(W'\right) \cdot F'$ provides a tricolimit by part (1). The trilimits case is dual. For part (3), observe that hom-wise cofibrancy of $\mathfrak{B}$ is only used in the second biequivalence. Instead, Lemma \ref{strictifying components of a trinatural transformation into Hom} can be used to give a different trinatural biequivalence, taking $G:= \mathfrak{B}\left(F-, Y\right)$.
\end{proof}

\noindent We note that part (2) of Theorem \ref{tricolimit reduction to strict 3 k transfors via LW} is similar to Proposition 3.3.1 of \cite{Campbell PhD}, and the proof techniques are also similar. The main difference is that Proposition 3.3.1 of \cite{Campbell PhD} makes a milder assumption on $\mathfrak{A}$, namely that the hom-$2$-categories are cofibrant, and makes a stronger assumption on $\mathfrak{B}$, namely that all $\mathbf{Gray}$-enriched limits exist.
\\
\\
\noindent Section \ref{Section reduction to strict pointwise cofibrant weights on cofibrant Gray categories} reduced the study of tricategorical limits and colimits to certain types of weights and diagrams on cofibrant shapes. This Subsection has shown that if the weight is adjusted further, by taking its image under $L$, then the hom out of a tricolimit (resp. into a trilimit) is biequivalent to a $\mathbf{Gray}$-category of $\mathbf{Gray}$-enriched transfors rather than weaker ones. It remains to characterise the kinds of weights in the triessential image of $L$, and to further simplify the biequivalence  $\mathfrak{B}\left(IW\odot F, ? \right) \rightarrow [\mathfrak{A}^\text{op}, \mathbf{Gray}]\left(LW, \mathfrak{B}\left(F-, ?\right)\right)$ from one internal to $\mathbf{TRICAT}_{s}\left(\mathfrak{A}^\text{op}, \mathbf{Gray}\right)$ to one internal to $[\mathfrak{A}^\text{op}, \mathbf{Gray}]$. This will be the aim of Subsection \ref{explicit Gray monad T and weak morphism classifier L} to follow.

\subsection{Projective cofibrancy and flexibility}\label{explicit Gray monad T and weak morphism classifier L}

Remark \ref{Explicit description of Gray monad} gives an explicit description of the $\mathbf{Gray}$-monad on $[\mathfrak{A}_{0}^\text{op}, \mathbf{Gray}]$ whose strict algebras are $\mathbf{Gray}$-enriched presheaves on $\mathfrak{A}$. After this we give an explicit description of the reflection  $L:\mathbf{Tricat}_\text{ls}\left(\mathfrak{A}^\text{op}, \mathbf{Gray}\right) \rightarrow [\mathfrak{A}^\text{op}, \mathbf{Gray}]$ of pseudoalgebras into strict ones from Lemma \ref{Tricat reflection into GrayHom} part (3). Part of the aim of this Subsection is to understand this reflection, when applied to objects which are already $\mathbf{Gray}$-functors, as a projective cofibrant replacement. Theorem \ref{trinatural transformation classifier cofibrancy} relates this homotopical property to an algebraic property, namely flexibility, and also relates general weights for tricategorical limits and colimits to flexible $\mathbf{Gray}$-enriched weights on strictified shapes. Finally, Corollary \ref{Reduction to Gray natural biequivalence} strengthens the sense in which flexible weights on hom-wise cofibrant $\mathbf{Gray}$-categories suffice to model tricategorical limits and colimits, using only notions expressible in terms of enrichment over $\mathbf{Gray}$ as a monoidal model category.

\begin{definition}\label{definition flexible weight}
	Let $\left(\mathfrak{B}, T\right)$ be a $\mathbf{Gray}$-monad for which the $\mathbf{Gray}$-adjunction $L\dashv i$ of \cite{Gurski Coherence in Three Dimensional Category Theory} Theorem 14.10 exists. Then a strict algebra $W \in \mathfrak{B}^{T}$ will be called \emph{flexible} if the counit $\varepsilon_{W}: Li\left(W\right) \rightarrow W$ has a section in the $\mathbf{Gray}$-category $\mathfrak{B}^{T}$.
\end{definition}

\noindent We begin with some examples of projective cofibrant weights in Example \ref{representables are flexible}, and relate projective cofibrancy to pointwise cofibrancy in Proposition \ref{projective and homwise implies pointwise}. Following this, we give the explicit descriptions of $T$ and $L$.

\begin{example}\label{representables are flexible}
	Representables $\mathfrak{A}\left(-, X\right): \mathfrak{A}^\text{op} \rightarrow \mathbf{Gray}$ are projective with respect to any $\mathbf{Gray}$-natural transformation that is pointwise surjective on objects. In particular, they are projective cofibrant. Indeed, by the $\mathbf{Gray}$-enriched Yoneda lemma, the lifting problem on the left is equivalent to the lifting problem on the right. Such a lifting exists if $t_{X}$ is surjective on objects.

	$$\begin{tikzcd}
		&E\arrow[d, "t"]
		\\
		\mathfrak{A}\left(-, X\right)\arrow[r, "v"'] \arrow[ru, dashed, bend left]
		&B&{}
	\end{tikzcd}\begin{tikzcd}
		&&E\left(X\right)
		\arrow[d, "t_{X}"]
		\\
		\mathbf{1}
		\arrow[rru, dashed]
		\arrow[r, "``1_{X}{"}"'] &\mathfrak{A}\left(X{,}X\right)
		\arrow[r, "v_{X}"']
		&B\left(X\right)
	\end{tikzcd}$$
	
	\noindent It is immediate that the weights for coproducts, $2$-cell coinserters, $3$-cell coinserters and $3$-cell coequifiers are projective cofibrant, since these are copowers of representables by generating cofibrations. The proof of Lemma 5.2.1 of \cite{Miranda Enriched Kleisli objects for pseudomonads} adapts to show that the weight for $3$-cell coinverters is also projective cofibrant, while by Corollary 5.2.5 of \cite{Miranda Enriched Kleisli objects for pseudomonads}, the comparatively more complicated weight of $\mathbf{Gray}$-enriched Kleisli objects is also projective cofibrant. Section 12.3 of \cite{Gurski Coherence in Three Dimensional Category Theory} can be understood as saying not only that $\mathbf{Gray}$-enriched codescent objects are certain weighted colimits, but that moreover their weight is projective cofibrant.
\end{example}

\begin{proposition}\label{projective and homwise implies pointwise}
	Let $\mathfrak{A}$ be a hom-wise cofibrant $\mathbf{Gray}$-category. If $W: \mathfrak{A}^\text{op} \rightarrow \mathbf{Gray}$ is projective cofibrant then it is also poinwise cofibrant. 
\end{proposition}

\begin{proof}
	Projective cofibrant objects are constructed via certain retracts and colimits of copowers of representables by generating cofibrations. These colimits are computed pointwise in $\mathbf{Gray}$, and preserve cofibrancy.
\end{proof}

\begin{remark}\label{Explicit description of Gray monad}
	We give an explicit description of the $\mathbf{Gray}$-monad $\left(T, \eta, \mu\right)$ on $[\mathfrak{A}_{0}, \mathbf{Gray}]$ induced by the $\mathbf{Gray}$-adjunction of Lemma \ref{Tricat reflection into GrayHom} part (1), and further unwind the description of the reflection $L: \mathbf{TRICAT}_\text{ls}\left(\mathfrak{A}^\text{op}, \mathbf{Gray}\right) \rightarrow [\mathfrak{A}^\text{op}, \mathbf{Gray}]$ of Lemma \ref{Tricat reflection into GrayHom} parts (3) and (4). These calculations will help show that the $\mathbf{Gray}$-adjunction $L \dashv i$ restricts to a triequivalence by taking weights which are projective cofibrant, when $\mathfrak{A}$ is hom-wise cofibrant. We will use \textbf{bold} to distinguish data in the $\mathbf{Gray}$-category $\mathfrak{A}$ from data in $2$-categories in the image of $W_{0}: \mathfrak{A}_{0} \rightarrow \mathbf{Gray}$.
	\\
	\\
	\noindent Recall that $F^{T}$ is given by left Kan extension along $O: \mathfrak{A}_{0} \rightarrow \mathfrak{A}$. Thus, given a family of $2$-categories $W_{0}: \mathfrak{A}_{0} \rightarrow \mathbf{Gray}$, the underlying family of $2$-categories of the free algebra $F^{T}\left(W_{0}\right)$ is given as depicted below. 
	
	$$T\left(W_{0}\right):=\begin{tikzcd}
		\mathfrak{A}_{0}\arrow[r, "O"]
		& \mathfrak{A}^\text{op} \arrow[r, "\mathbf{Lan}_{O}\left(W_{0}\right)"]
		& \mathbf{Gray}
		\\
		\mathbf{A} \arrow[rr, mapsto, "F^{T}\left(W_{0}\right)\left(\mathbf{A}\right)"]
		&& \sum\limits_{\mathbf{B} \in \mathfrak{A}_{0}} \mathfrak{A}\left(\mathbf{A}{,} \mathbf{B}\right) \otimes W_{0}\left(\mathbf{B}\right) 
	\end{tikzcd}$$  
	
	\noindent For example, if $W_{0}=\mathbf{1}$ on an object $\mathbf{A} \in \mathfrak{A}$ but is $\mathbf{0}$ on all other objects, then $T\left(W_{0}\right)$ is the representable on $\mathbf{A}$. The behaviour of $T\left(W\right)$  between hom-$2$-categories is adjoint transpose to the $2$-functor depicted below.
	
	$$\begin{tikzcd}
		\sum\limits_{\mathbf{C}\in \mathfrak{A}_{0}} \mathfrak{A}\left(\mathbf{B}{,}\mathbf{C}\right)\otimes \mathfrak{A}\left(\mathbf{A}{,}\mathbf{B}\right)\otimes W\left(\mathbf{C}\right) \arrow[rrrr, "\sum\limits_{\mathbf{C}\in \mathfrak{A}_{0}}\big( \circ_{\mathbf{A}{,}\mathbf{B}{,}\mathbf{C}}\otimes W\left(\mathbf{C}\right)\big)"]
		&&&& \sum\limits_{\mathbf{C}\in \mathfrak{A}_{0}} \mathfrak{A}\left(\mathbf{A}{,}\mathbf{C}\right)\otimes W\left(\mathbf{C}\right)
	\end{tikzcd}$$
	
	\noindent The unit and multiplication of $\left(T, \eta, \mu\right)$ are given in terms of the identities and composition in $\mathfrak{A}$. The unit is displayed below as an example, where we denote the coproduct coprojection by $\pi$. It is clear from this description that if $\mathfrak{A}$ has cofibrant homs each $W_{0}\left(\mathbf{A}\right)$ is cofibrant then $F^T \left(W_{0}\right)$ is also pointwise cofibrant, since cofibrant $2$-categories are closed under $\mathbf{Gray}$-tensor products and coproducts. 
	
	$$\begin{tikzcd}
		W\left(\mathbf{A}\right) \arrow[rr, "\left(1_{\mathbf{A}}{,}-\right)"]
		&& \mathfrak{A}\left(\mathbf{A}{,}\mathbf{A}\right)\otimes W\left(\mathbf{A}\right) \arrow[rr, "\pi_{\mathbf{A}}"]
		&&\sum\limits_{\mathbf{B} \in \mathfrak{A}} \mathfrak{A}\left(\mathbf{A}{,}\mathbf{B}\right)\otimes W\left(\mathbf{B}\right)
	\end{tikzcd}$$
\end{remark}

\noindent A similar analysis can be made with $\mathbf{Gray}$ replaced with any $\mathbf{Gray}$-category $\mathfrak{B}$ admitting coproducts and coowers by the hom-$2$-categories of $\mathfrak{A}$. Each instance of a $\mathbf{Gray}$-tensor is replaced by the copower in $\mathfrak{B}$.

\begin{corollary}\label{Gray monad for functors on codomain with copowers by homs of domain}
	Let $\mathfrak{A}$ and $\mathfrak{B}$ be $\mathbf{Gray}$-categories with $\mathfrak{A}$ small. Suppose $\mathfrak{B}$ has coproducts and copowers by all hom-$2$-categories of $\mathfrak{A}$.
	\begin{enumerate}
		\item The restriction $[O, \mathfrak{B}]: [\mathfrak{A}, \mathfrak{B}] \rightarrow [\mathfrak{A}_{0}, \mathfrak{B}]$ has a left adjoint in the $2$-category $\mathbf{Gray}$-$\mathbf{Cat}_{2}$.
		\item The adjunction in part (1) is moreover monadic in $\mathbf{Gray}\text{-}\mathbf{Cat}_{2}$, and the $\mathbf{Gray}$-category of pseudoalgebras is given by $\mathbf{TRICAT}_\text{ls}\left(\mathfrak{A}, \mathfrak{B}\right)$. 
		\item Suppose $\mathfrak{B}$ moreover has $\mathbf{Gray}$-enriched codescent objects of all $[n] \mapsto \left({T}^{n}X, {\mu}_{{T}^{n-1}X}\right)$. Let $i: [\mathfrak{A}, \mathfrak{B}] \rightarrow \mathbf{TRICAT}_{\text{ls}}\left(\mathfrak{A}, \mathfrak{B}\right)$ denote the inclusion of strict $\Gr$-natural transformations and $\Gr$-modifications. Then there is an adjunction $L \dashv i$ in $\mathbf{Gray}$-$\mathbf{CAT}_{2}$, with $L$ given on a pseudoalgebra $X$ by taking the $\mathbf{Gray}$-enriched codescent object of $[n] \mapsto \left({T}^{n}X, {\mu}_{{T}^{n-1}X}\right)$.
	\end{enumerate} 
\end{corollary}

\noindent Now suppose $W: \mathfrak{A}^\text{op} \rightsquigarrow \mathbf{Gray}$ is a locally strict trihomomorphism with underlying family of $2$-categories $W_{0}: \mathfrak{A}_{0}^\text{op} \rightarrow \mathbf{Gray}$. Let $\underline{W}_\text{ps}: \Delta_{G}^\text{op} \rightarrow [\mathfrak{A}^\text{op}, \mathbf{Gray}]$ be the $\mathbf{Gray}$-categorical codescent data corresponding to the $T$-pseudoalgebra $W$. We fix notation for data in the image of $\underline{W}_\text{ps}:\Delta_{G}^\text{op} \rightarrow [\mathfrak{A}^\text{op}, \mathbf{Gray}]$ including the $\mathbf{Gray}$-natural transformation $\phi: TW_{0} \rightarrow W_{0}$, and the $\mathbf{Gray}$-modifications $m: \phi.T\phi \Rightarrow \phi.\mu_{W_{0}}$ and $u: 1_{W_{0}} \Rightarrow \phi.\eta_{W_{0}}$. There are other $\mathbf{Gray}$-natural transformations, $\mathbf{Gray}$-modifications and also perturbations in the image of $\underline{W}_\text{ps}: \Delta_{G}^\text{op} \rightarrow [\mathfrak{A}^\text{op}, \mathbf{Gray}]$ which we have not listed but which are listed in the general setting in Section 11.1 of \cite{Gurski Coherence in Three Dimensional Category Theory}.
\\
\\
\noindent We give an explicit description for the reflection $L\left(W\right)$ into $[\mathfrak{A}^\text{op}, \mathbf{Gray}]$ of Lemma \ref{Tricat reflection into GrayHom} part (3). This means giving a presentation in terms of generators and relations of $2$-categories $L\left(W\right)\left(\mathbf{A}\right)$, where $L\left(W\right)$ is the $\mathbf{Gray}$-enriched codescent object of $\underline{W}_\text{ps}: \Delta_\text{ps}^{G} \rightarrow [\mathfrak{A}^\text{op}, \mathbf{Gray}]$. This colimit can be described pointwise in $\mathbf{Gray}$, and hence a presentation can be read off of the universal cocone defining $\mathbf{Gray}$-enriched codescent objects as per Definition 11.7 of \cite{Gurski Coherence in Three Dimensional Category Theory}, which builds upon their Definition 11.6. This cocone involves data as listed below, subject to axioms. 
\begin{itemize}
	\item A $\mathbf{Gray}$-natural transformation $w\in [\mathfrak{A}^\text{op}, \mathbf{Gray}] \left(T\left(W_{0}\right), LW\right)$,
	\item An adjoint equivalence $\zeta \dashv \zeta^{*}$ whose left adjoint is given by a $\mathbf{Gray}$-modification $\zeta: w.T\phi \rightarrow w.\mu_{W_{0}}$,
	\item Two invertible perturbations as depicted below. Here the equalities hold due to $\left(T, \eta, \mu\right)$ being a monad in the $2$-category $\mathbf{Gray}\text{-}\mathbf{CAT}_{2}$.
\end{itemize}

$$\begin{tikzcd}[font=\fontsize{9}{6}]
	&w.T\phi.T\eta_{W_{0}}
	\arrow[rr, "\zeta.1_{T\eta_{W_{0}}}"]
	&{}\arrow[dd, Rightarrow, shorten = 15, "U"]
	&w.\mu_{W_{0}}.T\eta_{W_{0}}
	\arrow[rdd, equal]
	\\
	\\
	w.1_{TW_{0}}
	\arrow[ruu, "Tu"]
	\arrow[rrrr, "1_{w}"']
	&&{}
	&&w.1_{TW_{0}}
\end{tikzcd}
$$

$$\begin{tikzcd}[font=\fontsize{9}{6}, column sep = 30]
	&w.T\phi.\mu_{TW_{0}}
	\arrow[dd, Rightarrow, shorten = 15, shift left = 15, "M"]
	\arrow[r, equal]
	&w.\mu_{W_{0}}.T^{2}\phi
	\arrow[rd, "\zeta.1_{T^{2}\phi}"]
	\\
	w.\mu_{W_{0}}.\mu_{TW_{0}}
	\arrow[ru, "\zeta.1_{\mu_{TW_{0}}}"]
	\arrow[rd, equal]
	&&&w.T\phi.T^{2}\phi
	\\
	&w.\mu_{W}.\mu_{TW_{0}}
	\arrow[r, "\zeta.1_{T\mu_{W_{0}}}"']
	&w.\mu_{W_{0}}.T\mu_{W_{0}}
	\arrow[ru, "1_{w}.Tm"']
\end{tikzcd}$$

\noindent We refer the reader to Section 11.3 of \cite{Gurski Coherence in Three Dimensional Category Theory} for the conditions that the data described above must satisfy. However, note that those conditions simplify considerably due to various coherences in a $\mathbf{Gray}$-monad being identities. These simplifications are similar to the description of the axioms for pseudoalgebras of $\mathbf{Gray}$-monads, as given in Definitions 13.4 and 13.8 of \cite{Gurski Coherence in Three Dimensional Category Theory}. 
\\
\\
\noindent On an object $\mathbf{A} \in \mathfrak{A}$, the $2$-category $L\left(W\right)\left(\mathbf{A}\right)$ has objects given by pairs $\left(\mathbf{f}, B\right)$ where $\mathbf{f} \in \mathfrak{A}\left(\mathbf{A}, \mathbf{B}\right)$ and $B \in W\left(\mathbf{B}\right)$. There are four types of generating morphisms in $L\left(W\right)\left(\mathbf{A}\right)$, two from $T\left(W_{0}\right)$ and two given as the components of $\zeta$ and $\zeta^{*}$. We list these below.

\begin{enumerate}
	\item For any $\boldsymbol\Phi: \mathbf{f} \Rightarrow \mathbf{f}'$ in $\mathfrak{A}$ and any $B \in W\left(\mathbf{B}\right)$, there is a generating morphism $\left(\boldsymbol\Phi, B\right): \left(\mathbf{f}, B\right) \rightarrow \left(\mathbf{f}', B\right)$ in the $2$-category $L\left(W\right)\left(\mathbf{A}\right)$,
	\item For any $\mathbf{f}\in \mathfrak{A}\left(\mathbf{A}, \mathbf{B}\right)$ and any $b: B\rightarrow B'$ in $W\left(\mathbf{B}\right)$ there is a generating morphism $\left(\mathbf{f}, b\right): \left(\mathbf{f}, B\right) \rightarrow \left(\mathbf{f}, B'\right)$ in the $2$-category $L\left(W\right)\left(\mathbf{A}\right)$.
	\item For any triple $\left(C,\mathbf{g}, \mathbf{f}\right)$ consisting of a composable pair \begin{tikzcd}
		\mathbf{A} \arrow[r, "\mathbf{f}"] & \mathbf{B} \arrow[r, "\mathbf{g}"] & \mathbf{C}
	\end{tikzcd} in $\mathfrak{A}$ and $C \in W\left(\mathbf{C}\right)$ there is a generating morphism as depicted below.
	
	$$\begin{tikzcd}
		W\left(\mathbf{f}\right)W\left(\mathbf{g}\right)\left(C\right)
		\arrow[rr, "\zeta_{C}^{\mathbf{g}{,}\mathbf{f}}"] &&	W\left(\mathbf{g}\mathbf{f}\right)\left(C\right)
	\end{tikzcd}$$
	\item For data as in (3), there is also a generating morphism ${\left(\zeta_{C}^{\mathbf{g}{,}\mathbf{f}}\right)}^{*}:   W\left(\mathbf{gf}\right)\left(C\right)\rightarrow W\left(\mathbf{f}\right)W\left(\mathbf{g}\right)\left(C\right)$ in the $2$-category $L\left(W\right)\left(\mathbf{A}\right)$.
\end{enumerate}
\noindent The underlying category of $L\left(W\right)\left(\mathbf{A}\right)$ is generated by these morphisms subject to the composition structure in the underlying category of the $2$-category $W\left(\mathbf{A}\right)$. The different types of $2$-cells in $L\left(W\right)\left(\mathbf{A}\right)$ are listed below.
\begin{enumerate}
	\item For each such triple $\left(C, \mathbf{g}, \mathbf{f}\right)$ formal units and counits for $\zeta \dashv \zeta^{*}$,
	\item There are generators from the $\mathbf{Gray}$-tensor products $\sum\limits_{\mathbf{B}\in \mathfrak{A}_{0}}\mathfrak{A}\left(\mathbf{A}, \mathbf{B}\right) \otimes W\left(\mathbf{B}\right)$. These include $2$-cells from each factor as well as formal interchangers. The formal interchangers are depicted below top right.
	\item There are $2$-cell components for the pseudonatural transformations $\zeta$ and $\zeta^{*}$, on generating morphisms in the $\mathbf{Gray}$-tensor product $$\sum\limits_{\mathbf{B}\in \mathfrak{A}_{0}}\sum\limits_{\mathbf{C}\in \mathfrak{A}_{0}} \mathfrak{A}\left(\mathbf{B}, \mathbf{C}\right)\otimes\mathfrak{A}\left(\mathbf{A}, \mathbf{B}\right)\otimes W\left(\mathbf{C}\right)$$
	
	\begin{tikzcd}
		W\left(\mathbf{f}\right)\left(B\right) \arrow[rr, "W\left(\boldsymbol\Phi\right)_{B}"]
		\arrow[dd, "W\left(\mathbf{f}\right)\left(b\right)"'] &{}\arrow[dd, Rightarrow, shorten = 10, "\left(\boldsymbol{\Phi}{,}b\right)"]
		& W\left(\mathbf{f}'\right)\left(B\right)
		\arrow[dd, "W\left(\mathbf{f}'\right)\left(b\right)"]
		\\
		\\
		W\left(\mathbf{f}\right)\left(B'\right) \arrow[rr, "W\left(\Phi\right)_{B'}"'] &{}& W\left(\mathbf{f}'\right)\left(B'\right)&{}
	\end{tikzcd}\begin{tikzcd}[column sep = 12]
		W\left(\mathbf{f}\right)W\left(\mathbf{g}\right)\left(C\right)
		\arrow[dd, "W{\left(\mathbf{f}\right)}W{\left(\mathbf{g}\right)}\left(c\right)"']
		\arrow[rr, "\zeta_{C{,}B}^{\mathbf{g}{,}\mathbf{f}}"] &{}\arrow[dd, Rightarrow, shorten = 10, "\zeta_{c}^{\mathbf{g}{,}\mathbf{f}}"]
		& W\left(\mathbf{g}\mathbf{f}\right)\left(C\right)
		\arrow[dd, "W\left(\mathbf{gf}\right)\left(c\right)"]
		\\
		\\
		W\left(\mathbf{f}\right)W\left(\mathbf{g}\right)\left(C'\right)
		\arrow[rr, "\zeta_{C'}^{\mathbf{g}{,}\mathbf{f}}"']
		&{}&W\left(\mathbf{g}\mathbf{f}\right)\left(C'\right)
	\end{tikzcd}
	\\
	\begin{tikzcd}
		W\left(\mathbf{f}\right)W\left(\mathbf{g}\right)\left(C\right)
		\arrow[dd, "W\left(\boldsymbol{\Phi}_{W\left(\mathbf{g}\right)}\right)"']
		\arrow[rr, "\zeta_{C}^{\mathbf{g}{,}\mathbf{f}}"] &{}\arrow[dd, Rightarrow, shorten = 10, "\zeta_{C}^{\mathbf{g}{,}\boldsymbol{\Phi}}"]
		& W\left(\mathbf{g}\mathbf{f}\right)\left(C\right)
		\arrow[dd, "W{\mathbf{g}\left(\boldsymbol{\Phi}\right)}_{C}"]
		\\
		\\
		W\left(\mathbf{f}'\right)W\left(\mathbf{g}\right)\left(C\right) 
		\arrow[rr, "\zeta_{C}^{\mathbf{g}{,}\mathbf{f}'}"'] &{}& W\left(\mathbf{g}\mathbf{f}'\right)\left(C\right)
	\end{tikzcd}\begin{tikzcd}[column sep = 12]
		W\left(\mathbf{f}\right)W\left(\mathbf{g}\right)\left(C\right)
		\arrow[dd, "W\left(\mathbf{f}\right){W\left(\Psi\right)}_{C}"']
		\arrow[rr, "\zeta_{C}^{\mathbf{g}{,}\mathbf{f}}"] &{}\arrow[dd, Rightarrow, shorten = 10, "\zeta_{C}^{\boldsymbol{\Psi}{,}\mathbf{g}}"]
		& W\left(\mathbf{g}\mathbf{f}\right)\left(C\right)
		\arrow[dd, "W{\left(\boldsymbol{\Psi}\mathbf{f}\right)}_{C}"]&{}
		\\
		\\
		W\left(\mathbf{f}\right)W\left(\mathbf{g}'\right)\left(C\right) 
		\arrow[rr, "\zeta_{C}^{\mathbf{g}'{,}\mathbf{f}}"'] &{}& W\left(\mathbf{g}'\mathbf{f}\right)\left(C\right)
	\end{tikzcd}
	\item For each $\mathbf{f}\in \mathfrak{A}\left( \mathbf{A},  \mathbf{B}\right)$ and $B \in W\left(\mathbf{B}\right)$ there is a generating invertible $2$-cell in the $2$-category $L\left(W\right)\left(\mathbf{A}\right)$ as depicted below left.
	\item For any quadruple $\left(D, \mathbf{h}, \mathbf{g}, \mathbf{f}\right)$ consisting of a composable triple \begin{tikzcd}
		\mathbf{A} \arrow[r, "\mathbf{f}"] & \mathbf{B} \arrow[r, "\mathbf{g}"] & \mathbf{C} \arrow[r, "\mathbf{h}"] & \mathbf{D}
	\end{tikzcd} in $\mathfrak{A}$ and $D \in W\left(\mathbf{D}\right)$ there is a generating invertible $2$-cell as depicted below right.
\end{enumerate} 

$$\begin{tikzcd}[font=\fontsize{9}{6}]
	W\left(\mathbf{f}\right)\left(B\right) \arrow[rr, "W\left(\mathbf{f}\right)\left(u_\mathbf{B}\right)_{B}"]
	\arrow[rrdd, bend right = 30, "1_{W\left(\mathbf{f}\right)\left(B\right)}"']&{}\arrow[dd, shorten = 10, Rightarrow, "U_\mathbf{f}^{B}"]& W\left(\mathbf{f}\right)W\left(1_\mathbf{B}\right)\left(B\right)
	\arrow[dd, "\zeta_{1_\mathbf{B}{,}\mathbf{f}}^{B}"]
	\\
	\\
	&{}&W\left(\mathbf{f}\right)\left(B\right)&{}
\end{tikzcd}\begin{tikzcd}[font=\fontsize{9}{6}]
	W\left(\mathbf{f}\right)W\left(\mathbf{g}\right)W\left(\mathbf{h}\right)\left(D\right)
	\arrow[rr, "\zeta_{\mathbf{g}{,}\mathbf{f}}^{W\left(\mathbf{h}\right)\left(D\right)}"]
	\arrow[dd, "W\left(\mathbf{f}\right)\left(\zeta_{\mathbf{h}{,}\mathbf{g}}^{D}\right)"']
	&{}\arrow[dd, Rightarrow, shorten = 10, shift right = 5, "M_{\mathbf{h}{,}\mathbf{g}{,}\mathbf{f}}^{D}"]
	&
	W\left(\mathbf{g}\mathbf{f}\right)W\left(\mathbf{h}\right)\left(D\right)
	\arrow[dd, "\zeta_{\mathbf{h}{,}\mathbf{gf}}^{D}"]
	\\
	\\
	W\left(\mathbf{f}\right)W\left(\mathbf{hg}\right)\left(D\right)
	\arrow[rr, "\zeta_{\mathbf{hg}{,}\mathbf{f}}^{D}"']
	&{}&
	W\left(\mathbf{hgf}\right)\left(D\right)
\end{tikzcd}$$
\noindent These data are subject to various relations specifying pseudonaturality of the components of $\zeta$, the modification condition for the components of $M$ and $U$, and the axioms for a $\mathbf{Gray}$-enriched codescent object. We do not list these as they are easily inferred from Definition 11.6 of \cite{Gurski Coherence in Three Dimensional Category Theory}.
\\
\\
\noindent Now assume that $W$ is a $\mathbf{Gray}$-functor. The counit $\varepsilon_{W}: Li\left(W\right) \rightarrow W$ is given on objects by $\left(\mathbf{f}, B\right) \mapsto W\left(\mathbf{f}\right)\left(B\right)$ and on generating morphisms as depicted below.

\begin{align*}
	&\left(\boldsymbol{\Phi}{,}B\right) &\mapsto &&W\left(\boldsymbol\Phi\right)_{B}: W\left(\mathbf{f}\right)\left(B\right) \rightarrow W\left(\mathbf{f}'\right)\left(B\right)
	\\
	&\left(\mathbf{f}{,}b\right) & \mapsto&& W\left(\mathbf{f}\right)\left(b\right): W\left(\mathbf{f}\right)\left(B\right) \rightarrow W\left(\mathbf{f}\right)\left(B'\right)
	\\
	&\zeta_{C}^{\mathbf{g}{,}\mathbf{f}} & \mapsto &&\chi_{\mathbf{g}{,}\mathbf{f}}^{C}: W\left(\mathbf{f}\right)W\left(\mathbf{g}\right) \left(C\right) \rightarrow W\left(\mathbf{gf}\right)\left(C\right)
\end{align*}

\noindent These prescriptions extend in the evident way to generating $2$-cells as data from $2$-categories $W\left(\mathbf{B}\right)$ and from the $\mathbf{Gray}$-category $\mathfrak{A}$ vary. Meanwhile, the generating $2$-cells listed above are mapped by $\varepsilon_{W\left(\mathbf{A}\right)}$ to the following data in $W\left(\mathbf{A}\right)$.
\begin{itemize}
	\item $\left(\boldsymbol\Phi, b\right)$ is mapped to the pseudonaturality constraint of ${W\left(\boldsymbol\Phi\right)}: W\left(\mathbf{f}\right) \rightarrow W\left(\mathbf{f}'\right)$ on $b \in W\left(\mathbf{B}\right)\left(B, B'\right)$.
	\item $\left(\zeta_{C}^{\mathbf{g}, \boldsymbol\Phi}\right)$ and $\left(\zeta_{C}^{\boldsymbol\Psi, \mathbf{f}}\right)$ are mapped to the modification components of $\chi_{\mathbf{g}, \boldsymbol{\Phi}}$ and $\chi_{\boldsymbol\Psi{,}\mathbf{f}}$ on $c\in W\left(\mathbf{C}\right)\left(C, C'\right)$.
	\item $\zeta_{c}^{\mathbf{g}, \mathbf{f}}$ is mapped to the pseudonaturality constraint of $\chi_{\mathbf{g}{,}\mathbf{f}}$ on $c\in W\left(\mathbf{C}\right)\left(C, C'\right)$.
	\item $U_{\mathbf{f}}^{B}$ and $M_{\mathbf{h}{,}\mathbf{g}{,}\mathbf{f}}^{D}$ are mapped to the modification component of $\delta_{\mathbf{f}}$ on $B$ and the modification component of $\omega_{\mathbf{h}{,}\mathbf{g}{,}\mathbf{f}}$ on $D$ respectively.
\end{itemize}

\begin{notation}
	For $\mathfrak{A}$ a $\mathbf{Gray}$-category, the full-sub $\mathbf{Gray}$-category of $[\mathfrak{A}, \mathbf{Gray}]$ on $\mathbf{Gray}$-functors which are flexible as algebras over $\left([\mathfrak{A}_{0}, \mathbf{Gray}], T\right)$ will be written as ${[\mathfrak{A}, \mathbf{Gray}]}_\text{flex}$.
\end{notation}

\begin{theorem}\label{trinatural transformation classifier cofibrancy}
	Let $\mathfrak{A}$ be a $\mathbf{Gray}$-category.
	\begin{enumerate}
		\item If $\mathfrak{A}$ has cofibrant hom $2$-categories then $L:\mathbf{TRICAT}_\text{ls}\left(\mathfrak{A}, \mathbf{Gray}\right) \rightarrow  [\mathfrak{A}^\text{op}, \mathbf{Gray}]$ preserves and reflects the property of being pointwise cofibrant.
		\item The counit $\varepsilon_{W}: Li\left(W\right) \rightarrow W$ is a pointwise (hence projective) trivial fibration.
		\item If $W: \mathfrak{A}^\text{op} \rightarrow \mathbf{Gray}$ is a projective cofibrant $\mathbf{Gray}$-functor then there is an equivalence $e: \eta_{iW} \Rightarrow i\left(s\right)$ in $\mathbf{TRICAT}_\text{ls}\left(\mathfrak{A}, \mathbf{Gray}\right)$, with $s: W \rightarrow Li\left(W\right)$ providing a biequivalence section to $\varepsilon_{W}$ in $[\mathfrak{A}^\text{op}, \mathbf{Gray}]$.
		\item For any locally strict trihomomorphism $W: \mathfrak{A}^\text{op}\rightsquigarrow \mathbf{Gray}$, the $\mathbf{Gray}$-functor $L\left(W\right):\mathfrak{A}^\text{op}\rightarrow \mathbf{Gray}$ is projective cofibrant.
		\item A $\mathbf{Gray}$-functor $W: \mathfrak{A}^\text{op} \rightarrow \mathbf{Gray}$ is a flexible algebra for the $\mathbf{Gray}$-monad $\big([\mathfrak{A}_{0}, \mathbf{Gray}], T\big)$ if and only if it is projective cofibrant.
		\item The $\mathbf{Gray}$ adjunction $L \dashv i$ restricts to a triequivalence between $\mathbf{Tricat}_\text{ls}\left(\mathfrak{A}^\text{op}, \mathbf{Gray}\right)$ and projective cofibrant objects (equivalently, ${[\mathfrak{A}^\text{op}, \mathbf{Gray}]}_{\text{flex}}$).
		\item If $W: \mathfrak{A}^\text{op} \rightarrow \mathbf{Gray}$ is projective cofibrant and $V: \mathfrak{A}^\text{op} \rightarrow \mathbf{Gray}$ is any $\mathbf{Gray}$-functor then \begin{enumerate}
			\item Any trinatural transformation $p: iW \rightarrow iV$ is equivalent in the hom-$2$-category displayed below, to a $\mathbf{Gray}$-natural transformation.
			
			$$\mathbf{TRICAT}_\text{ls}\left(\mathfrak{A}^\text{op}, \mathbf{Gray}\right)\left(iW, iV\right)$$
			\item If $q, r:W \rightarrow V$ are $\mathbf{Gray}$-natural transformations then any trimodification $\sigma: q \Rightarrow r$ is isomorphic to a $\mathbf{Gray}$-modification.
		\end{enumerate}
		\item Let $\mathfrak{B}$ be a tricategory. Then $\mathbf{TRICAT}\big(\mathfrak{B}, \mathbf{Bicat}\big)$ is triequivalent to ${[\mathbf{st}_{3}\left(\mathfrak{B}\right), \mathbf{Gray}]}_\text{flex}$.
	\end{enumerate}
\end{theorem}

\begin{proof}
	Part (1) is clear from the presentation of $L\left(W\right)\left(\mathbf{A}\right)$ described above, since the only relations in the $2$-categories $\left(LW\right)\left(\mathbf{A}\right)$ are from $W\left(\mathbf{A}\right)$. For part (2), $i\left(\varepsilon_{W}\right): iLi\left(W\right) \rightarrow iW$ is a retraction to $\eta_{iW}$ by the triangle identity $i\left(\varepsilon_{W}\right).\eta_{iW} = 1_{iW}$ for the $\mathbf{Gray}$-adjunction $L \dashv i$. The components of $\varepsilon_{W}$ are themselves also retractions and are hence surjective on objects and full on $1$-cells and $2$-cells. Moreover by two-out-of-three and Lemma \ref{Tricat reflection into GrayHom} part(5), $i\left(\varepsilon_{W}\right)$ is moreover a biequivalence. This means that its components are also faithful, and are hence trivial fibrations.
	\\
	\\
	\noindent For part (3), projective cofibrancy of $W$ combined with the fact that $\varepsilon_{W}$ is a projective trivial fibration from part (2), implies that there is a $\mathbf{Gray}$-natural transformation $s: W \rightarrow Li\left(W\right)$ satisfying $\varepsilon_{W}.s = 1_{W}$. But $s$ and $\eta_{iW}$ are both sections to $i\left(\varepsilon_{W}\right)$. Therefore, the strict trinatural transformation $i\left(s\right)$ is equivalent to $\eta_{iW}$ in $\mathbf{TRICAT}_{\text{ls}}\left(\mathfrak{A}^\text{op}, \mathbf{Gray}\right)$. By the universal property of the $\mathbf{Gray}$ adjunction $L \dashv i$, this equivalence in $\mathbf{TRICAT}_\text{ls}\left(\mathfrak{A}^\text{op}, \mathbf{Gray}\right)$ determines an equivalence $\overline{e}: s.\varepsilon_{W} \rightarrow 1_{Li\left(W\right)}$ in $[\mathfrak{A}^\text{op}, \mathbf{Gray}]$. The rest of the biequivalence structure for $s \dashv \varepsilon_{W}$ follows similarly, using the $\mathbf{Gray}$-adjunction $L \dashv i$.
	\\
	\\
	\noindent The argument for part (4) is entirely analogous to the proof of Theorem 5.12 of \cite{Lack Homotopy Theoretic Aspects of 2-monads}. We summarise this argument for completeness. Let $t: E \rightarrow B$ be a trivial fibration in $[\mathfrak{A}^\text{op}, \mathbf{Gray}]$ and let $v: L\left(W\right) \rightarrow B$ be $\mathbf{Gray}$-natural. For cofibrancy, a lifting $L\left(W\right) \rightarrow E$ is required. But since trivial fibrations are in particular biequivalences in $\mathbf{TRICAT}_\text{ls}\left(\mathfrak{A}^\text{op}, \mathbf{Gray}\right)$, there is a trinatural transformation $s: B \rightarrow E$ providing a section to $t$. Then by the universal property of $L \dashv i$, there is a unique $\mathbf{Gray}$-natural transformation $u: L\left(W\right) \rightarrow E$ satisfying $u\eta_{W} = su\eta_{W}$. Composing with $t$, we have that $tu\eta_{W} = tsv\eta_{W}$, which is in turn equal to $v\eta_{W}$ since $ts=1_{B}$. Finally, the universal property of $\eta_{W}$ means that $tu=v$, as desired. These data in $\mathbf{TRICAT}_\text{ls}\left(\mathfrak{A}^\text{op}, \mathbf{Gray}\right)$ are depicted below, with trinatural transformations that need not be $\mathbf{Gray}$-natural depicted in red.
	
	$$\begin{tikzcd}[font=\fontsize{9}{6}]
		&&&&E\arrow[dd, two heads, "t"]
		\\
		\\
		W \arrow[rr, red, "\eta_{W}"']
		&& L\left(W\right)
		\arrow[rruu, bend left = 30, dashed, "u"]
		\arrow[rr, "v"']
		&& B \arrow[rr, equal]
		&&B\arrow[lluu, red, "s"']
	\end{tikzcd}$$
	
	\noindent For part (5), projective cofibrant $\implies$ flexible algebra follows from part (3), while flexible algebra $\implies$ projective cofibrant follows from part (4) since retracts of cofibrant objects are cofibrant, and flexibility means that $W$ is a retract of $Li\left(W\right)$.
	\\
	\\
	\noindent We now consider part (6). Since $\eta: 1 \Rightarrow iL$ is a biequivalence from Lemma \ref{Tricat reflection into GrayHom} part (5), we just need to characterise the objects $W \in [\mathfrak{A}^\text{op}, \mathbf{Gray}]$ for which $\varepsilon_{W}: Li\left(W\right) \rightarrow W$ is a biequivalence internal to $[\mathfrak{A}^\text{op}, \mathbf{Gray}]$. By part (3), projective cofibrancy is sufficient. The converse follows by a similar argument to the one given for part (4). Given a trivial fibration $t: E \rightarrow B$ with section $p$ in $\mathbf{TRICAT}_\text{ls}\left(\mathfrak{A}^\text{op}, \mathbf{Gray}\right)$ and a $\mathbf{Gray}$-natural $v: W \rightarrow B$, consider the map $u: Li\left(W\right) \rightarrow E$ induced by the trinatural transformation \begin{tikzcd}
		W \arrow[r, "v"] & B \arrow[r, "p"] & E
	\end{tikzcd}. This will satisfy $tu = tpv\varepsilon_{W}$, and hence $tu = v\varepsilon_{W}$ since $p$ is a section to $t$. But now if $\varepsilon_{W}$ also has a section $s: W \rightarrow Li\left(W\right)$ in $[\mathfrak{A}^\text{op}, \mathbf{Gray}]$ then we also have $tus = v\varepsilon_{W}s= v$. Thus $us$ is the desired lifting of $v$ along $t$, and hence $W$ is projective cofibrant.
	\\
	\\
	\noindent Part (7) follows from part (4) using the universal property of the unit for $L \dashv i$. This is similar to the proof of $(d) \implies (e)$ in Theorem 4.14 of \cite{Two Dimensional Monad Theory}, but once again we give details for completeness. By part (4), projective cofibrancy means that there is a $\mathbf{Gray}$-natural $s: W \rightarrow Li\left(W\right)$ and an equivalence $e: i\left(s\right) \rightarrow \eta_{iW}$ in $\mathbf{TRICAT}_\text{ls}\left(\mathfrak{A}^\text{op}, \mathbf{Gray}\right)$. The weak $\left(3, k\right)$-transfors in parts (a) and (b) factor through $\eta$ via strict $\left(3, k\right)$-transfors, by the universal property. Composing with $e$ gives the desired equivalence (resp. isomorphism).
	\\
	\\
	\noindent Finally, the triequivalence for part (8) is constructed as the composite below. Here the last step uses part (1) and the fact that $\mathbf{st}_{3}\left(\mathfrak{B}\right)$ is hom-wise cofibrant. 
	
	{\small
		\begin{align*}
			\mathbf{TRICAT}\left(\mathfrak{B}{,}\mathbf{Bicat}\right) 
			&\sim \mathbf{TRICAT}\left(\mathfrak{B}{,}\mathbf{Gray}_{\text{c}}\right)
			&\text{Extending along } \mathbf{st}_{2}
			\\
			&&\text{ of Proposition} \ref{Hom triequivalent to Gray c} \text{ part (3)}
			\\
			&\sim \mathbf{TRICAT}\left(\mathbf{st}_{3}\left(\mathfrak{B}\right){,}\mathbf{Gray}_{\text{c}}\right)
			&\text{Restricting along the triequivalence inverse}
			\\
			&&E_{\mathfrak{B}} \text{ of Proposition } \ref{Proposition 3D strictification tetraadjunction} \text{ part (1)}
			\\
			&\sim \mathbf{TRICAT}_{\text{ls}}\left(\mathbf{st}_{3}\left(\mathfrak{B}\right){,}\mathbf{Gray}_{\text{c}}\right)
			&\text{Cofibrancy of } \mathbf{st}_{3}(\mathfrak{B})
			\\
			&\sim
			{[}\mathbf{st}_{3}\left(\mathfrak{B}\right){,}\mathbf{Gray}{]}_\text{flex}
			& \text{ Proposition }\ref{Proposition 3D strictification tetraadjunction}\text{ part (6)}
		\end{align*} 
	}
\end{proof}

\begin{remark}\label{generalities of flexibility and cofibrancy}
	Theorem \ref{trinatural transformation classifier cofibrancy} shows that projective cofibrancy and flexibility are equivalent properties, and that restricting to such weights suffices to capture $\mathbf{TRICAT}_\text{ls}\left(\mathfrak{A}^\text{op}, \mathbf{Gray}\right)$ up to triequivalence. Such a relationship is familiar in the $2$-categorical settings of \cite{Lack Homotopy Theoretic Aspects of 2-monads} and \cite{Two Dimensional Monad Theory}. It may also hold in even greater generality, such as a $\mathbf{Gray}$-enriched version of the settings in \cite{Lack Homotopy Theoretic Aspects of 2-monads} and \cite{Two Dimensional Monad Theory}, or perhaps even over sufficiently well-behaved monoidal model categories $\mathcal{V}$. We leave the investigations of such generalities to future research.
\end{remark}

\begin{corollary}\label{Reduction to Gray natural biequivalence}
	Let $\mathfrak{B}$ be a $\mathbf{Gray}$-category.
	
	\begin{enumerate}
		\item Suppose $\mathfrak{B}$ is hom-wise cofibrant. Consider $W$ and $F$ as in Theorem \ref{tricolimit reduction to strict 3 k transfors via LW} part (1). Then $IW\odot F$ is a tricolimit if and only if there is a projective weak equivalence $w: \mathfrak{B}\left(IW\odot F, ? \right) \rightarrow [\mathfrak{A}^\text{op}, \mathbf{Gray}]\left(LW, \mathfrak{B}\left(F-, ?\right)\right)$ in $[\mathfrak{B}, \mathbf{Gray}]$, such that the biequivalence $iw \dashv {\left(iw\right)}^{*}$ in the $\Gr$-category $\mathbf{Tricat}_\text{ls}\left(\mathfrak{B}, \mathbf{Gray}\right)$ has a strict counit.
		\item Suppose $\mathfrak{B}$ is hom-wise cofibrant. If the $\mathbf{Gray}$-functor depicted below is projective cofibrant, then the projective weak equivalence $w$ in part (1) is moreover a biequivalence internal to $[\mathfrak{B}, \mathbf{Gray}]$.
		
		$$\begin{tikzcd}
			\mathfrak{B} 
			\arrow[rr, "\mathfrak{Y}"] && {[}\mathfrak{B}^\text{op}{,}\mathbf{Gray}{]} \arrow[rr, "{F}^\text{op}"]
			&& {[}\mathfrak{A}^\text{op}{,}\mathbf{Gray}{]}
			\arrow[rrrr, "{[}\mathfrak{A}^\text{op}{,}\mathbf{Gray}{]}\left(W{,}-\right)"]
			&&&& \mathbf{Gray}
		\end{tikzcd}$$
		
		\item Let $W: \mathfrak{A}^\text{op} \rightarrow \mathbf{Gray}$ be projective cofibrant on a hom-wise cofibrant $\mathfrak{A}$ and consider a $\mathbf{Gray}$-functor $F: \mathfrak{A} \rightarrow \mathfrak{B}$. If the $\mathbf{Gray}$-enriched colimit $W\cdot F$ exists then it is also a tricolimit.
		\item If $\mathfrak{B}$ has all flexibly weighted $\mathbf{Gray}$-enriched limits (resp. colimits) on hom-wise cofibrant $\mathbf{Gray}$-categories then it has all trilimits (resp. tricolimits).
		\item Suppose $\mathfrak{B}$ and $\mathfrak{C}$ are $\mathbf{Gray}$-categories with small flexibly weighted $\mathbf{Gray}$-enriched limits (resp. colimits), and let $F: \mathfrak{A} \rightarrow \mathfrak{B}$ be a $\mathbf{Gray}$-functor that preserves flexibly weighted $\mathbf{Gray}$-enriched limits (resp. colimits). Then $F$ also preserves trilimits (resp. tricolimits).
	\end{enumerate}
\end{corollary}

\begin{proof}
	For part (1), by Theorem \ref{tricolimit reduction to strict 3 k transfors via LW} part (1) there is already a biequivalence internal to the $\mathbf{Gray}$-category $\mathbf{TRICAT}_\text{ls}\left(\mathfrak{A}^\text{op}, \mathbf{Gray}\right)$. We need to show that $w$ is $\mathbf{Gray}$-natural and that the counit of $w \dashv w^{*}$ is a $\mathbf{Gray}$-modification. But these follow from Theorem \ref{trinatural transformation classifier cofibrancy} part (7 a) and (7 b) respectively. This is because $\mathfrak{B}\left(IW\cdot F, ?\right)$ is representable, and is in particular projective cofibrant. Similarly, part (2) uses projective cofibrancy to apply Theorem \ref{tricolimit reduction to strict 3 k transfors via LW} part (7 a) and (7 b) to replace the unit by a $\mathbf{Gray}$-modification.
	\\
	\\
	\noindent For part (3), first observe that since $\mathfrak{A}$ is hom-wise cofibrant and $W$ is projective cofibrant, $W$ is also pointwise cofibrant by Proposition \ref{projective and homwise implies pointwise}. Hence by Lemma \ref{strictifying components of a trinatural transformation into Hom}, trinatural transformations out of $IW$ factor through $I: \mathbf{Gray} \rightarrow \mathbf{Bicat}$. Thus similar arguments using Theorem \ref{trinatural transformation classifier cofibrancy} parts (6) and (7) apply. Part (4) follows from Theorem \ref{tricolimit reduction to strict 3 k transfors via LW} parts (2) and (3) using the fact that $\mathbf{Gray}$-functors of the form $L\left(W\right)$ are flexible, from Theorem \ref{trinatural transformation classifier cofibrancy} parts (4) and (5). Part (5) follows similarly.
\end{proof}

\begin{remark}\label{Remark enriched weakness}
	The results thus far have effectively reduced the study of tricategorical universal properties to the setting of enriched weakness considered in \cite{Enriched Weakness} where $\mathcal{V}$ is the (category of fibrant objects in the) monoidal model category $\mathbf{Gray}$, and where the class of weak maps $\mathcal{E}$ are taken to be the biequivalences, which are the weak equivalences in this model structure. Note that although this general setting is not treated in \cite{Enriched Weakness}, in the specific case considered here all objects are fibrant, and indeed fibrant objects are closed under limits, colimits, and tensor product. Our setting is also considered in Example 3.10 of \cite{Bourke Lack Vokrinek Adjoint Functor Theorems for Homotopically Enriched Categories} where adjoint functor theorems expressed in terms of biequivalence retractions are developed.
\end{remark}

\section{Examples}\label{Examples trilimits and tricolimits}

\noindent We are now in a position to show that our main examples of tricategories, and certain functor tricategories, have tricategorical limits and colimits. This is done in Proposition \ref{main examples have all trilimits and tricolimits}, and an explicit formula for tricategorical limits in $\mathbf{Bicat}$ is given in Proposition \ref{trilimits in Hom}. This is applied to analyse how representables and the Yoneda embedding interact with trilimits in Proposition \ref{representables preserve trilimits}, and provide a coherence result for tricategories with finite trilimits in Theorem \ref{coherence for tricategories with trilimits}. Then in Subsection \ref{subsection specific weights and diagrams} we will re-examine various constructions of interest in the literature and show that they are examples of tricategorical limits and colimits. These include strictification of bicategories, which will be seen as a tricategorical copower in Example \ref{strictification as a tricolimit}, and a variant of the centre construction of \cite{Crans Generalised Centers} which will be seen as a descent object. Finally, Theorem \ref{biessential surjectivity on objects characterises trikleisli pseudoadjunctions in Gray or Hom} will characterise trikleisli pseudoadjunctions in $\mathbf{Bicat}$ and $\mathbf{Gray}$ as being ones whose left pseudoadjoints are biessentially surjective on objects.

\subsection{Trilimits in $\mathbf{Bicat}$, and in functor tricategories}\label{subsection trilimits in hom and functor tricategories}

\begin{proposition}\label{main examples have all trilimits and tricolimits}
	\hspace{1mm}
	\begin{enumerate}
		\item 
		The tricategories $\mathbf{Gray}$, $\mathbf{Gray}_{c}$, and $\mathbf{Bicat}$ all have tricategorical limits and colimits.
		\item If $\mathfrak{B}$ is a $\mathbf{Gray}$-category with flexible colimits and $\mathfrak{A}$ is another $\mathbf{Gray}$-category then $\mathbf{Tricat}_\text{ls}\left(\mathfrak{A}, \mathfrak{B}\right)$ has tricategorical colimits.
		\item If $\mathfrak{A}$ and $\mathfrak{B}$ are as in part (2) but now $\mathfrak{B}$ has flexible limits instead of flexible colimits then  $\mathbf{Tricat}_\text{ls}\left(\mathfrak{A}, \mathfrak{B}\right)$ has tricategorical limits.
		\item If $\mathfrak{A}$ is a tricategory then the tricategory $\mathbf{TRICAT}\left(\mathfrak{A}, \mathbf{Bicat}\right)$ has all tricategorical limits and colimits.
	\end{enumerate}
\end{proposition}

\begin{proof}
	Since $\mathbf{Gray}$ is cocomplete as a category enriched over itself, the existence of tricategorical limits and colimits follows from Theorem \ref{tricolimit reduction to strict 3 k transfors via LW} parts (2) and (3). For $\mathbf{Gray}_{c}$, the existence of tricategorical limits and colimits also follows from Theorem \ref{tricolimit reduction to strict 3 k transfors via LW} part (3) and Theorem \ref{trinatural transformation classifier cofibrancy} part (4), since cofibrant objects are closed under projective cofibrant weights on hom-wise cofibrant $\mathbf{Gray}$-categories. For $\mathbf{Bicat}$, this follows from the fact for $\mathbf{Gray}_{c}$, since by Lemma \ref{left triadjoint preserves tricolimits} the triequivalence between them will preserve tricategorical limits and colimits. This proves part (1).
	\\
	\\
	\noindent For part (2), first observe that by Corollary \ref{Gray monad for functors on codomain with copowers by homs of domain}, the restriction $O: [\mathfrak{A}, \mathfrak{B}] \rightarrow [\mathfrak{A}_{0}, \mathfrak{B}]$ is $\mathbf{Gray}$-monadic and $\mathbf{PsAlg}\left([\mathfrak{A}_{0}^\text{op}, \mathfrak{B}], T\right)\cong \mathbf{Tricat}_\text{ls}\left(\mathfrak{A}, \mathfrak{B}\right)$. The result then follows from Theorem 7.5 of \cite{Three dimensional monad theory}. Part (3) follows from part (2) by taking $\mathfrak{A}^\text{op}$ and $\mathfrak{B}^\text{op}$ in place of $\mathfrak{A}$ and $\mathfrak{B}$.
	\\
	\\
	\noindent For part (4), $\mathbf{TRICAT}\left(\mathfrak{A}, \mathbf{Bicat}\right)$ is triequivalent to $\mathbf{TRICAT}_\text{ls}\left(\mathbf{st}_{3}\left(\mathfrak{A}\right), \mathbf{Gray}_{c}\right)$ by restricting along $E^{*}: \mathbf{st}_{3}\left(\mathfrak{A}\right) \rightsquigarrow \mathfrak{A}$, extending along $\mathbf{st}_{2}': \mathbf{Bicat} \rightsquigarrow \mathbf{Gray}_{c}$, and using cofibrancy of $\mathbf{st}_{3}\left(\mathfrak{A}\right)$ as per Proposition \ref{Proposition cofibrancy localises} to replace general trihomomorphisms with locally strict ones. The result then follows from parts (2) and (3), since $\mathbf{Gray}_{c}$ has all flexibly weighted limits and colimits.
\end{proof}

\begin{lemma}\label{evaluation from 1 is biequivalence}
	For any bicategory $\mathcal{B}$, evaluation defines a biequivalence $\mathbf{Bicat}\left(\mathbf{1}, \mathcal{B}\right) \rightarrow \mathcal{B}$, which is moreover trinatural in $\mathcal{B}$.
\end{lemma}

\begin{proof}
	Observe that evaluation is surjective on objects and arrows, with the constant pseudofunctor at $X \in \mathcal{B}$ being sent to $X$ and the pseudonatural transformation with component $f: X \rightarrow Y$ being sent to $Y$. Finally, observe that modifications between such pseudonatural transformations are precisely $2$-cells of $\mathcal{B}$, so evaluation is also locally fully faithful. Trinaturality is also an easy inspection. See also (1.34) of \cite{Street Fibrations in bicategories}.
\end{proof}

\noindent Proposition \ref{trilimits in Hom}, to follow, categorifies Proposition 1.15 of \cite{Street Fibrations in bicategories}. It is mentioned in Remark 3.1 of \cite{Buhne PhD}.

\begin{notation}
	For a tricategory $\mathfrak{B}$, we will denote the tricategory $\mathbf{TRICAT}\left(\mathfrak{B}^\text{op}, \mathbf{Bicat}\right)$ by $\widehat{\mathfrak{B}}$.
\end{notation}

\begin{proposition}\label{trilimits in Hom}
	Consider a general pair of trihomomorphisms $W: \mathfrak{A} \rightsquigarrow \mathbf{Bicat}$ and \\$F:\mathfrak{A} \rightsquigarrow \mathbf{Bicat}$. A trilimit can be given as depicted below.
	
	$$<W, F> = \widehat{\mathfrak{A}^\text{op}}\left(W, F\right)$$
	
\end{proposition}

\begin{proof}
	By part (1) of Proposition \ref{main examples have all trilimits and tricolimits}, the trilimit exists. Let $<W, F>$ denote the trilimit. The claim follows from the following chain of trinatural biequivalences.
	
	\begin{align*}
		<W{,}F>
		&\sim\mathbf{Bicat}\left(1{,}<W{,}F>\right)
		&\text{Lemma \ref{evaluation from 1 is biequivalence}} 
		\\
		&\sim\widehat{\mathfrak{A}^\text{op}}\left(W{,}\mathbf{Bicat}\left(1{,}<W{,}F>\right) \right)
		& \text{Definition \ref{trilimits and tricolimits definition}}
		\\
		&\sim \widehat{\mathfrak{A}^\text{op}}\left(W{,}F \right)
		&\text{Lemma \ref{evaluation from 1 is biequivalence}} 
	\end{align*}

\end{proof}

\begin{proposition}\label{representables preserve trilimits}
	Let $W: \mathfrak{A} \rightsquigarrow \mathbf{Bicat}$ and $F: \mathfrak{A} \rightsquigarrow \mathfrak{B}$ be trihomomorphisms. \begin{enumerate}
		\item Consider an object $X \in \mathfrak{B}$. Then the representable $\mathfrak{B}\left(X, -\right): \mathfrak{B} \rightsquigarrow \mathbf{Bicat}$ preserves the trilimit $<W, F>$ if it exists.
		\item If there is an object $Y \in \mathfrak{B}$ so that for all $X \in \mathfrak{B}$, the bicategory $\mathfrak{B}\left(X, Y\right)$ is a trilimit of $\mathfrak{B}\left(X, F-\right)$ weighted by $W$, then $Y$ is a trilimit of $F$ weighted by $W$.
		\item The Yoneda embedding $\mathfrak{Y}_{\mathfrak{B}}:\mathfrak{B} \rightsquigarrow \widehat{\mathfrak{B}}$ preserves trilimits. 
	\end{enumerate}
\end{proposition}

\begin{proof}
	For a bicategory $\mathcal{C}$, consider the following chain of biequivalences.
	{\small
		\begin{align*}
			\mathbf{Bicat}\big(\mathcal{C}{,}\mathfrak{B}\left(X{,}<W{,}F>\right)\big) 
			&\sim \mathbf{Bicat}\big(\mathcal{C}{,}\widehat{\mathfrak{A}^\text{op}}\left(W{,}\mathfrak{B}\left(X{,}F-\right)\right)\big)
			& \text{Definition of }<W{,}F>
			\\
			&\sim \mathbf{Bicat}\big(\mathcal{C}{,}<W{,}\mathfrak{B}\left(X{,}F-\right)>\big)
			& \text{Proposition \ref{trilimits in Hom}}
			\\
			& \sim \widehat{\mathfrak{A}^\text{op}}\big(W{,}\mathbf{Bicat}\left(\mathcal{C}{,}\mathfrak{B}\left(X{,}F-\right)\right)\big)
			&\text{Definition of }
			\\
			&& <W{,}\mathfrak{B}\left(X{,}F-\right)>
		\end{align*}
		\small}%
	\\
	\noindent The proof for part (1) follows by noting that each biequivalence is trinatural in $\mathcal{C}\in \mathbf{Bicat}$ and $X \in \mathfrak{B}$, while the proof for part (2) follows by taking the particular case $\mathcal{C}:= \mathbf{1}$ and using Lemma \ref{evaluation from 1 is biequivalence}. For part (3), the proof follows from the following chain of biequivalences which are trinatural in $G \in \widehat{\mathfrak{B}}$.
	
	{\small
		\begin{align*}
			\widehat{\mathfrak{B}}\left(G{,}\mathfrak{B}\left(-{,}<W{,}F>\right)\right)
			&\sim \widehat{\mathfrak{B}}\left(G{,}\widehat{\mathfrak{A}^\text{op}}\left(W{,}\mathfrak{B}\left(-{,}F\right)\right)\right)
			&\text{Definition of } <W{,}F>
			\\
			&\sim
			\widehat{\mathfrak{A}^\text{op}}\left(W{,}\widehat{\mathfrak{B}}\left(G{,}\mathfrak{B}\left(-{,}F\right)\right)\right)
			& \widehat{\mathfrak{A}^\text{op}}\left(G{,}-\right)\text{ preserves }
			\\
			&&<W{,}F>\text{ by part (1)}
		\end{align*}
		\small}
\end{proof}

\noindent Theorem \ref{coherence for tricategories with trilimits}, to follow, categorifies Theorem 4.1 of \cite{Power Coherence for categories with finite bilimits}. It refers to the definitions of finite trilimits and finite enriched limits given before Theorem 3.3.4 of \cite{Campbell PhD}, which we recall below.

\begin{definition}\label{finite weights}
	Let $W: \mathfrak{A} \rightsquigarrow \mathbf{Bicat}$ be a trihomomorphism, considered as a weight for trilimits, and let $W': \mathfrak{A}' \rightarrow \mathbf{Gray}$ be a $\Gr$-functor.
	\begin{enumerate}
		\item $W$ is called a \emph{finite weight} if the tricategory $\mathfrak{A}$ has finitely many $3$-cells, and for each $X \in \mathfrak{A}$ the bicategory $WX$ has finitely many $2$-cells.
		\item $W'$ is called a \emph{finite weight} if $\mathfrak{A}'$ is a finitely presented $\mathbf{Gray}$-category, and for each $X \in \mathfrak{A}'$ the $2$-category $W'X$ is finitely presented.
	\end{enumerate} 
\end{definition}

\noindent Our proof is similar in spirit to the one given in \cite{Power Coherence for categories with finite bilimits}, except for two crucial points. Firstly, we need to replace $\mathbf{TRICAT}\left(\mathfrak{B}^\text{op}, \mathbf{Bicat}\right)$ with a triequivalent $\mathbf{Gray}$-category in which finite flexible $\mathbf{Gray}$-enriched limits can be considered. Secondly, a cofibrancy assumption is also required on the hom-$2$-categories of the shapes for these enriched limits.

\begin{theorem}\label{coherence for tricategories with trilimits}
	Let $\mathfrak{B}$ be a small tricategory with finite trilimits. Then there is a triequivalence $G: \mathfrak{B} \rightsquigarrow \overline{\mathfrak{B}}$ with $\overline{\mathfrak{B}}$ a (possibly large) $\mathbf{Gray}$-category that has $\mathbf{Gray}$-enriched limits for finite, flexible weights on hom-wise cofibrant $\mathbf{Gray}$-categories.
\end{theorem}

\begin{proof}
	Define $\overline{\mathfrak{B}}$ to be the triessential image of $B\circ \mathfrak{Y}_\mathfrak{B}$ below, where $\mathfrak{Y}_{\mathfrak{B}}$ is the Yoneda embedding and $B$ is the triequivalence of Theorem \ref{trinatural transformation classifier cofibrancy} part (8). That is, $\overline{\mathfrak{B}}$ is the full sub-$\mathbf{Gray}$-category on those objects in ${[}{\mathbf{st}_{3}\left(\mathfrak{B}\right)}^\text{op}{,}\mathbf{Gray}{]}_\text{flex}$ which are biequivalent to a $\mathbf{Gray}$-functor ${\mathbf{st}_{3}\left(\mathfrak{B}\right)}^\text{op} \rightarrow \mathbf{Gray}$ of the form $L\left(\mathbf{st}_{2}\mathfrak{B}\left(E^{*}-, X\right)\right)$ for some object $X \in \mathfrak{B}$.

	$$\begin{tikzcd}[font=\fontsize{9}{6}]
		\mathfrak{B} 
		\arrow[rrdd, "G"']
		\arrow[rr, hookrightarrow,"\mathfrak{Y}_{\mathfrak{B}}"]
		&& \mathbf{TRICAT}\left(\mathfrak{B}^\text{op}{,}\mathbf{Bicat}\right)
		\arrow[rr, "B"]
		&& {[}{\mathbf{st}_{3}\left(\mathfrak{B}\right)}^\text{op}{,}\mathbf{Gray}{]}_\text{flex}
		\\
		\\
		&&
		\overline{\mathfrak{B}}
		\arrow[rruu, hookrightarrow, "H"']
	\end{tikzcd}$$
	
	\noindent $G: \mathfrak{B} \rightsquigarrow \underline{\mathfrak{B}}$ is triessentially surjective by construction, and it is tri-fully faithful by the cancellation property for tri-fully faithful trihomomorphisms since $H$ is tri-fully faithful by construction, $\mathfrak{Y}_\mathfrak{B}$ is tri-fully faithful by Lemma \ref{Yoneda lemma for tricategories}, and $B$ is the biequivalence of Theorem \ref{trinatural transformation classifier cofibrancy} part (8). $G$ is therefore a triequivalence. Let $G^{*}: \overline{B} \rightsquigarrow \mathfrak{B}$ denote the triequivalence inverse. It suffices to show that $\overline{\mathfrak{B}}$ has finite $\mathbf{Gray}$-enriched flexible limits. Suppose $\mathfrak{A}$ is a hom-wise cofibrant $\mathbf{Gray}$-category, let $W: \mathfrak{A} \rightarrow \mathbf{Gray}$ be a finite flexible weight, and let $F: \mathfrak{A} \rightarrow \overline{\mathfrak{B}}$ be a diagram. Consider the $\mathbf{Gray}$-enriched limit $\{W, HF\}$ of $HF$ with respect to $W$. By Corollary \ref{Reduction to Gray natural biequivalence} part (3), this is also a trilimit since the hom-$2$-categories of $\mathfrak{A}$ are assumed to be cofibrant. Indeed, it is a trilimit of the composite $B \circ \mathfrak{Y}_{\mathfrak{B}}\circ G^{*}\circ F$, where $B \circ \mathfrak{Y}_{\mathfrak{B}}$ is as depicted above. But by assumption the trilimit $<W, G^{*}F>$ exists in $\mathfrak{B}$ and by Proposition \ref{representables preserve trilimits} part (3), it is preserved by $\mathfrak{Y}_\mathfrak{B}$. Hence $B\circ\mathfrak{Y}_\mathfrak{B}\left(<W, G^{*}F>\right) \sim \{W, HF\}$. By the definition of $\overline{\mathfrak{B}}$ as a triessential image, it therefore has this enriched limit.
\end{proof}

\begin{remark}\label{size issues with coherence for tricategories with finite trilimits}
	We have chosen to allow $\overline{\mathfrak{B}}$ to possibly be large, since our focus is on coherence issues for tricategorical limits and this size issue is of a different nature. Nonetheless, a small choice for $\overline{\mathfrak{B}}$ can also be made similarly to what is described in the proof of Theorem 3.1 of \cite{Power Coherence for categories with finite bilimits} by closing $B\circ\mathfrak{Y}_{\mathfrak{B}}$ under finite flexible limits rather than under biequivalences. Moreover, Remark 3.2 of \cite{Power Coherence for categories with finite bilimits} applies here as well; the argument does not use the axiom of choice and can be adapted to regular cardinals $\kappa$, with the statement for finite trilimits being the case $\kappa = \aleph_{0}$.
\end{remark}

\begin{remark}\label{Remark left adjoint of flex into trilim}
	Suppose a tetracategory structure $\mathbf{TRILIM}$ was known in which the objects are tricategories with finite trilimits, and the morphisms are trihomomorphisms preserving finite trilimits. Then Theorem 3.3.4 and Corollary 3.3.5 of \cite{Campbell PhD} would be part of a tetra-homomorphism $\overline{\left(-\right)}:\mathbf{FLEX} \to \mathbf{TRILIM}$ from a four dimensional structure of $\mathbf{Gray}$-categories with finite flexibly weighted limits on hom-wise cofibrant shapes, $\Gr$-functors which preserve these limits up to isomorphism, and strict higher transfors between these. The construction of Theorem \ref{coherence for tricategories with trilimits} may be the action on objects of a tetra-adjoint to $\overline{\left(-\right)}$.	
\end{remark} 

\begin{corollary}\label{tricolimits in functor tricategories}
	Let $\mathfrak{A}$ be a tricategory and $\mathfrak{B}$ be a $\Gr$-category with tricategorical colimits. Then the $\Gr$-category $\mathbf{Tricat}(\mathfrak{A}, \mathfrak{B})$ also has tricategorical colimits.
\end{corollary}

\begin{proof}
	The proof is analogous to the one given for Proposition \ref{main examples have all trilimits and tricolimits} part (4), with the triequivalence $G: \mathfrak{B} \to \overline{\mathfrak{B}}$ of Theorem \ref{coherence for tricategories with trilimits} in place of $\mathbf{st}_{2}: \mathbf{Bicat} \to \mathbf{Gray}_{c}$.
\end{proof}

\subsection{Specific weights and diagrams}\label{subsection specific weights and diagrams}

Example \ref{center as descent object} shows that a variant of the centre construction for $\mathbf{Gray}$-monoids can be viewed as a descent object, and Example \ref{Example strictification of pseudo-double categories} views strictification of pseudo-double categories as a codescent object. Both of these are examples of flexibly weighted $\Gr$-enriched (co)limits on hom-wise cofibrant $\Gr$-categories. By Corollary \ref{Reduction to Gray natural biequivalence} part (3) these are also examples of tricategorical (co)limits. In Example \ref{strictification as a tricolimit} we see that strictification of bicategories is also a tricategorical copower, and finally Theorem \ref{biessential surjectivity on objects characterises trikleisli pseudoadjunctions in Gray or Hom} shows that trikleisli pseudoadjunctions in $\mathbf{Gray}$ and $\mathbf{Bicat}$ are characterised by having biessentially surjective on objects left pseudoadjoints.

\begin{example}\label{center as descent object}
	In this example we describe how a slight modification of the centre construction for $\mathbf{Gray}$-monoids, considered in \cite{Crans Generalised Centers}, can be seen as a tricategorical descent object. We first recall from \cite{Street monoidal centre as a limit} the analogous construction for pseudomonoids in $\left(\mathbf{Cat}, \times, \mathbf{1}\right)$, or monoidal categories. Indeed, the centre construction for monoidal categories can be seen as a bilimit. In particular, given a monoidal category $\mathcal{V}$, one first builds a truncated pseudo-cosimplicial object in $\mathbf{Cat}$ as below.
	
	$$\begin{tikzcd}
		\mathcal{V}\arrow[rr, shift left = 2, "\delta_{0}"]
		\arrow[rr, shift right = 2, "\delta_{1}"']
		&&
		{[}\mathcal{V}{,}\mathcal{V}{]}
		\arrow[rr, shift left = 5, "\delta_{0}"]
		\arrow[rr, "\delta_{1}"description]
		\arrow[rr, shift right = 5, "\delta_{2}"']
		&&
		{[}\mathcal{V}\times \mathcal{V}{,}\mathcal{V}{]}
	\end{tikzcd}$$
	
	\noindent Here $\delta_{0}: \mathcal{V} \rightarrow [\mathcal{V}, \mathcal{V}]$ corresponds via adjointness to the tensor product on $\mathcal{V}$, while $\delta_{1}:\mathcal{V} \rightarrow [\mathcal{V}, \mathcal{V}]$ corresponds to the tensor restricted along the twist map $c_{\mathcal{C}, \mathcal{D}}: \mathcal{C} \times \mathcal{D} \rightarrow \mathcal{D} \times \mathcal{C}$. The three maps from ${[}\mathcal{V}{,}\mathcal{V}{]}$ to ${[}\mathcal{V}\times \mathcal{V}{,}\mathcal{V}{]}$ correspond via adjointness to the maps depicted below.
	
	$$\begin{tikzcd}[font=\fontsize{9}{6}]
		{[}\mathcal{V}{,}\mathcal{V}{]} \times \mathcal{V}  \times \mathcal{V}
		\arrow[rd, "\overline{\delta}_{0}"']
		\arrow[r, "\text{ev}\times 1"] &	\mathcal{V} \times\mathcal{V} \arrow[d, "\bigoplus"] 
		\\
		& \mathcal{V}
		&{}
	\end{tikzcd}\begin{tikzcd}[font=\fontsize{9}{6}]
		{[}\mathcal{V}{,}\mathcal{V}{]} \times \mathcal{V}  \times \mathcal{V}
		\arrow[rd, "\overline{\delta}_{1}"']
		\arrow[r, "1\times \bigoplus"] &	{[}\mathcal{V}{,}\mathcal{V}{]} \times\mathcal{V} \arrow[d, "\text{ev}"] 
		\\
		& \mathcal{V}
		&{}
	\end{tikzcd}\begin{tikzcd}[font=\fontsize{9}{6}]
		{[}\mathcal{V}{,}\mathcal{V}{]} \times \mathcal{V}  \times \mathcal{V}
		\arrow[d, "\overline{\delta}_{2}"']
		\arrow[r, "c_{{[}\mathcal{V}{,}\mathcal{V}{]}{,}\mathcal{V}}"] 
		&	\mathcal{V} \times  {[}\mathcal{V}{,}\mathcal{V}{]} \times \mathcal{V} 
		\arrow[d, "1\times \text{ev}"] 
		\\
		\mathcal{V} 
		& \mathcal{V} \times\mathcal{V} \arrow[l ,"\bigoplus"]
	\end{tikzcd}$$
	
	\noindent There is also a functor $s_{0}: [\mathcal{V}, \mathcal{V}] \rightarrow \mathcal{V}$ given by evaluation at the unit object $I \in \mathcal{V}$, and natural isomorphisms $s_{0}d_{0} \cong 1_{V} \cong s_{0}d_{1}$ given by the left and right unitors.
	\\
	\\
	\noindent We now describe a new construction for cubical tricategories with one object, which is a greater level of generality to what is considered in \cite{Baez Neuchl Braided Monoidal 2-categories} and \cite{Crans Generalised Centers}. We will then specialise it to $\mathbf{Gray}$-monoids, recovering not $\mathbf{Tricat}_\text{su}\left(\Sigma\mathcal{V},\Sigma \mathcal{V}\right)\left(1_{\Sigma\mathcal{V}}, 1_{\Sigma\mathcal{V}}\right)$ but rather $\mathbf{Tricat}_\text{s}\left(\Sigma\mathcal{V},\Sigma \mathcal{V}\right)\left(1_{\Sigma\mathcal{V}}, 1_{\Sigma\mathcal{V}}\right)$, in which the objects are general trinatural transformations rather than unital ones. The construction we are about to describe will have the universal property of a $\mathbf{Gray}$-enriched limit, but since the weight for descent objects is flexible and the $\Gr$-category $\Delta^{G}$ described in Definition 11.1 of \cite{Gurski Coherence in Three Dimensional Category Theory} is hom-wise cofibrant, it will also have the universal property of a tricategorical limit.
	\\
	\\
	\noindent Let $\mathcal{V}$ be such a cubical tricategory with one object, considered as a monoidal bicategory. Form a truncated $\mathbf{Gray}$-enriched cosimplicial object $\underline{\mathcal{V}}:\Delta^{G} \rightarrow \mathbf{Gray}$ by similarly mapping $[n]$ to $\mathbf{Gray}\left(\otimes_{n}\mathcal{V}, \mathcal{V}\right)$, where the domain is the $n$-fold $\mathbf{Gray}$-tensor product. The $1$-cells in $\Delta^{G}$ are mapped similarly to what is described in \cite{Street monoidal centre as a limit}. The $2$-cells $A$, $L$ and $R$ are also mapped to the analogous associator and left and right unitors, which are now adjoint equivalences rather than invertible $2$-cells in $\mathbf{Gray}$. The $2$-cells $N^{ds}$, $N^{sd}$ and $N^{s}$ are mapped to identity pseudonatural transformations, viewed as $2$-cells in $\mathbf{Gray}$. Of the $3$-cells in $\Delta^{G}$,
	
	\begin{itemize}
		\item The three $3$-cells of the form $\pi_{i, j, k}$ are mapped to $\mathfrak{p}$,
		\item The two $3$-cells of the form $\mu_{i}$ are mapped to $\mathfrak{m}$,
		\item The two $3$-cells of the form $\nu_{i}^{l}$ are mapped to $\mathfrak{l}$,
		\item The two $3$-cells of the form $\nu_{i}^{r}$ are mapped to $\mathfrak{r}$,
		\item The $3$-cell $\nu_{s}$ is mapped to $\mathfrak{i}$.
	\end{itemize} 
	\noindent In the case where $\mathcal{V}$ is a $\mathbf{Gray}$-monoid, all $k$-cells in $\Delta^{G}$ for $k \geq 2$ are mapped to identities. We describe the descent object of $\underline{V}$ admits a simpler description in this setting. Its objects consist of $\left(X, \beta_{X, -}, r_{-, ?}^{X}, u^{X}\right)$ where the first three pieces of data are as described in Section 6 of \cite{Miranda semi-strictly generated closed structure on Gray-Cat}. These are given as components of the data $x$, $\varepsilon$, and $M$ of Definition 11.6 from \cite{Gurski Coherence in Three Dimensional Category Theory}. Meanwhile, the fourth datum is an invertible $2$-cell $u^{X}: \beta_{X, I} \cong 1_{\mathcal{V}}$, given as the component of the datum $U$ from Definition 11.6 of \cite{Gurski Coherence in Three Dimensional Category Theory}. This is in place of the condition written as $R_{A, I} = \text{id}_{A}$ in Section 3.1.1 of \cite{Crans Generalised Centers}, and means that the resulting $\mathbf{Gray}$-monoid will only satisfy axiom 2.9 of their Definition 2.2 up to a coherent isomorphism. It can be seen as a unitor, in reference to the trinatural transformations perspective of Section 6 of \cite{Miranda semi-strictly generated closed structure on Gray-Cat}; the axioms for codescent objects correspond to the associativity and left and right unit constraints for trinatural transformations. Similarly, the two-dimensional aspect of the universal property described in Definition 11.6 of \cite{Gurski Coherence in Three Dimensional Category Theory} may be used to describe the $1$-cells of $\mathbf{Z}\left(\mathcal{V}\right)$. In doing this, one finds that the $2$-cell and $3$-cell components written as $\alpha: g_{1}\otimes x \Rightarrow g_{2}\otimes x$ and $\Gamma$ in \cite{Gurski Coherence in Three Dimensional Category Theory} correspond to a $1$-cell and $2$-cell in $\mathcal{V}$, while the axioms specified for these data in \cite{Gurski Coherence in Three Dimensional Category Theory} correspond to the unit and composition axioms for trimodifications from $\Sigma \mathcal{V}$ to $\Sigma \mathcal{V}$.
	\\
	\\
	\noindent A similar but more complicated description can be given for the tricategorical limit in the setting of cubical tricategories with one object. Note that in this setting, $u^{X}$ will involve the left and right unitors of $\mathcal{V}$. Various $2$-cells and $3$-cells involved in the cosimplicial object corresponding to $\mathcal{V}$ are given by identities, even in this more general setting. We highlight that if a different shape $\Delta_{\text{st}}^{G}$ to $\Delta_{\text{st}}^{G}$ is used in which these $2$-cells are already identities, then $\Delta_{\text{st}}^{G}$ will fail to be hom-wise cofibrant. This means that the similarly defined weighted limit of such a diagram $\underline{\mathcal{V}}^\text{st}: \Delta_{\text{st}}^{G}$ corresponding to $\mathcal{V}$ may no longer also be a tricategorical limit.
\end{example}

\begin{example}\label{Example strictification of pseudo-double categories}
	Recall (7.5 of \cite{Limits in Double Categories}) that in a pseudo-double category, associativity and unit laws for the composition of vertical morphisms only holds up to suitably coherent double-cells which are globular in the sense that their boundary horizontal morphisms are identities. Nonetheless, every pseudo-double category is biequivalent to a strict double category i.e. one in which these laws hold on the nose. Moreover, the biequivalence is given by the unit of the left triadjoint to the inclusion $\mathbf{DblCat} \to \mathbf{PsDblCat}$ of the $\mathbf{Gray}$-category of strict double categories, strict double functors into the tricategory of pseudo-double categories and pseudo-double functors \cite{Campbell Strictification}. Given a pseudo-double category $\mathbb{A}$, this universal property suggests an alternative model for its strictification $\mathbf{st}_{2}(\mathbb{A})$ to the one described in \cite{Limits in Double Categories}. This strict double category has the same category $\mathbb{A}_{0}$ of objects and horizontal morphisms as $\mathbb{A}$, and its category $\mathbf{st}_{2}(\mathbb{A})_{1}$ of objects and vertical morphisms is freely generated from the graph of objects and vertical morphisms in $\mathbb{A}$. There are five different kinds of generating double-cells in $\mathbf{st}_{2}(\mathbb{A})$.
	
\begin{enumerate}
	\item For each double-cell in $\mathbb{A}$ there is a generating double-cell in $\mathbf{st}_{2}(\mathbb{A})$ with the same boundary.
	\item For every object $X \in \mathbb{A}$ there is a globular generating double-cell $\iota_{X}: (-)_{X} \to 1_{X}$ in $\mathbf{st}_{2}(\mathbb{A})$, corresponding to a morphism in $\mathbf{st}_{2}(\mathbb{A})_{1}$ from the empty vertical path on $X$ to the path of length one consisting of the vertical identity on $X$ in $\mathbb{A}$.
	\item For every composable pair of vertical morphisms \begin{tikzcd}
		X \arrow[r, "h"] & Y \arrow[r, "k"] & Z
	\end{tikzcd} in $\mathbb{A}$ there is a generating globular double-cell $\chi_{k, h}: (k, h) \to k\circ h$ in $\mathbf{st}_{2}(\mathbb{A})$, corresponding to a morphism in $\mathbf{st}_{2}(\mathbb{A})_{1}$ from the vertical path of length two to the vertical path of length one given by the composite in $\mathbb{A}$.
	\item There are generating double-cells given by formal inverses $\chi_{k, h}^{-1}: k\circ h \to (k, h)$ and $\iota_{X}^{-1}: 1_{X} \to (-)_{X}$ to the generating double-cells just described in (2) and (3).
\end{enumerate}

\noindent Finally, there are relations in the presentation of $\mathbf{st}_{2}(\mathbb{A})$ corresponding to each of the equations needed for the evident prospective pseudo-double functor $\eta_\mathbb{A}: \mathbb{A} \to \mathbf{st}_{2}(\mathbb{A})$ to be well-defined. \footnote{Note that when $\mathbb{A}$ is a bicategory thought of as a pseudo-double category with only identity horizontal morphisms then the construction $\mathbf{st}_{2}(\mathbb{A})$ described here indeed provides a concrete left adjoint to the inclusion $I_{2}: \mathbf{Gray} \to \mathbf{Bicat}$, as mentioned in Proposition \ref{Hom triequivalent to Gray c} part (2).}
\\
\\
	\noindent We observe that this strictification $\mathbf{st}_{2}(\mathbb{A})$ can be described as a $\mathbf{Gray}$-enriched codescent object in $\mathbf{DblCat}$ over the diagram of discrete double categories \begin{tikzcd}
		{(\Delta_{\text{ps}}^{G})}^\text{op} \arrow[r, "\underline{\mathbb{A}}"] & \mathbf{Cat} \arrow[r, "\mathbf{disc}"] & \mathbf{DblCat}
	\end{tikzcd}. Here $\Delta_\text{ps}^{G}$ is the $\Gr$-category of Definition 11.4 of \cite{Gurski Coherence in Three Dimensional Category Theory}, $\underline{\mathbb{A}}$ is the pseudo-simplicial object in $\Cat$ corresponding to the pseudo-double category $\mathbb{A}$, and $\mathbf{disc}$ views categories as strict double categories with only identity vertical morphisms. Indeed, the colimit cocone consists of \begin{itemize}
	\item a double functor $e_\mathbb{A}:\mathbf{disc}(\mathbb{A})_{0} \to \mathbf{st}_{2}(\mathbb{A})$ which includes objects and horizontal morphisms. Since $\mathbf{st}_{2}(\mathbb{A})_{0} = \mathbb{A}_{0}$, this can be identified with the component of the counit for the adjunction $\mathbf{disc} \dashv (-)_{0}$ between $\Cat$ and $\mathbf{DblCat}$ at $\mathbf{st}_{2}(\mathbb{A})$.
	\item a vertical pseudonatural transformation $\phi: e_\mathbb{A}.d_{1} \to e_\mathbb{A}.d_{0}$ whose components on objects and morphisms are the generating vertical morphisms in $\mathbf{st}_{2}(\mathbb{A})$ and generating double-cells in $\mathbf{st}_{2}(\mathbb{A})$ of type (1) described above.
	\item Invertible double modifications as displayed below, where we omit notation for $\mathbf{disc}$ and the unlabelled regions commute.
\end{itemize}


$$\begin{tikzcd}[font=\fontsize{9}{6}]
	\mathbb{A}_{2} \arrow[dd, "d_{0}"']
	\arrow[rr, "d_{2}"]\arrow[rd, "d_{1}"description]
	&{}& \mathbb{A}_{1}\arrow[rd, "d_{1}"]
	\\
	&\mathbb{A}_{1}\arrow[rr, "d_{1}"]\arrow[dd, "d_{0}"]
	&{}\arrow[dd, Rightarrow, shorten = 15, "\phi"]
	&\mathbb{A}_{0}\arrow[dd, "e_\mathbb{A}"]
	\\
	\mathbb{A}_{1}\arrow[rd, "d_{0}"']&{}&&&\Lleftarrow_{\chi}
	\\
	&\mathbb{A}_{0} \arrow[rr, "e_\mathbb{A}"']
	&{}& \mathbf{st}_{2}(\mathbb{A})
\end{tikzcd}\begin{tikzcd}[font=\fontsize{9}{6}]
	\mathbb{A}_{2} \arrow[dd, "d_{0}"']
	\arrow[rr, "d_{2}"]
	&{}
	& \mathbb{A}_{1}\arrow[rd, "d_{1}"]\arrow[dd, "d_{0}"description]
	\\
	&&&\mathbb{A}_{0}\arrow[dd, "e_\mathbb{A}"]
	\arrow[d, Rightarrow, shorten = 5, shift right = 12, "\phi"]
	\\
	\mathbb{A}_{1}\arrow[rd, "d_{0}"']\arrow[rr, "d_{1}"]
	&{}\arrow[d, Rightarrow,"\phi"]
	&\mathbb{A}_{0}\arrow[rd, "e_\mathbb{A}"]
	&{}
	\\
	&\mathbb{A}_{0} \arrow[rr, "e_\mathbb{A}"']
	&{}& \mathbf{st}_{2}(\mathbb{A})
\end{tikzcd}$$

$$\begin{tikzcd}[font=\fontsize{9}{6}]
	&&{}
	&&\mathbb{A}_{0}\arrow[rrdd, "e_\mathbb{A}"]\arrow[dddd, Rightarrow, shorten = 15, "\phi"]
	\\
	\\
	\mathbb{A}_{0} \arrow[rrrruu, bend left = 30, "1_{\mathbb{A}_{0}}"]
	\arrow[rrrrdd, bend right = 30, "1_{\mathbb{A}_{0}}"']
	\arrow[rr, "s_{0}"]
	&& \mathbb{A}_{1} \arrow[rruu, "d_{1}"]\arrow[rrdd, "d_{0}"]
	&&&&\mathbf{st}_{2}(\mathbb{A})&\Lleftarrow_{\iota}&\mathbb{A}_{0}\arrow[rr, bend left = 20, "e_\mathbb{A}"name=A]\arrow[rr, bend right = 20, "e_\mathbb{A}"'name=B] && \mathbf{st}_{2}(\mathbb{A})
	\\
	\\
	&&{}&&\mathbb{A}_{0}\arrow[rruu, "e_\mathbb{A}"]
	\arrow[from=A, to=B, Rightarrow, shorten = 5, "1_{e_\mathbb{A}}"]
\end{tikzcd}$$

\noindent Finally, the equations that these data must satisfy in order to be a pseudo-codescent object are indeed precisely those required for $\mathbb{A} \to \mathbf{st}_{2}(\mathbb{A})$ to assemble into a pseudo-double functor. As in the previous example, since the weight is flexible and the shape has cofibrant hom-$2$-categories, this colimit also has a tricategorical universal property.
\end{example}

\begin{remark}\label{Remark lax double functors}
	In \cite{Limits in Double Categories} it is observed that the assignment of a pseudo-double category $\mathbb{A}$ to the strict double category $\mathbb{P}$ there described does not extend to lax double functors. This is because the generating double cells in $\mathbb{P}$ have boundary vertical paths of arbitrary length, and hence require both lax and oplax functoriality constraints. However, the category $\mathbf{DblCat}$ is also reflective inside the category $\mathbf{PsDblCat}_\text{lax}$ of pseudo-double categories and lax double functors, with the reflection $\mathbb{A} \mapsto \mathbf{Lax}(\mathbb{A})$ being given by a variant of what is described in Example \ref{Example strictification of pseudo-double categories} but with $\Delta^{G}$ in place of $\Delta_\text{ps}^{G}$. In particular, given a pseudo-double category $\mathbb{A}$, there is a universal lax double functor $\mathbb{A} \to \mathbf{Lax}(\mathbb{A})$ into a strict double category, and the assignment $\mathbb{A} \mapsto \mathbf{Lax}(\mathbb{A})$ is indeed functorial, sending lax double functors to strict ones.
\end{remark}

\begin{example}\label{strictification as a tricolimit}
	We describe how strictification of a bicategory $\mathcal{W}$ can be seen as a fairly simple tricolimit. This tricolimit is reasonable to call a `tricopower of the point by $\mathcal{W}$'. Let $W: \mathbf{1}^\text{op} \rightarrow \mathbf{Bicat}$ denote the trihomomorphism that is constant on a bicategory $\mathcal{W} \in \mathbf{Bicat}$, and let $T: \mathbf{1} \rightarrow \mathbf{Gray}$ denote the $\mathbf{Gray}$-functor that is constant on the terminal object. Then a tricolimit $\mathcal{W} \odot T$ is given by its strictification, as described for example in Section 2.2.3 of \cite{Gurski Coherence in Three Dimensional Category Theory}.
	
	\begin{align*}
		\mathbf{Gray}\left(\mathbf{st}_{2}(\mathcal{W}), \mathcal{B}\right)
		&\cong \mathbf{Bicat}\left(\mathcal{W}, I\mathcal{B}\right)
		&\text{Proposition} \ref{Hom triequivalent to Gray c}
		\\
		& \sim \widehat{{\mathbf{1}}^\text{op}}\left(W,I\mathcal{B}\right)
		&\text{Proposition \ref{trilimits in Hom}}
		\\
		&=  \widehat{{\mathbf{1}}^\text{op}}\left(W,\mathbf{Gray}\left(T, I\mathcal{B}\right)\right)
		& T\text{ terminal} 
	\end{align*}

\end{example}

\noindent We now define the tricategorical version of Kleisli objects for pseudomonads.

\begin{definition}\label{trikleisli object definition}
	Let $\left(A, S\right): \mathbf{Psmnd} \rightarrow \mathfrak{K}$ be a strict trihomomorphism into a tricategory $\mathfrak{K}$. A \emph{trikleisli object} of $\left(A, S\right)$, if it exists, will be a tricolimit of $\left(A, S\right)$ weighted by the trihomomorphism depicted below.
	
	$$\begin{tikzcd}
		\mathbf{Psmnd}^\text{op} \arrow[rr, "W"] && \mathbf{Gray} \arrow[rr, "I"] && \mathbf{Bicat}
	\end{tikzcd}$$

 \noindent Here $W$ is the weight for $\mathbf{Gray}$-categorical Kleisli objects, described in Remark 3.2.5 of \cite{Miranda Enriched Kleisli objects for pseudomonads}, and $I_{2}: \Gr \to \mathbf{Bicat}$ is the inclusion.
\end{definition}

\begin{theorem}\label{biessential surjectivity on objects characterises trikleisli pseudoadjunctions in Gray or Hom}
	Let \begin{tikzcd}
		\mathcal{B}
		\arrow[rr, shift right = 2,"U"'] 
		&\bot
		& \mathcal{A}
		\arrow[ll, shift right = 2, "F"']
	\end{tikzcd} be a pseudoadjunction in a tricategory $\mathfrak{K}$ where $\mathfrak{K}$ is either $\mathbf{Gray}$ or $\mathbf{Bicat}$. Then $\mathcal{B}$ is trikleisli for the induced pseudomonad on $\mathcal{A}$ if and only if $F$ is biessentially surjective on objects.
\end{theorem}

\begin{proof}
	For $\mathfrak{K} = \mathbf{Gray}$, first recall that $\mathbf{Psmnd}$ is hom-wise cofibrant and that the weight $W:\mathbf{Psmnd}^\text{op} \rightarrow \mathbf{Gray}$ is projective cofibrant, by Corollary 5.2.5 of \cite{Miranda Enriched Kleisli objects for pseudomonads}. Then the $\mathbf{Gray}$-enriched colimit is trikleisli by Corollary \ref{Reduction to Gray natural biequivalence} part (3). But, by $(1) \iff (2)$ of Theorem 4.2.2 of \cite{Miranda Enriched Kleisli objects for pseudomonads}, the canonical comparison from the $\mathbf{Gray}$-enriched Kleisli object is a biequivalence if and only if $F$ is biessentially surjective on objects. The result then follows from Corollary \ref{trilimits stable under biequivalence} part (2).
	\\
	\\
	\noindent For part (2), first apply $\mathbf{st}_{2}: \mathbf{Bicat} \rightsquigarrow \mathbf{Gray}$. Observe that $\mathbf{st}_{2}$ is tricocontinuous by Lemma \ref{left triadjoint preserves tricolimits}, and it preserves and reflects the property that $1$-cells may have of being biessentially surjective on objects. Then apply part (1) to see that the pseudoadjunction is trikleisli in $\mathbf{Gray}$ if and only if $\mathbf{st}_{2}(F)$ is biessentially surjective on objects. Finally, $I_{2}: \mathbf{Gray} \rightarrow \mathbf{Bicat}$ also preserves and reflects biessential surjectivity on objects, and by Proposition \ref{tricolimits preserves by I: Gray --> Hom} it also preserves the trikleisli object of the pseudomonad $\left(\mathbf{st}_{2}(\mathcal{A}), \mathbf{st}_{2}\left(UF\right)\right)$. To see that the assumptions of Proposition \ref{tricolimits preserves by I: Gray --> Hom} apply, first observe that $\mathbf{st}_{2}(\mathcal{A})$ is cofibrant by Proposition \ref{Hom triequivalent to Gray c} part (3), hence the $\mathbf{Gray}$-enriched Kleisli object is also cofibrant by Proposition 5.1.4 of \cite{Miranda Enriched Kleisli objects for pseudomonads}, and finally the $\mathbf{Gray}$-enriched Kleisli object is also a trikleisli object as discussed in the first paragraph of this proof.
\end{proof}

\begin{remark}
	Tricategorical colimits $W \odot F$ in $\Gr$ may be given as flexibly weighted colimits $W' \cdot F$ enriched over $\Gr$. Since $\Gr$ is locally presentable, these enriched colimits may in turn be described via presentations involving generators and relations. It is typically harder to infer from such a presentation a more explicit description of the $2$-category structure of $W' \cdot F$. When $W$ is the weight for trikleisli objects then we may take $W'$ to be $W$ itself. A description in terms of generators and relations in given in Section 3.3 of \cite{Miranda Enriched Kleisli objects for pseudomonads}. In contrast, Theorem \ref{biessential surjectivity on objects characterises trikleisli pseudoadjunctions in Gray or Hom} shows that the tricategorical colimit $W \odot F$ may indeed be described explicitly, rather than just via a presentation. In particular, the $2$-category of free pseudoalgebras for the pseudomonad $F: \mathbf{Psmnd} \to \Gr$ is such a description. Moreover, if the tricolimit is being considered in $\mathbf{Bicat}$ rather than in $\Gr$ then even the Kleisli bicategory, as described in Section 2.2 of \cite{Miranda Enriched Kleisli objects for pseudomonads}, provides such as explicit description.
\end{remark}

\begin{example}\label{tricategorical eilenberg-moore for opmonoidal pseudomonads}
	It is well-known that there is a free monoid monad on the monoidal category $(\Gr, \otimes, \mathbf{1})$, but that this monad does enrich over $\Gr$ as it is given by pseudofunctors between hom-$2$-categories. Nonetheless, strict and pseudo algebras for the resulting three-dimensional monad can be considered, as per Chapter 13 of \cite{Gurski Coherence in Three Dimensional Category Theory}. These correspond to $\mathbf{Gray}$-monoids and one-object cubical tricategories respectively. One can consider these as objects of tricategories determined by oplax morphisms and oplax transformations of algebras, dually to what is described in Section 13.3 of \cite{Gurski Coherence in Three Dimensional Category Theory}. The strict algebras form the $\Gr$-category $\mathbf{Gray}$-$\text{monoids}_\text{oplax}$ described in Section 3.1 of \cite{Miranda Opmonoidal Pseudomonads}, while the pseudoalgebras form the sub-tricategory $\mathbf{cubMonBicat}_\text{oplax} \hookrightarrow \mathbf{MonBicat}_\text{oplax}$ spanned by those monoidal bicategories that are cubical as tricategories with one object, and opmonoidal $2$-functors between them. The results of Section 5 of \cite{Miranda Opmonoidal Pseudomonads} may then be interpreted as saying that
	
	\begin{itemize}
		\item $\mathbf{Gray}$-$\mathbf{monoids}_\text{oplax}$ has a $\mathbf{Gray}$-enriched Eilenberg-Moore object for a pseudomonad \cite{Coherent Approach to Pseudomonads} if and only if the opmonoidal $2$-functor $S: \A \to \A$ has identity modifications $\gamma$, $\omega$, $\delta$ mediating the usual associativity and unit laws for opmonoidal functors.
		\item Nonetheless, both $\mathbf{Gray}$-$\mathbf{monoids}_\text{oplax}$ and $\mathbf{cubMonBicat}_\text{oplax}$ have tricategorical limits weighted by the usual weight for Eilenberg-Moore objects for pseudomonads. For $\mathbf{Gray}$-$\mathbf{monoids}_\text{oplax}$ this is given by the strictification of the monoidal bicategory $(\A^S, \overline{\otimes}, \overline{I}, \gamma, \omega, \delta)$ of Theorem 5.2.5 of \cite{Miranda Opmonoidal Pseudomonads}.
		\item A tricategorical Eilenberg-Moore object for a pseudomonad on a strict algebra can be chosen in $\mathbf{cubMonBicat}_\text{oplax}$ such that the trihomomorphism $\mathbf{cubMonBicat}_\text{oplax} \to \mathbf{Gray}$ which forgets the monoidal structure sends tricategorical Eilenberg-Moore objects to $\Gr$-enriched Eilenberg-Moore objects, and hence preserves tricategorical Eilenberg-Moore objects.
	\end{itemize}  

\noindent Dual considerations also apply to the results of \cite{Miranda Kleisli Bicategories for Symmetric Monoidal Pseudomonads} for monoidal pseudomonads. In particular, extended monoidal structures on Kleisli bicategories can be seen as trikleisli objects for pseudomonads in appropriate tricategories of monoidal bicategories.
\end{example}

\section{Conclusions and future directions}

\noindent We have described how the data for a tricategorical (co)limit can be transformed into the data for a flexibly, or projective cofibrantly, weighted $\Gr$-enriched (co)limit on a shape whose hom-$2$-categories are cofibrant. Moreover, if the $\Gr$-enriched (co)limit of such data exists then it is also a tricategorical (co)limit of the same data, while conversely a tricategorical (co)limit of such data is analogous to the $\Gr$-enriched notion but with biequivalences in place of isomorphisms of $2$-categories. We have used these results to give a three-dimensional analogue of Power's coherence theorem for bicategorical limits. We have also shown that leading examples of tricategories such as $\mathbf{Bicat}$, $\Gr$ and certain functor tricategories have tricategorical limits and colimits. Finally, we have shown that these tricategorical universal properties capture certain constructions that fail to have $\Gr$-enriched universal properties in their own right, such as Kleisli bicategories and lifted monoidal structures for opmonoidal pseudomonads.
\\
\\
\noindent Many other universal constructions with a tricategorical flavor have been described in the literature. See \cite{Walker universal property of bicategory of polynomials} for examples to do with bicategories of polynomials, $\mathbf{Poly}\left(\mathcal{E}\right)$, and sub-bicategories of spans $\mathbf{Span}\left(\mathcal{E}\right)$. It is plausible that these are also expressible as tricategorical colimits. Finally, as mentioned in Remark \ref{Remark enriched weakness}, we have effectively reduced the study of tricategorical universal properties to the setting of enriched weakness with respect to the class of weak equivalences $\mathcal{E} \subseteq \mathcal{V}:= \Gr$. The further development of weak enriched universal properties in such settings is hence motivated by examples from three-dimensional category theory. We leave such investigations and the further development of this theory to future research.
\\
\\
\noindent In a follow-up paper \cite{Miranda Pseudo-double categories as pseudomonads}\footnote{for now see Chapter 7 of \cite{Miranda PhD}}, we will develop some of the theory of free (co)completions under classes of weights for tri(co)limits $\Phi$, and consider in detail the case where $\Phi$ consists of just the weight for trikleisli objects. In this way we will develop an analogue of the theory of wreaths for pseudomonads, and apply this to the theory of pseudo-double categories.

\end{document}